\documentclass[UTF-8,reqno]{amsart}
\usepackage{enumerate}
\usepackage{mhequ}
\usepackage{amssymb,url,color, booktabs,nccmath}
\usepackage[left=3.2cm,right=3.2cm,top=4cm,bottom=4cm]{geometry}
\usepackage{mathrsfs}
\usepackage{enumitem,dsfont}

\usepackage{color}
\usepackage[colorlinks=true]{hyperref}
\hypersetup{
    linkcolor=blue,          
    citecolor=red,        
    filecolor=blue,      
    urlcolor=cyan
}

\definecolor{darkergreen}{rgb}{0.0, 0.5, 0.0}

\numberwithin{equation}{section}
\def\theequation{\arabic{section}.\arabic{equation}}
\newcommand{\be}{\begin{eqnarray}}
\newcommand{\ee}{\end{eqnarray}}
\newcommand{\ce}{\begin{eqnarray*}}
\newcommand{\de}{\end{eqnarray*}}
\newtheorem{theorem}{Theorem}[section]
\newtheorem{lemma}[theorem]{Lemma}
\newtheorem{proposition}[theorem]{Proposition}
\newtheorem{Examples}[theorem]{Example}
\newtheorem{corollary}[theorem]{Corollary}

\newtheorem{definition}[theorem]{Definition}
\theoremstyle{definition}
\newtheorem{remark}[theorem]{Remark}

\DeclareMathOperator{\supp}{supp}

\def\[{{\Big[}}
\def\]{{\Big]}}
\def\<{{\langle}}
\def\>{{\rangle}}
\def\({{\Big(}}
\def\){{\Big)}}

\def\bx{{\mathbf{x}}}
\def\tr{\mathrm {tr}}

\def\dif{{\mathord{{\rm d}}}}

\def\={&\!\!=\!\!&}

\def\cH{{\mathcal H}}

\def\cL{{\mathcal L}}

\def\1{{\mathbf{1}}}

\def\geq{\geqslant}
\def\leq{\leqslant}

\def\div{\mathord{{\rm div}}}

\def\[{{\Big[}}
\def\]{{\Big]}}
\def\<{{\langle}}
\def\>{{\rangle}}
\def\({{\Big(}}
\def\){{\Big)}}

\def\bx{{\mathbf{x}}}
\def\tr{\mathrm {tr}}

\def\dif{{\mathord{{\rm d}}}}

\def\={&\!\!=\!\!&}
\def\bt{\begin{theorem}}
\def\et{\end{theorem}}
\def\bl{\begin{lemma}}
\def\el{\end{lemma}}
\def\br{\begin{remark}}
\def\er{\end{remark}}
\def\bx{\begin{Examples}}
\def\ex{\end{Examples}}
\def\bd{\begin{definition}}
\def\ed{\end{definition}}
\def\bp{\begin{proposition}}
\def\ep{\end{proposition}}
\def\bc{\begin{corollary}}
\def\ec{\end{corollary}}

\def\fa{{\mathfrak a}}
\def\fc{{\mathfrak c}}
\def\fm{{\mathfrak m}}
\def\fC{{\mathfrak C}}

\def\geq{\geqslant}
\def\leq{\leqslant}

\def\div{\mathord{{\rm div}}}

\def\Id{\textrm{Id}}

 \def\R{\mathbb R}
 \def\R{\mathbb R}    
\def\N{\mathbb N} \def\kk{\kappa} 
   
\def\<{\langle} \def\>{\rangle}

\def\${|\!|\!|}

\allowdisplaybreaks

\begin{document}

\title[Nonuniqueness for 2D NSE Driven by a Space-Time White Noise]{ Non-unique Ergodicity for the 2D Stochastic  Navier-Stokes Equations with  Derivative of Space-Time White Noise}

\author{Huaxiang L\"u}
\address[H. L\"u]{Academy of Mathematics and Systems Science,
Chinese Academy of Sciences, Beijing 100190, China}
\email{lvhuaxiang22@mails.ucas.ac.cn }

\author{Xiangchan Zhu}
\address[X. Zhu]{ Academy of Mathematics and Systems Science,
Chinese Academy of Sciences, Beijing 100190, China}
\email{zhuxiangchan@126.com}

\thanks{
Research  supported   by National Key R\&D Program of China (No. 2022YFA1006300, 2020YFA0712700) and the NSFC (No.  12090014, 12288201) and
  the support by key Lab of Random Complex Structures and Data Science, Chinese Academy of Science. The financial support by the DFG through the CRC 1283 "Taming uncertainty and profiting
 from randomness and low regularity in analysis, stochastics and their applications" is greatly acknowledged.
}

\begin{abstract}
We prove existence of infinitely many stationary solutions as well as ergodic stationary solutions for  the stochastic Navier-Stokes equations  on $\mathbb{T}^2$
\begin{align*}
\dif u+\div(u\otimes u)\dif 
t+\nabla p\dif t&=\Delta u\dif t
+ (-\Delta)^{\fa/2}\dif B_t,\ \ \ \ \div u=0,\notag
\end{align*}
driven by   derivative of space-time white noise,  where $\fa\in[0,\frac13)$. In this setting, the solutions are not function valued and probabilistic renormalization is required to give a meaning to the equations. Finally, we show that the stationary distributions are  not Gaussian distribution $N(0,\frac12(-\Delta)^{\fa-1})$.  The proof relies on a  time-dependent decomposition and a stochastic version of the convex integration method which provides uniform moment bounds in some function spaces.
\end{abstract}


\date{\today}

\maketitle

\tableofcontents
\section{Introduction}

In the present paper, we consider the following 
 two dimensional stochastic Navier–Stokes systems on $\mathbb{T}^2=\mathbb{R}^2/\mathbb{Z}^2$  
\begin{align}\label{1.1}
\dif u+\div(u\otimes u)\dif 
t+\nabla p\dif t&=\Delta u\dif t+(-\Delta)^{\fa/2}\dif B_t,\\
\div u&=0,\notag
\end{align}
where $u$ is the fluid
velocity, $p$ is the associated pressure, $\fa\in[0,\frac13)$. Here $B$ is a cylindrical Wiener process on some stochastic basis $(\Omega, \mathcal{F}, (\mathcal{F}_t)_{t\in\mathbb{R}}, \mathbf{P})$. 
The time
derivative of $B$  is the delta correlated  space-time white noise. Due to irregularity of noise, the solution $u$ is not function valued and the non-linear term is not well-defined in the classical sense.

When $\fa=0$, Da Prato and Debussche \cite{DPD02} viewed \eqref{1.1} as perturbations of linear equation and established the local well-posedness of \eqref{1.1}   as well as global-in-time existence for almost every initial condition with respect to the Gaussian invariant measure $N(0,\frac12(-\Delta)^{-1})$. This result can be strengthened to all initial data by the strong Feller property (see \cite{ZZ17}). More recently, Hairer and Rosati \cite{HR23} employed a PDE argument to prove the existence and uniqueness of a global-in-time solution, while sharp non-uniqueness was demonstrated through the convex integration method in \cite{LZ23}. Furthermore, utilizing the support theorem (see \cite{HS21}), the uniqueness of the invariant measure with support in $C^{-\kappa}$ for $\kappa>0$ small enough was established.  However,  when $\fa>0$, existence of invariant measure is unknown and  it seems difficult to use the  method in \cite{HR23}  in this  case.

The more irregular 3D case necessitates innovative approaches, i.e. the theory of regularity structure \cite{Hai14} or paracontrolled calculus \cite{GIP15}. These frameworks have proven effective in handling various singular subcritical SPDEs, as evidenced by studies on the Kardar-Parisi-Zhang (KPZ) equation and stochastic quantization equations for quantum fields (cf. \cite{BHZ19, CH16, BCCH21, OSSW21, HS23}). In particular, they led to a local well-posedness theory for \eqref{1.1} in the 3D case (see \cite{ZZ15}). However, global existence of solutions and existence of invariant measure is more challenging and it seems out of reach by using classical PDE argument. Recently, Hofmanov\'{a}, Zhu, and the second named author \cite{HZZ21b} constrcuted non-unique global solutions  in the  3D case  via combination of paracontrolled calculus and convex integration. The existence of invariant measures in the 3D case still remains an open problem. The main motivation of this work is to deduce existence of global stationary solutions in the 3D case following the convex integration method, which is more complicated  as it requires  regularity structure  or paracontrolled calculus to give a meaning to the solutions.

In this paper, our focus is on situations where  Da Prato-Debussche's trick is sufficient to confer meaningful interpretations to solutions: the case where $\fa\in[0,\frac13)$; additional clarification can be found in Remark \ref{rmk:fa}. Subsequently, we employ the stochastic convex integration method to construct infinitely many solutions possessing uniform in time moment bounds. This construction, in turn, yields to infinitely many  (ergodic) stationary solutions by Krylov-Bogoliubov's argument. The more complicated 3D case will be considered in our future work.

  As mentioned above, the main tool in this paper is the convex integration method.  Since the seminal works by De Lellis and Sz\'ekelyhidi \cite{DLS09, DLS10, DLS13}, convex integration has become a powerful tool in the context of fluid dynamics. It has yielded numerous remarkable results pertaining to the Navier–Stokes equations and Euler equations in the deterministic setting  \cite{BDLIS15, BCV18,BDLSV19, BV19b,CL20, CL22,GR23,GKN23}.
Notably, the convex integration method has found successful application in the stochastic setting as well, see \cite{BFH20,CFF19, CDLDR18,HZZ19,HZZ22,HZZ,HZZ21b, RS21,  Yam21b, Yam22a, Yam22b}.

\subsection{Main result}
Using Da Prato-Debussche's trick, we divide the equation (\ref{1.1}) into two parts  $u=v+z$:

\begin{align}
    \dif z-\Delta z \dif t+  z \dif t+\nabla p_z\dif t&=(-\Delta)^{\fa/2}\dif B_t,\notag\\
    \div z&=0,\label{3.1}
\end{align}
 and
  $v$ solves the nonlinear equation
\begin{align}\label{3.2}
    \dif v-\Delta v\dif t-  z\dif t+\div((v+z)\otimes(v+z))\dif t+\nabla p_v\dif t&=0,\\
    \div v&=0.\notag
\end{align}
 The stationary 
solution to (\ref{3.1}) can be written as
\begin{align}
z(t)=\int_{-\infty}^te^{(t-s)(\Delta-I)}(-\Delta)^{\fa/2}\dif B_s,\notag
\end{align}
where $e^{t\Delta}$ is the heat semigroup.  Note that $z\otimes z$ is not well-defined in the classical sense. As $z$ is Gaussian, we can understand $z\otimes z$ as a Wick product using renormalisation and denote it as $z^{:2:}$. We give more details in Section \ref{sec:pe}.

\br\label{rmk:fa}
 By Proposition \ref{prop:Z},  $z$ is a random distribution of space-time regularity $-\fa-\kk$ for any $\kk>0$ and $z^{:2:}$ is of regularity $-2\fa-\kk$. Applying Schauder's estimate,  $v$ is a distribution of regularity no more than $1-2\fa-\kk$. Therefore, the term $v\otimes z+z\otimes v$ in (\ref{3.2}) is well-defined only when $0\leq\fa<\frac13$. When $\fa\geq\frac13$, a further decomposition is necessary and we refer to  \cite{HZZ21b} for further details.
\er

Within this study, we focus on analytically weak solutions which satisfy the equations in the
following sense.
\bd\label{def:weaksol}
 We say that $((\Omega, \mathcal{F},(\mathcal{F}_t)_{t\in\mathbb{R}}, \mathbf{P}), u, B)$ is an analytically weak solution to the Navier-Stokes system (\ref{1.1}) provided
 
(1) $(\Omega, \mathcal{F},(\mathcal{F}_t)_{t\in\mathbb{R}}, \mathbf{P})$ is a stochastic basis with a complete right-continuous filtration;

(2) $B$ is an $\mathbb{R}^2$-valued, spatial mean and divergence free, cylindrical Wiener process with respect to the filtration $(\mathcal{F}_t)_{t\in\mathbb{R}}$;

(3)  $ u-z \in L^2_{\rm loc}(\mathbb{R}; L^2)\cap C(\mathbb{R}; W^{\frac13,\frac65})$ $\mathbf{P}$-a.s. and is $(\mathcal{F}_t)_{ t\in\mathbb{R}}$-adapted, where  $z(t)=\int_{-\infty}^te^{(t-s)(\Delta-I)}$ $(-\Delta)^{\fa/2}\dif B_s\in C(\mathbb{R};C^{-\fa-\kappa})$  $\mathbf{P}$-a.s.  $z^{:2:}$ is defined in the sense of Wick product belongs to $C(\mathbb{R};C^{-2\fa-\kappa}) $  $\mathbf{P}$-a.s. for $\kappa>0$ small enough;

(4)  for every $ -\infty < s \leq t < \infty$ it holds $\mathbf{P}$-a.s.
\begin{align}
    \langle (u-z)(t), \psi\rangle& + \int_s^t \langle
\div\((u-z)\otimes (u-z)+(u-z)\otimes z+z\otimes (u-z)+z^{:2:}\), \psi\rangle \dif r\\
&= \langle (u-z)(s), \psi\rangle+  \int_s^t \langle 
z, \psi\rangle \dif r + \int_s^t  \langle \Delta (u-z),\psi\rangle \dif r\label{3.2*}
\end{align} 
 for all $\psi \in C^\infty(\mathbb{T}^2), \div\psi = 0.$ 
\ed

As pathwise non-uniqueness and non-uniqueness in law 
have been established in our previous work \cite{LZ23}, we understand stationarity
in the sense of shift invariance of laws of solutions on the space of trajectories, see also \cite{BFHM19,
BFH20, FFH21, HZZ22, HZZ22c,HLZZ23}.
More precisely, we define the joint trajectory space for the regular part of the solution and the driving Brownian motion as
\begin{align*}
 \mathcal{T}:=\(C(\mathbb{R}; W^{\frac13,\frac65})\cap L^2_{\rm loc}(\mathbb{R};L^2) \) 
  \times C(\mathbb{R};C^{- 1-\kappa})
\end{align*}
for some $\kappa>0$ small enough.
Let $S_t, t \in\mathbb{R}$ be shifts on trajectories given by
\begin{align*}
S_t(v,B)(\cdot) = (v(\cdot + t), B(\cdot + t) - B(t)),\ \ t\in\mathbb{R},\ \ (v,B)\in\mathcal{T}.
\end{align*}
We note that the shift in the last component acts differently in order to guarantee that for a
Brownian motion $B$ the shift $S_tB$ is again a Brownian motion.

Stationary solutions to the stochastic Navier–Stokes equations (\ref{1.1}) are defined as follows.

\bd \label{def:statsol}
We say that $((\Omega, \mathcal{F},(\mathcal{F}_t)_{t \in\mathbb{R}}, \mathbf{P}), u, B)$ is a stationary solution to the stochastic
Navier–Stokes equations provided it satisfies (\ref{1.1}) in the sense of Definition \ref{def:weaksol} and its law is shift invariant, that is,
\begin{align*}
\mathcal{L}[S_t(u-z,B)] = \mathcal{L}[u-z,B]
\end{align*}
for all $t \in\mathbb{R}$, where  $z(t)=\int_{-\infty}^te^{(t-s)(\Delta-I)}(-\Delta)^{\fa/2}\dif B_s$.
\ed

We proceed by defining the ergodic stationary solution, which aligns with the concept of ergodicity for the dynamical system $(\mathcal{T} , \mathcal{B}(\mathcal{T} ), (S_t, t \geq0), \mathcal{L}[u-z, B])$, where $\mathcal{B}(\mathcal{T}  )$ denotes the $\sigma$-algebra of Borel sets on $\mathcal{T} $ . Accordingly,
we may formulate ergodicity of stationary solutions as ergodicity of the associated dynamical system.

\bd\label{def:ergstatsol}
A stationary solution $((\Omega, \mathcal{F},(\mathcal{F}_t)_{ t\in\mathbb{R}}, \mathbf{P}), u, B)$ is ergodic provided
$$\mathcal{L}[u-z, B](A) = 1\ or\ \mathcal{L}[u-z, B](A) = 0$$
for all $A \subset \mathcal{T}$ Borel and shift invariant, where  $z(t)=\int_{-\infty}^te^{(t-s)(\Delta-I)}(-\Delta)^{\fa/2}\dif B_s$.
\ed

Our main result reads as follows:
\bt\label{thm:exist}
There exist

$(1)$ infinitely many stationary solutions;

$(2)$ infinitely many ergodic stationary solutions;\\
to the stochastic Navier–Stokes equation (\ref{1.1}). \\
Moreover, these stationary distributions are  not  Gaussian distribution  $N(0,\frac12(-\Delta)^{ \fa-1})$.
\et

In this paper, we follow the  stochastic convex integration in \cite{HZZ22,CDZ23,HZZ22c,HLZZ23} to construct solutions with uniform in time bound by introducing expectations during the iterative estimates in convex integration. As in \cite{LZ23,HZZ21b}, we use a decomposition of the system \eqref{3.2}, which makes this singular setting amenable to convex integration. More precisely, we split $v=v^1+v^2$ where $v^1$ represents the irregular part and $v^2$ the regular part. Then we put forward a joint iteration procedure approximating both equations.  Our primary objective is to establish a uniform moment bound for both $v^1$ and $v^2$. However, due to the nonlinearity of the equation, estimating any $p$th-moment involves higher moments.  For $v_q^2$ this could be achieved by choosing parameters during convex integration carefully  as in \cite{HZZ22}. But for $v_q^1$, as we use Schauder's estimate, this increasing moment makes  impossible to close the estimates for $v_q^1$. In this paper, we address this issue by choosing a time-dependent localizer $R$ such that the  high frequency parts of $z$ can be bounded by a constant. Subsequently, we have to bound the time regularity of these terms in the estimate of $v_q^2$ ( see Proposition \ref{deltaz} and Section \ref{votiefr}). Furthermore, as we aim at a uniform in time bound for $v^1_q$,  we employ an induction on time intervals and introduce an extra damping term in the equation of $v^1$ to achieve the required estimate (see Section \ref{4.13-14}).

We recall that in the previous paper \cite{HZZ22, HZZ}, the authors demonstrated non-uniqueness through controlling the kinetic energy of the constructed solutions. While this argument could be applied to infer non-uniqueness for $v$, it doesn't imply non-uniqueness in the law of $u$ since the joint law of $(v,z)$ remains unknown.  In the present paper, non-uniqueness follows from different starting iterations.

\noindent{\bf Organization of the paper.}
In Section~\ref{s:not} we collect the basic notations and preliminaries
used throughout the paper.  In Section \ref{cog}, we  present a formal decomposition
of the system into the system for $v^1$ and $v^2$, and give the estimate for stochastic objects.   Section \ref{sec:adc} sets up the iterative convex integration procedure. Estimates of $v^1_q$  are presented in Section \ref{est:vq1}.  Section \ref{proofof3.1} is devoted to the core of the convex integration construction, namely, the iteration Proposition \ref{prop:3.1}.  In Section \ref{sec6} we used a Krylov–Bogoliubov’s argument  to provide proofs of existence of non-unique stationary solutions, namely, Theorem \ref{thm:exist}. Finally in Appendix we collect several auxiliary results.

\section{Preliminaries}
\label{s:not}

  Throughout the paper, we employ the notation $a\lesssim b$ if there exists a constant $c>0$ such that $a\leq cb$.  We write $a \simeq b$ if $a\lesssim b$ and $b\lesssim a$.
 \subsection{Function spaces} 
  Given a Banach space $E$ with a norm $\|\cdot\|_E$ and a domain $ D\subset\mathbb{R}$, we write $C_DE=C(D;E)$ for the space of continuous functions from $D$ to $E$, equipped with the supremum norm $\|f\|_{C_DE}=\sup_{t\in D}\|f(t)\|_{E}$. 
  For $\alpha\in(0,1)$ we  define $C^\alpha_DE$ as the space of $\alpha$-H\"{o}lder continuous functions from $D$ to $E$, endowed with the norm $\|f\|_{C^\alpha_DE}=\sup_{s,t\in D,s\neq t}\frac{\|f(s)-f(t)\|_E}{|t-s|^\alpha}+\sup_{t\in D}\|f(t)\|_{E}.$ Here we use $C_D^\alpha$ to denote the case when $E=\mathbb{R}$. 
  For $p\in [1,\infty]$ we write $L^p_DE=L^p(D;E)$ for the space of $L^p$-integrable functions from $D$ to $E$, equipped with the usual $L^p$-norm. 
    We use $L^p$ to denote the set of  standard $L^p$-integrable functions from $\mathbb{T}^2$ to $\mathbb{R}^2$. For $s>0$, $p>1$ we set $W^{s,p}:=\{f\in L^p; \|(I-\Delta)^{\frac{s}{2}}f\|_{L^p}<\infty\}$ with the norm  $\|f\|_{W^{s,p}}=\|(I-\Delta)^{\frac{s}{2}}f\|_{L^p}$. 
 We define $\mathcal{S}(\mathbb{R}^2)$ as the Schwartz space on $\mathbb{R}^2$, and the distribution space $\mathcal{S}'(\mathbb{R}^2)$ is defined by duality. Similarly we define the space of $\R^2$-valued distributions as $ \mathcal{S}'(\mathbb{R}^2; \mathbb{R}^2)$.
    
 We denote by  $C^{N}_{D,x}$ the space of $C^{N}$-functions on $D\times\mathbb{T}^{2}$ with $N\in\N_0:=\N\cup\{0\}$. The spaces are equipped with the norms
$$
\|f\|_{C^N_{D,x}}=\sum_{\substack{0\leq n+|\alpha|\leq N\\ n\in\N_{0},\alpha\in\N^{2}_{0} }}\sup_{t\in D}\|\partial_t^n \partial_x^\alpha f\|_{ L^\infty}.
$$

We define the projection onto null-mean functions
 as $\mathbb{P}_{\neq0} f := f -\int_{\mathbb{T}^d}\!\!\!\!\!\!\!\!\!\!\!\!\!\; {}-{}\ \ f\dif x$. For a tensor $T$ , we denote its traceless part by $\mathring{T} := T- \frac{1}{2} \tr (T )\Id$. By 
$\mathcal{S}^{2\times 2}$ we denote the space of symmetric matrix and by $\mathcal{S}^{2\times 2}_0$ the space of symmetric trace-free matrix.

We also introduce the  the Fourier and the inverse Fourier transform on $\mathbb{R}^2$. For $f\in \mathcal{S}(\mathbb{R}^2)$
$$\mathcal{F} f (\zeta) :=\int_{\mathbb{R}^2}e^{-2\pi ix\cdot \zeta} f (x)\dif x, \ \ \ \mathcal{F}^{-1}f (x) := \int_{\mathbb{R}^2}e^{2\pi ix\cdot \zeta} f  (\zeta)\dif \zeta.$$
For $f\in \mathcal{S}'(\mathbb{R}^2)$, the Fourier and the inverse Fourier transform are defined by duality. We could view functions on the torus as periodic functions on the full sapce.

We use $(\Delta_i)_{i\geq -1}$ to denote the Littlewood–Paley blocks corresponding to a dyadic partition of unity. Besov spaces on the torus with general indices $\alpha \in \mathbb{R}, p, q \in[1, \infty]$ are defined as the completion of $C^\infty(\mathbb{T}^2)$ with respect to the norm
$$\|u\|_{B_{p,q}^\alpha}:=\Big{(}\sum_{j\geq-1}2^{j\alpha q}\|\Delta_ju\|_{L^p}^q\Big{)}^{1/q}.$$
The H\"{o}lder–Besov space $C^\alpha$ is given by $C^\alpha = B^\alpha_{\infty,\infty}$, and we also set $H^\alpha = B_{2,2}^\alpha,\ \alpha \in \mathbb{R}$.

The following embedding results will be frequently used on $\mathbb{T}^2$.
\bl\label{lem:2.2}
$(1).($\cite[Lemma A.2]{GIP15}$)$  Let $1\leq p_1\leq p_2\leq \infty$ and $1\leq q_1\leq q_2\leq \infty$, and let $\alpha\in\mathbb{R}$. Then $B_{p_1,q_1}^\alpha\subset B_{p_2,q_2}^{\alpha-2(1/p_1-1/p_2)}.$\\
$(2).($\cite[Theorem 4.6.1]{Tri78}$)$Let $s\in\mathbb{R},1<p<\infty,\epsilon>0$. Then $W^{s,2}=B_{2,2}^s=H^s$, and $B_{p,1}^s\subset W^{s,p}\subset B_{p,\infty}^s\subset B_{p,1}^{s-\epsilon}.$
\el

\subsection{Paraproducts}
Paraproducts were introduced by Bony in
\cite{Bon81} and they permit to decompose a product of two distributions into three parts which behave differently in terms of regularity. More precisely, using the Littlewood-Paley blocks, the product $fg$ of two Schwartz distributions $f,g \in \mathcal{S}' (\mathbb{R}^2)$ can be formally decomposed as
$$fg=f\prec g+f\succ g+f\circ g,$$
with 
$$f\prec g=g\succ f=\sum_{i\geq-1}\sum_{i<j-1}\Delta_if\Delta_jg,\ \ f\circ g=\sum_{|i-j|\leq1}\Delta_if\Delta_jg.$$
Here, the paraproducts $\prec$ and $\succ$ are always well-defined and critical is the resonant product denoted by $\circ$. In general, it is only well-defined provided the sum of the regularities of $f$ and $g$ in terms of Besov spaces is strictly positive. Moreover, we have the following paraproduct estimates. 
\bl\label{lem:2.3}$($\cite[Lemma 2.1]{GIP15},\cite[Proposition A.7]{MW17a}$)$
Let $\beta\in\mathbb{R},p,p_1,p_2,q\in[1,\infty]$ such that $\frac{1}{p}=\frac{1}{p_1}+\frac{1}{p_2}$. Then it holds 
$$\|f\prec g\|_{B_{p,q}^\beta}\lesssim\|f\|_{L^{p_1}}\|g\|_{B_{p_2,q}^\beta},$$
and if $\alpha<0$ then 
$$\|f\prec g\|_{B_{p,q}^{\alpha+\beta}}\lesssim\|f\|_{B_{p_1,q}^\alpha}\|g\|_{B_{p_2,q}^\beta}.$$
If $\alpha+\beta>0$ then it holds 
$$\|f\circ g\|_{B_{p,q}^{\alpha+\beta}}\lesssim\|f\|_{B_{p_1,q}^\alpha}\|g\|_{B_{p_2,q}^\beta}.$$
\el
We denote $\succeq=\circ+\succ,\ \preceq=\circ+\prec$.

Analogously to the real-valued case, we may define paraproducts for vector-valued distributions. More precisely, for two vector-valued distributions $f, g \in \mathcal{S}'(\mathbb{R}^2; \mathbb{R}^2)$, we use the following
tensor paraproduct notation
\begin{align*}
f \otimes g = (f_ig_j)^m_{i,j=1} = f \prec\!\!\!\!\!\!\!\bigcirc g + f \circledcirc\ g + f \succ\!\!\!\!\!\!\!\bigcirc \ g
= (f_i \prec g_j)_{i,j=1}^m + (f_i \circ g_j )^m_{i,j=1} + (f_i \succ g_j )^m_{i,j=1},    
\end{align*}
and note that Lemma \ref{lem:2.3} carries over mutatis mutandis. We also denote $\succcurlyeq\!\!\!\!\!\!\!\bigcirc=\circledcirc+\succ\!\!\!\!\!\!\!\bigcirc,\ \preccurlyeq\!\!\!\!\!\!\!\bigcirc=\circledcirc+\prec\!\!\!\!\!\!\!\bigcirc$.

We also recall the following lemma for the Helmholtz projection  $\mathbb{P}_{\rm H}=\Id-\nabla\Delta^{-1}\div$:
\bl\label{lem:est:leray}$($\cite[Lemma 3.6]{HZZ21b}$)$
Assume that $\alpha\in\mathbb{R}$ and $p \in [1, \infty]$. Then for every $k, l = 1, 2$ 
\begin{align*}
    \|\mathbb{P}_{\rm H}^{kl}f\|_{B_{p.\infty}^\alpha}\lesssim\|f\|_{B_{p.\infty}^\alpha}.
\end{align*}
\el

Finally, we introduce high and low frequency projections from \cite{HR23}. Compared with localizers in \cite{HZZ21b}, the following projections allows $J \in\mathbb{R}^+$. Let $J >0$. For
$f\in \mathcal{S}'(\mathbb{R}^2)$ we define the projections
 $$ \mathcal{H}_{J}f=\check{h}_J*f,\ \ \mathcal{L}_{ J}f=\check{l}_J*f,$$
where

\begin{align*}
\check{h}_J=\mathcal{F}^{-1}(h(|\cdot|/J)), \ \ \check{l}_J=\mathcal{F}^{-1}(l(|\cdot|/J)),  
\end{align*}
for smooth function $h:[0,\infty)\to[0,\infty)$ satisfying 
$$h(r)=1,\ \ {\rm{if}} \ r\geq1,\ \ \ h(r)=0,\ \ {\rm{if}}\ r\leq\frac12, \ \ l(r)=1-h(r).$$
 Here we view $f$ as periodic function on whole space.

Then we have the following result.
\bl$($\cite[Lemmas 4.2-4.3]{HR23}$)$
For $\alpha\leq\beta\leq \gamma$\begin{align}\|\cH_{J}f\|_{C^\alpha}\lesssim J^{-(\beta-\alpha)}\|f\|_{C^\beta},\ \ \|\cL_{ J}f\|_{C^\gamma}\lesssim J^{\gamma-\beta}\|f\|_{C^\beta}.   \label{deltaf} \end{align}
\el

\subsection{Anti-divergence operators}

 We recall the following anti-divergence operator  
 $\mathcal{R}$ from \cite[Appendix D]{CL20}, which acts on vector fields $v$ as
$$(\mathcal{R}v)_{ij}=\mathcal{R}_{ijk}v_k$$
where
$$\mathcal{R}_{i j k}=- \Delta^{-1} \partial_{k} \delta_{i j}+\Delta^{-1} \partial_{i} \delta_{j k}+\Delta^{-1} \partial_{j} \delta_{i k}.$$
Then $\mathcal{R}v(x)$ is a symmetric trace-free matrix for each $x \in \mathbb{T}^2$, and $\mathcal{R}$ is a right
inverse of the $\div$ operator, i.e. $\div(\mathcal{R}v) = v-\int_{\mathbb{T}^2}\!\!\!\!\!\!\!\!\!\!\!\!\!\!\; {}-{}\ v\dif x$
By a direct computation we have for any divergence-free $v\in C^\infty(\mathbb{T}^2;\mathbb{R}^2)$
\begin{align}\label{rdeltav}
\mathcal{R}\Delta v=\nabla v+\nabla^{T} v. 
\end{align}

Let $C_0^\infty(\mathbb{T}^2;\mathbb{R}^{2})$ be the set of smooth functions with zero mean. By \cite[Theorem D.2]{CL20} we know for any $1\leq p \leq\infty$, $\sigma\in\mathbb{N}$, $f\in C_0^\infty(\mathbb{T}^2;\mathbb{R}^{2})$  
\begin{align}
\|\mathcal{R}f(\sigma\cdot)\|_{ L^p}\lesssim\sigma^{-1} \| f\|_{ L^p}.\label{bb1}
\end{align}
 
Let $C_0^\infty(\mathbb{T}^2;\mathbb{R}^{2\times 2})$ be the set of smooth matrix valued functions with zero mean. We also introduce the bilinear version $\mathcal{B}: C^\infty(\mathbb{T}^2;\mathbb{R}^2) \times C^\infty_0(\mathbb{T}^2; \mathbb{R}^{2\times 2})\to C^\infty(\mathbb{T}^2; \mathcal{S}_0^{2\times 2})$ by $$(\mathcal{B}(v,A))_{ij}=v_l\mathcal{R}_{ijk}A_{lk}-\mathcal{R}(\partial_iv_l\mathcal{R}_{ijk}A_{lk}).$$

Then by \cite[Theorem D.3]{CL20} we have $\div(\mathcal{B}(v,A))=vA-{\int_{\mathbb{T}^2}\!\!\!\!\!\!\!\!\!\!\!\!\; {}-{}\ vA\dif x}$ and for any $1 \leq p \leq\infty$ 
\begin{align}
\|\mathcal{B}(v,A)\|_{L^p}\lesssim\|v\|_{C_x^1}\|\mathcal{R}A\|_{L^p}.    \label{bb}
\end{align}

\subsection{Probabilistic elements}\label{sec:pe}
Regarding the driving noise, we assume that $B$ is a  vector-valued  two-sided $L^2$-cylindrical Wiener process  with zero spatial mean and zero divergence on some stochastic basis   $(\Omega, \mathcal{F}, (\mathcal{F}_t)_{t\in\mathbb{R}}, \mathbf{P})$ (see e.g. \cite[page 99]{PR07}). The time derivative of $B$ is the space-time white noise. 

For $p\in[1,\infty),q\in[1,\infty],\alpha\in(0,\infty)$  and a Banach space $B$, we denote
\begin{align*}
    \$u \$ _{B,p}^p:=\sup_{t\in\mathbb{R}}\mathbf{E}\left[\sup_{s\in[t,t+1]}\|u(s)\|_{B}^p\right],\  \$ u \$ _{C_{t,x}^1,p}^p:=\sup_{t\in\mathbb{R}}\mathbf{E}\left[\|u\|_{C_{[t,t+1],x}^1}^p\right],
\end{align*}
\begin{align*}
    \$u \$ _{C_t^\alpha B,p}^p:=\sup_{t\in\mathbb{R}}\mathbf{E}\left[\|u\|_{C_{[t,t+1]}^\alpha B}^p\right],\  \$u \$ _{L_t^q B,p}^p:=\sup_{t\in\mathbb{R}}\mathbf{E}\left[\|u\|_{L^q_{[t,t+1]}B}^p\right].
\end{align*}

We  recall the stochastic objects needed for the rigorous formulation of \eqref{1.1}. Consider the linear equation \eqref{3.1} and we perform renormalization in order to define the term $z\otimes z$ in \eqref{3.2}. We denote by $z_{\epsilon}$ the spatial
mollification of $z$ and $z^{:2:}$ can be defined as the limit of 
 \begin{align}
     z^{:2:}_{\epsilon}= z_{\epsilon} \otimes z_{\epsilon}-\mathbf{E}[z_{\epsilon} \otimes z_{\epsilon}]  \label{z:2:}
 \end{align}
 in $L^p(\Omega;C(\mathbb{R};C^{-2\fa-\kappa}))$ for every $p\geq1,\kappa>0$.
 
 We refer to the following result that is analogous to \cite[Lemma 3.2, Theorem 5.1]{DPD02}.

\bp\label{prop:Z}
 Let $\fa_0:=\frac{1-3\fa}6\wedge \frac1{12}$. For any $ \kappa>0, q\geq2$ we have

\begin{align}
& \$ z \$ _{C^{-\fa-\kappa},q}+ \$ z \$ _{C_t^{\frac{\fa_0}{2}}C^{-\fa-\fa_0-\kappa},q}
\leq (q-1)^{1/2}L,\label{3.7}
\end{align}
\begin{align}
 \$ z^{:2:} \$ _{C^{-2\fa-\kappa},q}\leq (2q-1)L,\label{3.7*}
\end{align}
where $L$ is a constant.
\ep
 
In the following,  we fix $L\geq 1$ large enough such that Proposition \ref{prop:Z} holds.

\section{Analytic decomposition and estimate for stochastic objects} \label{cog}
In this section, we introduce a time-dependent decomposition of the equation \eqref{3.2}, which makes it possible to deduce a uniform moment bound for the solutions to \eqref{3.2}. More specifically let the equation of $v^1$ contains  the irregular terms  in \eqref{3.2} while  the regular terms are put in $v^2$: 
 
\begin{align}\label{3.3}
(\partial_t-\Delta+\mathfrak{c})v^1-z+\nabla p^1+\div(z^{:2:}+(v^1+v^2)\prec\!\!\!\!\!\!\!\bigcirc\cH_{R}z+\cH_{R}z\succ\!\!\!\!\!\!\!\bigcirc(v^1+v^2&))=0,\\
\div v^1=0.\notag
\end{align}
\begin{align}\label{3.4}
(\partial_t-\Delta)v^2-\mathfrak{c}v^1&+\nabla p^2+\div((v^1+v^2)\otimes\cL_{R}z+\cL_{R}z\otimes(v^1+v^2))\notag\\
=&-\div((v^1+v^2)\otimes(v^1+v^2)+(v^1+v^2)\succcurlyeq\!\!\!\!\!\!\!\bigcirc\cH_{R}z+\cH_{R}z\preccurlyeq\!\!\!\!\!\!\!\bigcirc(v^1+v^2)),\\
\div v^2&=0.\notag
\end{align}

It is easy to see that $v=v^1+v^2$ solves \eqref{3.2}.
 Compared with \cite{LZ23}, here we add a new damping term $\fc v^1$ in the first equation  with $\fc$ being a large constant given in \eqref{def:fC} below  and $R$
 in \eqref{3.3}, \eqref{3.4} depends on time, which  plays an important role in the estimate of $v^1$ (see Section \ref{est:vq1} for more explanation).  More precisely,
 \begin{align}
 R(t):=  (1+ \mathfrak{C}\|z(t)\|_{C^{-\fa-\fa_0-\kk}}) ^{6}\label{def:R}
\end{align}
with $\fC$ being a large constant given in \eqref{def:fC}.  Here we recall $\fa_0=\frac{1-3\fa}6\wedge \frac1{12}$ and $0<\kk<\frac13-\fa-\fa_0$. 
 The choice of $R(t)$ ensures the $\cH_{R}z$ part could be bounded by a constant, which avoids increasing of moments during the estimate of $v^1$. However, we have to pay a price to bound commutator errors from  $\cH_{R(s)}$ and $\cL_{ R(s)}$ in the estimate of  $v^2$ (see Section \ref{votiefr}). 
 To this end, we estimate the H\"older regularity of these stochastic objects with time-depend cut-off.
 \bp\label{deltaz}
For any  $0<\kappa<\frac13-\fa-\fa_0$, it holds for any $t\in\mathbb{R}$
\begin{align}
&\sup_{s\in[t,t+1]}\|\cL_{ R(s)}z(s)\|_{L^\infty}
\lesssim   \fC^3(1+\sup_{s\in[t,t+1]}\|z(s)\|_{C^{-\fa-\kappa}})^4,\notag\\
&  \sup_{s\in[t,t+1]}\|\cH_{ R(s)}z(s)\|_{C^{-\frac16-\fa-\fa_0-\kk}}\lesssim \frac1{\mathfrak{C}},\notag\\
 & \|\cL_{R}z\|_{C_{[t-1,t+1]}^{\frac{\fa_0}2}C^{\frac13}}\lesssim  \fC^{6}(1+\sup_{s\in[t-1,t+1]}\|z(s)\|_{C^{-\fa-\kappa}}) ^{6} \|z\|_{C_{[t-1,t+1]}^{\frac{\fa_0}{2}}C^{-\fa-\fa_0-\kk}},\notag
\end{align}
and
\begin{align*}
     \|\cH_{R}z\|_{C_{[t-1,t+1]}^{\frac{\fa_0}2}C^{-\fa-\fa_0-\kk}}\lesssim    \fC^{6}(1+\sup_{s\in[t-1,t+1]}\|z(s)\|_{C^{-\fa-\kappa}}) ^{6} \|z\|_{C_{[t-1,t+1]}^{\frac{\fa_0}{2}}C^{-\fa-\fa_0-\kk}}.
\end{align*}
 Here the implicit constants are universal and independent of $z$.
\ep

\begin{proof}
 By (\ref{deltaf}) and \eqref{def:R} we obtain for $0<\kk<\frac13-\fa$
  \begin{align}
\sup_{s\in[t,t+1]}\|\cL_{ R(s)}z(s)\|_{L^\infty}&\lesssim\sup_{s\in[t,t+1]}\|\cL_{ R(s)}z(s)\|_{C^{\frac16}}  \lesssim\sup_{s\in[t,t+1]} R(s)^{\frac12}\|z(s)\|_{C^{-\fa-\kk}}\notag\\
&\lesssim \fC^3(1+\sup_{s\in[t,t+1]}\|z(s)\|_{C^{-\fa-\kappa}})^4,\notag\\
  \sup_{s\in[t,t+1]}\|\cH_{ R(s)}z(s)\|_{C^{-\frac16-\fa-\fa_0-\kk}}&\lesssim   \sup_{s\in[t,t+1]} R(s)^{-\frac16}\|z(s)\|_{C^{-\fa-\fa_0-\kk}}\lesssim\frac1{\mathfrak{C}}.\notag
\end{align}
The last two inequalities need to be handled carefully. For any $t-1\leq s\leq u\leq t+1$  we have
\begin{align*}
   \|\cL_{ R(u)}z(u)-\cL_{ R(s)}z(s)\|_{C^{\frac13}}\leq\|\cL_{ R(u)}\big{(}z(u)-z(s)\big{)}\|_{C^{\frac13}}+\left\|\mathcal{F}^{-1}\(l(\frac{|\cdot|}{R(u)})-l(\frac{|\cdot|}{R(s)})\)*z(s)\right\|_{C^{\frac13}}.
\end{align*} 
For the first term  by \eqref{deltaf} and \eqref{def:R} we obtain for $0<\kappa<\frac13-\fa-\fa_0$
 \begin{align}
     \|\cL_{ R(u)}\big{(}z(u)-z(s)\big{)}\|_{C^{\frac13}}&\lesssim R(u)^{\frac23}\|z(u)-z(s)\|_{C^{-\fa-\fa_0-\kk}}\notag\\
   &\lesssim (1+\fC\sup_{s\in[t-1,t+1]}\|z(s)\|_{C^{-\fa-\fa_0-\kk}})^{4} \|z\|_{C_{[t-1,t+1]}^{\frac{\fa_0}{2}}C^{-\fa-\fa_0-\kk}}|u-s|^{\frac{\fa_0}{2}}.\label{est:deltaR1}
 \end{align}
For the second term we have
\begin{align*}
    l\(\frac{|\cdot|}{R(u)}\)-l\(\frac{|\cdot|}{R(s)}\)&=-\int_0^1l'\( \frac{|\cdot|}{\theta R(u)+(1-\theta)R(s)}\)\frac{|\cdot|(R(u)-R(s))}{(\theta R(u)+(1-\theta)R(s))^2}\dif \theta\\
    &=\int_0^1\frac{R(u)-R(s)}{\theta R(u)+(1-\theta)R(s)}l^*\(\frac{|\cdot|}{\theta R(u)+(1-\theta)R(s)}\)\dif \theta,
\end{align*}
where we define $l^*(x):=-l'(x)x$ with support on $[\frac12,1]$. Then using Young's inequality and \eqref{deltaf} we have for $0<\kappa<\frac13-\fa$
\begin{align}
&\left\|\mathcal{F}^{-1}\( l\(\frac{|\cdot|}{R(u)}\)-l\(\frac{|\cdot|}{R(s)}\)\)*z(s)\right\|_{C^{\frac13}}= \sup_{j} 2^{j/3}\left\|\mathcal{F}^{-1}\(l\(\frac{|\cdot|}{R(u)}\)-l\(\frac{|\cdot|}{R(s)}\)\)*\Delta_jz(s)\right\|_{L^\infty}\notag\\
  &\leq\sup_{j} 2^{ j/3}\int_0^1\frac{|R(u)-R(s)|}{\theta R(u)+(1-\theta)R(s)}\left\|\mathcal{F}^{-1}l^*\(\frac{|\cdot|}{\theta R(u)+(1-\theta)R(s)}\)*\Delta_jz(s)\right\|_{L^\infty}\dif \theta\notag\\
   &\lesssim  \sup_{j}2^{-(\kappa+\fa) j}\|\Delta_jz(s)\|_{L^\infty}|R(u)-R(s)|\notag\\
   &\ \ \ \ \ \ \ \ \ \ \times\int_0^1[\theta R(u)+(1-\theta)R(s)]^{\frac23-1}\left\|\mathcal{F}^{-1}\(l^*\(\frac{|\cdot|}{\theta R(u)+(1-\theta)R(s)}\)\)\right\|_{L^1}\dif \theta.\notag
   \end{align}
Here we observed that the support of $l^*(\frac{|\cdot|}{\theta R(u)+(1-\theta)R(s)})$ and $\mathcal{F}(\Delta_jz)$ intersect only when $\theta R(u)+(1-\theta)R(s)\simeq 2^j$. 
   Next, we have $$\left\|\mathcal{F}^{-1}\(l^*\(\frac{|\cdot|}{\theta R(u)+(1-\theta)R(s)}\)\)\right\|_{L^1}=\|\mathcal{F}^{-1}\big{(}l^*(|\cdot|)\big{)}\|_{L^1}<\infty,$$
   and
   \begin{align*}
       |R(u)-R(s)|&\lesssim \fC^{6}(1+\sup_{s\in[t-1,t+1]}\|z(s)\|_{C^{-\fa-\fa_0-\kappa}}) ^{5} \|z(u)-z(s)\|_{C^{-\fa-\fa_0-\kk}}\\
      & \lesssim \fC^{6}(1+\sup_{s\in[t-1,t+1]}\|z(s)\|_{C^{-\fa-\kappa}}) ^{5} \|z\|_{C_{[t-1,t+1]}^{\frac{\fa_0}{2}}C^{-\fa-\fa_0-\kk}}|u-s|^{\frac{\fa_0}{2}}.
   \end{align*}
 Hence by $R\geq1$ we deduce that
   \begin{align}
   & \left\|\mathcal{F}^{-1}\(l\(\frac{|\cdot|}{R(u)}\)-l\(\frac{|\cdot|}{R(s)}\)\)*z(s)\right\|_{C^{\frac13}}\notag\\ 
    &\ \ \ \ \ \ \ \lesssim  \fC^{6}(1+\sup_{s\in[t-1,t+1]}\|z(s)\|_{C^{-\fa-\kappa}}) ^{6} \|z\|_{C_{[t-1,t+1]}^{\frac{\fa_0}{2}}C^{-\fa-\fa_0-\kk}}|u-s|^{\frac{\fa_0}{2}},\label{meanvalue}
\end{align}
which combined with \eqref{est:deltaR1} implies the third inequality.

 The last inequlaity  follows from \eqref{deltaf} and \eqref{meanvalue}: for any $t-1\leq s\leq u\leq t+1$
\begin{align*}
    \|\cH_{ R(u)}z(u)&-\cH_{ R(s)}z(s)\|_{C^{-\fa-\fa_0-\kk}}\\
    &\leq\|\cH_{ R(u)}\big{(}z(u)-z(s)\big{)}\|_{C^{-\fa-\fa_0-\kk}}+\|(\cL_{ R(s)}-\cL_{ R(u)})z(s)\|_{ C^{\frac13}}\\
    &\leq \fC^{6}(1+\sup_{s\in[t-1,t+1]}\|z(s)\|_{C^{-\fa-\kappa}}) ^{6} \|z\|_{C_{[t-1,t+1]}^{\frac{\fa_0}{2}}C^{-\fa-\fa_0-\kk}}|u-s|^{\frac{\fa_0}{2}}.
\end{align*} 
\end{proof}

\section{ Analytic construction}\label{sec:adc}

In this section we apply the convex integration method to the nonlinear equation (\ref{3.4}) indexed by a parameter $q\in\mathbb{N}_0$, and  consider it in conjunction with equation \eqref{3.3} through a joint iterative process, following the methodology outlined in \cite{HZZ21b,LZ23}. We choose an increasing sequence $\{\lambda_q\}_{q\in\mathbb{N}_0}\subset\mathbb{N}$ which diverges to $\infty$,  alongside a sequence $\{\delta_q\}_{q\in\mathbb{N}}\subset (0,1)$
 which is decreasing to 0. For some $a\in\mathbb{N},b\in\mathbb{N},\beta\in(0,1)$ we let
$$\lambda_q=a^{(b^q)},\ \ q\in\mathbb{N}_0,\ \ \delta_q=\frac{1}{2}\lambda_1^{2\beta}\lambda_q^{-2\beta},\ q\in\mathbb{N}, \ \ \delta_0=1.$$
Here $\beta$ will be chosen sufficiently small and $a$ as well as $b$ will be chosen sufficiently large.  We introduce a constant $\aleph\in(0,1)$ to absorb the universal constants in \eqref{leqaleph} and \eqref{4.1*}  below, which  is also used  to prove the constructed solutions  are not Gaussian $N(0,\frac12(-\Delta)^{\fa-1})$  (see Proposition \ref{6.1} below).
 In addition, we introduce a small parameter $\alpha\in(0,1)$ and used $ \sum_{q\geq1} \delta_q^{1/2}\leq  \frac{1}{\sqrt2}\sum_{q\geq1}a^{b\beta-qb\beta}\leq\frac{1}{\sqrt2} \frac{1}{1-a^{-b\beta}}\leq1,\sum_{q\geq1}\lambda_{q}^{-\alpha/2}\leq\frac1{a^{\alpha b/2}-1}\leq \aleph$ 
which boils down to 
\begin{align}\label{ab2}
a^{b\beta}\geq2+\sqrt2,\ \ a^{\alpha b/2}\geq1+\frac1\aleph,
\end{align}
 which we assumed from now on.  More details on the choice of these parameters will be given below in the course of the construction. 

At each step $q$, a triple $(v_q^1, v_q^2,\mathring{R}_q)$ is constructed solving the following system on $[0,\infty)$:

\begin{align}
(\partial_t-\Delta+\mathfrak{c})v_q^1-z+\nabla p_q^1&+\div(z^{:2:}+(v_q^1+v_q^2)\prec\!\!\!\!\!\!\!\bigcirc\cL_{ f(q)}\cH_{R}z+\cL_{ f(q)}\cH_{R}z\succ\!\!\!\!\!\!\!\bigcirc(v_q^1+v_q^2))=0,\notag\\
\div v_q^1&=0,\ \ \ v_q^1(0)=0.\label{3.5}
\end{align}
\begin{align}
(\partial_t&-\Delta)v_q^2-\mathfrak{c}v_q^1+\nabla p_q^2+\div((v_q^1+v_q^2)\otimes\cL_{R}z+\cL_{R}z\otimes(v_q^1+v_q^2))\notag\\
&=-\div((v_q^1+v_q^2)\otimes(v_q^1+v_q^2)+(v_q^1+v_q^2)\succcurlyeq\!\!\!\!\!\!\!\bigcirc\cH_{R}z+\cH_{R}z\preccurlyeq\!\!\!\!\!\!\!\bigcirc(v_q^1+v_q^2))+\div \mathring{R}_q,\notag\\ 
&\div v_q^2=0,\label{3.6}
\end{align}
where we defined $f(q)=\lambda_q^{\alpha/2}$ for $\alpha\in(0,1)$ small enough. Here we add the localizers $\cL_{ f(q)}$ in the equation of $v_{q+1}^1$, which are
used to control the blow up of a certain norm of $v_{q}^1$ as $q\to\infty$. Note that the Reynolds stress $\mathring{R}_q$ is only included in the equation for $v_q^2$, and $v_q^2$ contains the quadratic nonlinearity  to reduce the  stress.  Furthermore, given a specific $v_q^2$, we establish the existence and uniqueness of a solution $v_q^1$ to the linear equation (\ref{3.5}) using a fixed point argument, together with the uniform estimate derived subsequently.

Similarly to \cite{HLZZ23,HZZ22c},  we work in the framework of stationarity understood with respect
to shifts on trajectories. In the upcoming discussion, we will construct the triple $(v_q^1, v_q^2, \mathring{R}_q)$ over the entire domain $\mathbb{R}$ in time. We define $v_q^1(t)=0$ for $t\in(-\infty,0)$ and  we only assume $(v_q^1,v_q^2)$  satisfies \eqref{3.5} on the internal $[0,\infty)$, not on $\mathbb{R}$. \eqref{3.6} is satisfied  for $t\in\mathbb{R}$ to handle the mollification in the convex integration below. 

 We intend to start the iteration from  a given divergence-free, mean-zero function  $v_0\in C^1(\mathbb{T}^2;\R^2)$ and denote
\begin{align}
\fm:=\|v_0\|_{C^1}\vee1.\label{omegan}
\end{align}
We introduce a parameter $r\geq1$ and in the following we will prove $2r{\rm th}$-moment of solutions are finite. Let $\bar{M}_L:= (M^{*}r(L\vee \fm)\aleph^{-1})^{\frac{3}{1-3\fa}}$ and $M_L:=2( \bar{M}_L r)^3$, where $M^*\geq1$ is the implicit constant from \eqref{est:r0t}  below. Then we  choose $\fc,\fC$ in \eqref{3.5}-\eqref{3.6} as follows:
\begin{align}
 \mathfrak{c}=(rL\aleph^{-1})^{ \frac6{1-3\fa}},\ \  \fC=\fm\aleph^{-1}.\label{def:fC} 
\end{align}

The key result is the following iterative proposition, which we prove in Section \ref{proofof3.1} below.

\bp\label{prop:3.1}
Let $ \fa\in[0,\frac13),r\geq1, \aleph\in(0,1)$ be fixed. For any given divergence-free, mean-zero $v_0\in C^1(\mathbb{T}^2;\R^2)$, there exists a choice of parameters $a,b,\alpha,\beta$ such
that the following holds true: Let $(v_q^1, v_q^2,\mathring{R}_q)$ be an $(\mathcal{F}_t)_{t\in\mathbb{R}}$-adapted solution to (\ref{3.5}) and (\ref{3.6})  on $[0,\infty)$ satisfying
\begin{align}
 \$ v_q^2 \$ _{L^2_tL^2,2r}\leq M_0M_L^{1/2}\sum_{i=0}^q\delta_i^{1/2},
 \label{3.8.1}
\end{align}
for a universal constant $M_0  \geq1$, and for $m\geq1$
 \begin{align}
 \$ v_q^2 \$ _{L^2_tL^2,m}\leq M_0(7^{q-1}\cdot 3 \bar{M}_L m)^{9\cdot {7}^{q-1}},
 \label{3.8.2}
\end{align}
 and
\begin{align}
& \$ v_q^2 \$ _{C_{t,x}^1,m}\leq\lambda_q^{7/2}(7^{q-1}\cdot 5 \bar{M}_L m)^{15\cdot {7}^{q-1}},\ \  \$ v_{q}^2 \$ _{C_{t,x}^1,2r}
\leq  \fm\lambda_{q}^4,
\label{3.9}\\
& \$ v_q^2 -v_0 \$ _{W^{\frac12,\frac65},2r}\leq \sum_{i=1}^q\lambda_{i}^{-\alpha/2}\leq  \aleph,\ \  \$ v_q^2 \$ _{W^{\frac12,\frac65},m}
\leq (7^{q-1}\cdot \frac{7}2 \bar{M}_L m)^{\frac{21}2\cdot {7}^{q-1}}, \label{3.10}\\
& \$ \mathring{R}_{q} \$ _{L^1_tL^1,r}\leq M_L\delta_{q+1},\ \ \
 \$ \mathring{R}_{q} \$ _{L^1_tL^1,m}\leq (7^{q}\cdot  \bar{M}_L m)^{3\cdot {7}^q}.\label{3.11}
\end{align}

Then there exists an $(\mathcal{F}_t)_{ t\in\mathbb{R}}$-adapted process $(v_{q+1}^1,v_{q+1}^2, \mathring{R}_{q+1})$ which solves  (\ref{3.5}) and (\ref{3.6}) at the level $q+1$  and obeys \eqref{3.8.1}-\eqref{3.11} at the level $q+1$. Moreover, it satisfies
\begin{align}
& \$ v_{q+1}^2-v_q^2 \$ _{L^2_tL^2,2r}\leq M_0\delta_{q+1}^{1/2}M_L^{1/2},
\label{3.15}\\
& \$ v_{q+1}^2-v_q^2 \$ _{W^{\frac12,\frac65},2r}\leq\lambda_{q+1}^{-\alpha/2}\leq\delta_{q+1}^{1/2},\ \ \$ v_{q+1}^2-v_q^2 \$ _{W^{\frac12,\frac65},m}\leq \lambda_{q+1}^{-\alpha/2}(7^{q}\cdot \frac{7}2 \bar{M}_L m)^{\frac{21}2\cdot {7}^q}.\label{3.16}
\end{align}
\ep

It's noteworthy that no bounds on $v_q^1,v_{q+1}^1$ were included in the statement of Proposition \ref{prop:3.1}. Indeed,
the definition of the new velocity $v_{q+1}^2$ does not require $v_{q+1}^1$. Then having $v_{q+1}^2$ at hand, all the
necessary bounds for $v_{q+1}^1$ can be deduced from  Section \ref{est:vq1} below. Specifically, we demonstrate the following.

\bp\label{prop:3.2}
Under the assumptions of Proposition \ref{prop:3.1},  it holds for $\kappa\in(0,\frac13-\fa-\fa_0),m\geq1$
\begin{align}
 \$ v_{q}^1 \$ _{L^2,2r}\lesssim \$ v_q^1 \$ _{B^{\frac23-\fa-\kappa}_{\frac32,\infty},2r}
&\lesssim \aleph\lesssim1,\label{3.13.1}\\
 \$ v_{q}^1 \$ _{L^2,m}\lesssim \$ v_q^1 \$ _{B^{\frac23-\fa-\kappa}_{\frac32,\infty},m}&\lesssim  \$ z \$ _{C^{-\fa-\kappa},m}+ \$ z^{:2:} \$ _{C^{-2\fa-\kappa},m}+ \$ v_q^2 \$ _{W^{\frac12,\frac65},m}\notag\\ 
&\lesssim(7^{q-1}\cdot \frac{7}2 \bar{M}_L m)^{\frac{21}2\cdot {7}^{q-1}},\label{3.13.2}\\
     \$ v_{q+1}^1-v_q^1 \$ _{L^2,2r}&\lesssim   \$ v_{q+1}^1-v_q^1 \$ _{B_{\frac32,\infty}^{\frac23-\fa-\kappa},2r}\lesssim {\lambda_{q}^{-\frac{\alpha}{24}}
    }\leq \delta_{q+1}^{1/2},\label{3.12}\\ 
  \$ v_q^1 \$ _{C_t^{\frac16}L^\infty,2r}+   \$ v_q^1 \$ _{C^{\frac13},2r}&\lesssim \lambda_q^{2\alpha}r^{1/2}L.\label{4.4.1}
\end{align}
\ep

Here the implicit constants are independent of $\fc,\fC,\fm$ and $L$.

 We start the iteration from $v_0^2(t)\equiv v_0$ for $t\in\R$. Then (\ref{3.8.1})-(\ref{3.10}) hold by the choice of parameters. In that case,  we  define $v_0^1(t)=0$ for $t<0$, and determine $v_0^1$ for $t\geq0$ by solving (\ref{3.5}). By \eqref{rdeltav}, $\mathring{R}_0$ is the trace-free part of the matrix

\begin{align}
-\nabla v_0^2-\nabla^Tv_0^2&-\fc\mathcal{R}{v_0^1}+(v_0^1+v_0^2){\otimes}\cL_{R}z+\cL_{R}z{\otimes}(v_0^1+v_0^2)\notag\\
&+(v_0^1+v_0^2){\otimes}(v_0^1+v_0^2)+(v_0^1+v_0^2)\ {\succcurlyeq\!\!\!\!\!\!\!\bigcirc}\cH_{R}z+\cH_{R}z\ {\preccurlyeq\!\!\!\!\!\!\!\bigcirc}(v_0^1+v_0^2).\notag
\end{align}

By \eqref{3.7}, \eqref{3.7*}, \eqref{omegan} and \eqref{3.13.2}  we obtain
 \begin{align*}
      \$ v_0^1 \$ _{L^2,2m}\lesssim \$ v_0^1 \$ _{B^{\frac23-\fa-\kappa}_{\frac32,\infty},2m}\lesssim  \$ z \$ _{C^{-\fa-\kappa},2m}+ \$ z^{:2:} \$ _{C^{-2\fa-\kappa},2m}+\$ v_0^2 \$ _{W^{\frac12,\frac65},m}\lesssim mL+ \fm.
 \end{align*}
Together with Lemmas \ref{lem:2.2}, \ref{lem:2.3} with $0<\kappa<\frac13-\fa-\fa_0$, \eqref{bb1} and Proposition \ref{deltaz} we have 
\begin{align}
 \$ \mathring{R}_{0} \$ _{L^1_tL^1,m}&\lesssim  \$ v_0^2 \$ _{C_{t,x}^1,m}+  \fc\$ v_0^1 \$ _{L^2,m}+ \$ v_0^1 +v_0^2\$ _{L^2,2m}^2\notag\\ 
 &\ \ \ \ \  +(\$ v_0^1 \$ _{B^{\frac23-\fa-\kappa}_{\frac32,\infty},2m}+\$ v_0^2\$_{C_{t,x}^1,2m})( \$ \cL_{R}z \$ _{L^\infty,2m}+ \$ z \$ _{C^{-\fa-\kappa},2m})\notag\\
 &\lesssim (\$ v_0^1 \$ _{B^{\frac23-\fa-\kappa}_{\frac32,\infty},2m}+\$ v_0^2\$_{C_{t,x}^1,2m}) (\fC^3\$ z \$ _{C^{-\fa-\kappa},2m}^4+\fC^3+\fc)+\$ v_0^1 +v_0^2\$ _{L^2,2m}^2\notag\\
&\lesssim ( mL+ \fm)(\fC^3 L^4m^{2}+\fc+ mL+ \fm)\leq  m^3(M^{*}r(L\vee \fm)\aleph^{-1})^{\frac{9}{1-3\fa}}=( \bar{M}_L m)^3,\label{est:r0t}
\end{align}
where $M^*\geq1$ is the implicit constant. Then by the choice of $M_L$ we deduce that
\begin{align*}
 \$ \mathring{R}_0 \$ _{L^1_tL^1,r}
&\leq ( \bar{M}_L r)^3= \frac1{2}M_L. 
\end{align*}
So (\ref{3.11}) is satisfied at the level $q = 0$, since $\delta_1 = \frac{1}{2}$.   

We deduce the following result by using Propositions \ref{prop:3.1} and \ref{prop:3.2}.

\bt\label{thm:guina}

Let $\fa\in[0,\frac13), r\geq1,\aleph\in(0,1)$.  For any divergence-free, mean-zero $v_0\in C^1(\mathbb{T}^2 ;\R^2)$,
there exist  $\kk>0$ small enough, $\zeta\in(0,1)$ and an $(\mathcal{F}_t)_{ t\in\mathbb{R}}$-adapted process $(v^1,v^2)$ such that 
\begin{align*}
v^1&\in L^{2r}(\Omega; C(\mathbb{R};H^{\zeta})\cap C^\zeta(\mathbb{R}; B_{\frac32,\infty}^{\frac23-\fa-\kappa})),\\
    v^2&\in L^{2r}(\Omega;
 L^{2}(\mathbb{R}; H^{\zeta} )  \cap C^\zeta(\mathbb{R}; W^{\frac12,\frac65}))
\end{align*}
 and is an analytically weak solution to (\ref{3.3}) and (\ref{3.4}) on $[0,\infty)$ with   $z(t)=\int_{-\infty}^te^{(t-s)(\Delta-I)}(-\Delta)^{\fa/2}\dif B_s$ and $z^{:2:}$ is defined through \eqref{z:2:}. Then $u:=v+z:=v^1+v^2+z$ is an analytically weak solution to (\ref{1.1})  on $[0,\infty)$. Moreover,
the solution satisfies
\begin{align}
 \$ v \$ _{L_t^{2}H^{\zeta},2r}+ \$ v \$ _{C_t^\zeta W^{(\frac23-\fa-\kappa)\wedge\frac12,\frac65},2r} <\infty,\label{eq:guina1}
\end{align}
and \begin{align}
 \$ v -v_0 \$ _{W^{(\frac23-\fa-\kappa)\wedge\frac12,\frac65},2r}\leq \sqrt{\aleph}.\label{eq:guina3}
\end{align}

\et

\begin{proof}
Starting from  $v_0^2(t)\equiv v_0$, we repeatedly apply Proposition \ref{prop:3.1} and obtain $(\mathcal{F}_t)_{ t\in\mathbb{R}}$-adapted processes $(v_q^1, v_q^2,\mathring{R}_q),q\in\mathbb{N}$. By (\ref{3.12}), (\ref{4.4.1}) and interpolation we obtain
\begin{align*}
    \sum_{q\geq0} \$ v_{q+1}^1-v_q^1 \$ _{H^{\zeta},2r}&\lesssim \sum_{q\geq0} \$ v_{q+1}^1-v_q^1 \$ _{L^2,2r}^{1-4\zeta} \$ v_{q+1}^1-v_q^1 \$ _{H^{1/4},2r}^{4\zeta}\\
    &\lesssim  \sum_{q\geq0}\delta_{q+1}^{\frac{1-4\zeta}2}\lambda_{q+1}^{8\zeta\alpha}(r^{1/2}L)^{4\zeta}\lesssim \sum_{q\geq0}\delta_{q+1}^{\frac{1-4\zeta-8\zeta/\beta}{2}}r^{1/2}L<\infty,
\end{align*}
where we chose $0<\zeta<\frac{\beta}{8+4\beta}$. Similarly we have
 \begin{align*}
    \sum_{q\geq0} \$ v_{q+1}^1-v_q^1 \$ _{C_t^{\zeta}{B_{\frac32,\infty}^{\frac23-\fa-\kappa}},2r}<\infty
\end{align*}
by choosing  $0<\zeta<\frac{\beta}{12+6\beta}$.
As a result, the limit $v^1=\lim_{q\to\infty}v_q^1$ exists and lies in $
L^{2r}(\Omega;C(\mathbb{R};H^{\zeta})\cap C^\zeta(\mathbb{R};B_{\frac32,\infty}^{\frac23-\fa-\kappa}))$.  Moreover together with (\ref{3.9}), (\ref{3.15}), \eqref{3.16}, interpolation  and H\"older's inequality we obtain 
\begin{align}
  \sum_{q\geq0} \$ v_{q+1}^2-v_q^2 \$ _{L^{2}_tH^{\zeta},2r}
 &\leq  \sum_{q\geq0} \$ v_{q+1}^2-v_q^2 \$ _{L^2_tL^2,2r}^{1-\zeta} \$ v_{q+1}^2-v_q^2 \$ _{ C_{t,x}^1,2r}^{\zeta}\notag\\
&\lesssim  \sum_{q\geq0}M_0M_L^{1/2}\delta_{q+1}^{1-\zeta}\lambda_{q+1}^{4\zeta}\lesssim M_0 M_L^{1/2} \sum_{q\geq0}\delta_{q+1}^{\frac{(1-\zeta)\beta-4\zeta}{2\beta}}<\infty, \notag
\end{align}
and 
\begin{align}
 \sum_{q\geq0} \$ v_{q+1}^2-v_q^2 \$ _{C_t^\zeta W^{\frac12,\frac65},2r}<\infty,\notag
\end{align}
where we chose $0<\zeta<\frac{\beta}{\beta+4 }$. As a result, we obtain the limit $v^2=\lim_{q\to\infty}v_q^2$ exists and lies in $L^{2r}(\Omega;L^{2}(\mathbb{R}; H^{\zeta} )\cap C^\zeta(\mathbb{R}; W^{\frac12,\frac65}))$. Then  by Sobolev embedding \eqref{eq:guina1} follows.

Furthermore by (\ref{3.11}) we obtain $\mathring{R}_{q}\to0$ in $L^{2r}(\Omega;L^1( \mathbb{R};L^1))$ as $q\to\infty$. Thus $(v^1,v^2)$ satisfies equations (\ref{3.3}) and (\ref{3.4})  on $[0,\infty)$ in the analytic weak sense.  Hence $u=v^1+v^2+z$ defined above solves (\ref{1.1}).
 
Moreover, by \eqref{3.10} and \eqref{3.13.1}  we obtain \eqref{eq:guina3}. In fact, by Sobolev embedding we have
\begin{align}
     \$ v -v_0\$ _{W^{(\frac23-\fa-\kappa)\wedge\frac12,\frac65},2r}\lesssim \$ v^1 \$ _{B_{\frac32,\infty}^{\frac23-\fa-\kk},2r}+ \$ v^2 -v_0\$ _{W^{\frac12,\frac65},2r}\lesssim\aleph\leq \sqrt{\aleph}.\label{leqaleph}
\end{align}
Here we choose $\aleph>0$ small enough to absorb the implicit constant.
\end{proof}

\section{Estimate of $v_q^1$}\label{est:vq1}
In this section, we prove Proposition \ref{prop:3.2} under the assumptions of Proposition \ref{prop:3.1}. Moreover, we recall that the equation for $v_q^1$ is linear. Consequently, for a given  $v_q^2$ we obtain the existence and uniqueness of solution $v_q^1$ to (\ref{3.5})  on $[0,\infty)$ with initial condition $v_q^1(0)=0$ by a fixed point argument together with the uniform estimate derived in the sequel.  Also, if $v_q^2$ is $(\mathcal{F}_t)_{ t\in\mathbb{R}}$-adapted, so is $v_q^1$. This in particular gives the existence of $v_{q+1}^1$ in Proposition \ref{prop:3.1}, once the new velocity $v_{q+1}^2$ was constructed in Section \ref{proofof3.1}.

In the following,  we make use of the localizers $\cH_{R}$ present in the equation for $v_q^1$ in (\ref{3.3}).
Namely, by selecting an appropriate choice of $R$ in \eqref{def:R} we can always apply Proposition \ref{deltaz} to absorb $z(s)$ in the
estimates and then take expectation  without increasing the moment when using H\"older's inequality. In the proof we estimate $v_q^1$ by induction on time intervals.
 This crucially relies on our  capacity  to choose $\aleph>0$ sufficiently small enough to ensure  $\fc,\fC$ large enough.

Now we only need to prove the estimate on $[0,\infty)$ as $v_q^1(t)=0$ for any $t<0$.

\subsection{Proof of  (\ref{3.13.1}) and (\ref{3.13.2})}\label{4.13-14}
We write (\ref{3.5}) in the mild sense. For any $N\in\mathbb{N}_0,s\in[N,N+2]$ we have

\begin{align}
v_q^1(s)&=e^{(s-N)(\Delta-\mathfrak{c})}v_q^1(N)+\int_{N}^se^{(s-u)(\Delta-\mathfrak{c})}\mathbb{P}_{\rm{H}}[z-\div(z^{:2:})]\dif u\notag\\
&\ \ \ -\int_{N}^se^{(s-u)(\Delta-\mathfrak{c})}\mathbb{P}_{\rm{H}}\div((v_q^1+v_q^2)\prec\!\!\!\!\!\!\!\bigcirc\cL_{ f(q)}\cH_{ R}z+\cL_{ f(q)}\cH_{R}z\succ\!\!\!\!\!\!\!\bigcirc(v_q^1+v_q^2))\dif u\notag\\ 
&=e^{(s-N)(\Delta-\mathfrak{c})}v_q^1(N)+\int_{0}^{s-N}e^{(s-N-u)(\Delta-\mathfrak{c})}\mathbb{P}_{\rm{H}}[z-\div(z^{:2:})](u+N)\dif u\notag\\
&\ \ -\int_{0}^{s-N}e^{(s-N-u)(\Delta-\mathfrak{c})}\mathbb{P}_{\rm{H}}\div[(v_q^1+v_q^2)\prec\!\!\!\!\!\!\!\bigcirc\cL_{ f(q)}\cH_{ R}z+\cL_{ f(q)}\cH_{R}z\succ\!\!\!\!\!\!\!\bigcirc(v_q^1+v_q^2)](u+N)\dif u,\label{4.0}
\end{align}
where $\mathbb{P}_{\rm{H}}$ is the Helmholtz projection.

By Lemma \ref{lem:2.3}, Lemma \ref{lem:est:leray}, Lemma \ref{C.3} and (\ref{deltaf})  we
have for $0<\kk<\frac13-\fa$  and $s\in[N,N+2]$
\begin{align}
\|v_q^1(s)\|_{B_{\frac32,\infty}^{\frac23-\fa-\kappa}}&\lesssim e^{-(s-N)\mathfrak{c}}\|v_q^1(N)\|_{B_{\frac32,\infty}^{\frac23-\fa-\kappa}}+\mathfrak{c}^{-\frac{1/3-\fa}{2}}\sup_{s\in[N,N+2]}(\|z(s)\|_{C^{-\fa-\kappa}}+\|z^{:2:}(s)\|_{C^{-2\fa-\kappa}})\notag\\
&+\sup_{s\in[N,N+2]}\|v_q^1(s)+v_q^2(s)\|_{L^\frac32}\|\cH_{ R(s)}z(s)\|_{C^{-\frac13-\fa-\kappa}},\notag
\end{align}
where the implicit constant is independent of $N$. Thus together with Proposition \ref{deltaz} and Sobolev embedding we have for $s\in[N,N+2]$
\begin{align}
\|v_q^1(s)\|_{B_{\frac32,\infty}^{\frac23-\fa-\kappa}}&\lesssim e^{-(s-N)\fc}\|v_q^1(N)\|_{B_{\frac32,\infty}^{\frac23-\fa-\kappa}}+\fc^{-\frac{1/3-\fa}{2}}\sup_{s\in[N,N+2]}(\|z(s)\|_{C^{-\fa-\kappa}}+\|z^{:2:}(s)\|_{C^{-2\fa-\kappa}})\notag\\
&+ \frac1{\mathfrak{C}}\sup_{s\in[N,N+2]}\|v_q^2(s)\|_{W^{\frac12,\frac65}}
+ \frac1{\mathfrak{C}}\sup_{s\in[N,N+2]}\|v_q^1(s)\|_{B_{\frac32,\infty}^{\frac23-\fa-\kappa}}.\notag
\end{align}
Therefore we obtain the following: for every   $m\geq1$

\begin{align}
\mathbf{E}[&\sup_{s\in[N+1,N+2]}\|v_q^1(s)\|_{B_{\frac32,\infty}^{\frac23-\fa-\kappa}}^m]^{\frac{1}{m}}\leq M_1\Big{(}(e^{-\fc}+\frac1{\mathfrak{C}})\mathbf{E}[\sup_{s\in[N,N+1]}\|v_q^1(s)\|_{B_{\frac32,\infty}^{\frac23-\fa-\kappa}}^{m}]^{\frac{1}{m}}\notag\\
&+\fc^{-\frac{1/3-\fa}{2}}(\$ z \$ _{C^{-\fa-\kappa},m}+ \$ z^{:2:} \$ _{C^{-2\fa-\kappa},m})+\frac1{\mathfrak{C}} \$ v_q^2 \$ _{W^{\frac12,\frac65},m}+ \frac1{\mathfrak{C}}\mathbf{E}[\sup_{s\in[N+1,N+2]}\|v_q^1(s)\|_{B_{\frac32,\infty}^{\frac23-\fa-\kappa}}^m]^{\frac{1}{m}}\Big{)}.\notag
\end{align}
where $M_1$ is a universal constant independent of $N$.
Then we choose  $\mathfrak{C},\fc$ large enough such that 
\begin{align}
\mathbf{E}[\sup_{s\in[N+1,N+2]}&\|v_q^1(s)\|_{B_{\frac32,\infty}^{\frac23-\fa-\kappa}}^m]^{\frac{1}{m}}
\leq \frac{1}{2} \mathbf{E}[\sup_{s\in[N,N+1]}\|v_q^1(s)\|_{B_{\frac32,\infty}^{\frac23-\fa-\kappa}}^{m}]^{\frac{1}{m}}\notag\\ 
&+\frac{8}7M_1\left(\fc^{-\frac{1/3-\fa}{2}}(\$ z \$ _{C^{-\fa-\kappa},m}+ \$ z^{:2:} \$ _{C^{-2\fa-\kappa},m})+\frac1{\mathfrak{C}} \$ v_q^2 \$ _{W^{\frac12,\frac65},m}\right),\label{4.1*}
\end{align}
 Moreover, following a similar argument as above, and given that $v_q^1(0)=0$, we can deduce that
\begin{align}
\mathbf{E}[\sup_{s\in[0,1]}&\|v_q^1(s)\|_{B_{\frac32,\infty}^{\frac23-\fa-\kappa}}^m]^{\frac{1}{m}}\leq \frac{8}7M_1\left(\fc^{-\frac{1/3-\fa}{2}}( \$ z \$ _{C^{-\fa-\kappa},m}+ \$ z^{:2:} \$ _{C^{-2\fa-\kappa},m})+ \frac1{\mathfrak{C}} \$ v_q^2 \$ _{W^{\frac12,\frac65},m}\right).\notag
\end{align}
Thus by induction we obtain  for any $N\in\mathbb{N}_0$
\begin{align}
     \mathbf{E}[\sup_{s\in[N,N+1]}\|v_q^1(s)\|_{B_{\frac32,\infty}^{\frac23-\fa-\kappa}}^m]^{\frac{1}{m}}
 \leq3M_1\left(\fc^{-\frac{1/3-\fa}{2}}( \$ z \$ _{C^{-\fa-\kappa},m}+ \$ z^{:2:} \$ _{C^{-2\fa-\kappa},m})+\frac1{\mathfrak{C}} \$ v_q^2 \$ _{W^{\frac12,\frac65},m}\right).\label{4.1}
\end{align}
 As $v_q^1(t)=0$ for $t<0$, together with the embedding $B_{\frac32,\infty}^{\frac23-\fa-\kappa}\subset L^2$ with $\frac23-\fa-\kappa>\frac13$,  (\ref{3.7}), \eqref{3.7*}, \eqref{omegan} and (\ref{3.10}) we get
\begin{align}
  \$ v_q^1 \$ _{L^2,m} \lesssim   \$ v_q^1 \$ _{B_{\frac32,\infty}^{\frac23-\fa-\kappa},m}
&\lesssim  \$ z \$ _{C^{-\fa-\kappa},m}+ \$ z^{:2:} \$ _{C^{-2\fa-\kappa},m}+ \$ v_q^2 \$ _{W^{\frac12,\frac65},m}\notag\\ 
 & \lesssim(7^{q-1}\cdot \frac{7}2 \bar{M}_L m)^{\frac{21}2\cdot {7}^{q-1}}.\notag
\end{align}
 Together with   (\ref{3.7}), \eqref{3.7*},
 \eqref{omegan}, \eqref{3.10}  and the choice of $\fc,\fC$ we obtain
\begin{align}
  \$ v_q^1 \$ _{ L^2,2r}\lesssim   \$ v_q^1 \$ _{B_{\frac32,\infty}^{\frac23-\fa-\kappa},2r}&\lesssim  \fc^{-\frac{1/3-\fa}{2}}( \$ z \$ _{C^{-\fa-\kappa},2r}+ \$ z^{:2:} \$ _{C^{-2\fa-\kappa},2r})+ \frac1{\mathfrak{C}} \$ v_q^2 \$ _{W^{\frac12,\frac65},2r}\notag\\ 
 &\lesssim  \fc^{-\frac{1/3-\fa}{2}}  rL+ \frac{\fm+1}{\mathfrak{C}}\lesssim \aleph.\notag
\end{align}

 Next, we present some estimates for time regularity, which will be used  later. We again apply Lemmas \ref{C.2} and \ref{C.3}  to \eqref{4.0} and get for $N\in\mathbb{N}_0$
 
\begin{align}
\|v_q^1\|_{C^{\frac13-\frac{\kappa+\fa}2}_{[N+1,N+2]}L^\frac32}&\lesssim e^{-\fc}\|v_q^1(N)\|_{B_{\frac32,\infty}^{\frac23-\fa-\kappa}}+\fc^{-\frac{1/3-\fa}{2}}\sup_{s\in[N,N+2]}(\|z(s)\|_{C^{-\fa-\kappa}}+\|z^{:2:}(s)\|_{C^{-2\fa-\kappa}})\notag\\
&+\sup_{s\in[N,N+2]}\|v_q^1(s)+v_q^2(s)\|_{L^\frac32}\|\cH_{ R(s)}z(s)\|_{C^{-\frac13-\fa-\kappa}},\notag\\
&\lesssim \fc^{-\frac{1/3-\fa}{2}}\sup_{s\in[N,N+2]}(\|z(s)\|_{C^{-\fa-\kappa}}+\|z^{:2:}(s)\|_{C^{-2\fa-\kappa}})\notag\\
&+\sup_{s\in[N,N+2]}\|v_q^1(s)\|_{B_{\frac32,\infty}^{\frac23-\fa-\kappa}}+ \frac1{\mathfrak{C}}\sup_{s\in[N,N+2]}\|v_q^2(s)\|_{W^{\frac12,\frac65}}.\notag
\end{align}
Then by  \eqref{3.7}, \eqref{3.7*}, \eqref{3.10}, \eqref{3.13.1}, \eqref{3.13.2} and interpolation we obtain

\begin{align}
 \$ v_q^1 \$ _{C_t^{\frac1{6}-\frac{\kappa+\fa}2}B^{\frac13}_{\frac32,\infty},2r}+ \$ v_q^1 \$ _{C_t^{\frac13-\frac{\kappa+\fa}2}L^\frac32,2r}&\lesssim1,\notag\\
 \$ v_q^1 \$ _{C_t^{\frac1{6}-\frac{\kappa+\fa}2}B^{\frac13}_{\frac32,\infty},m}+ \$ v_q^1 \$ _{C_t^{\frac13-\frac{\kappa+\fa}2}L^\frac32,m}&\lesssim(7^{q-1}\cdot \frac{7}2 \bar{M}_L m)^{\frac{21}2\cdot {7}^{q-1}}.\label{4.2}
\end{align}

\subsection{Proof of  (\ref{3.12})}
Additionally, by employing the same method as previously, we can derive for $s\in[N,N+2]$
 \begin{align*}
\|v_{q+1}^1(s)&-v_q^1(s)\|_{B_{\frac32,\infty}^{\frac23-\fa-\kappa}}\lesssim   e^{-\fc(s-N)}\|v_{q+1}^1(N)-v_q^1(N)\|_{B_{\frac32,\infty}^{\frac23-\fa-\kappa}}\notag\\
&+\sup_{s\in[N,N+2]}\|(v_{q+1}^1+v_{q+1}^2-v_q^1-v_q^2)(s)\|_{L^\frac32}\|\cH_{R}z(s)\|_{C^{-\fa-\kappa-\frac13}}\notag\\
&+\sup_{s\in[N,N+2]}\|v_q^1(s)+v_q^2(s)\|_{L^\frac32}\|(\cL_{ f(q+1)}-\cL_{ f(q)})\cH_{R}z(s)\|_{C^{-\fa-\kappa-\frac13}}.
\end{align*}
Then by (\ref{deltaf}), $ 0<\fa_0\leq \frac1{12}$ and Proposition \ref{deltaz} we have
\begin{align*}
 \sup_{s\in[N,N+2]}\|\cH_{R}z(s)\|_{C^{-\fa-\kappa-\frac13}}&\lesssim1,\\
\sup_{s\in[N,N+2]}\|(\cL_{ f(q+1)}-\cL_{ f(q)})\cH_{R}z\|_{C^{-\fa-\kappa-\frac13}}&\lesssim \lambda_{q}^{-\frac{\alpha}{2}(\frac16-\fa_0)}\sup_{s\in[N,N+2]}\|\cH_{R}z\|_{C^{-\frac16-\fa-\fa_0-\kappa}}\lesssim \frac1{\mathfrak{C}}\lambda_{q}^{-\frac{\alpha}{24}}.
\end{align*}
Also a similar argument as \eqref{4.1} implies that
\begin{align*}
\mathbf{E}[\sup_{s\in[N+1,N+2]}\|v_{q+1}^1(s)&-v_q^1(s)\|_{B_{\frac32,\infty}^{\frac23-\fa-\kappa}}^{2r}]^{\frac{1}{2r}}\leq \frac12\mathbf{E}[\sup_{s\in[N,N+1]}\|v_{q+1}^1(s)-v_q^1(s)\|_{B_{\frac32,\infty}^{\frac23-\fa-\kappa}}^{2r}]^{\frac{1}{2r}}\notag\\ 
&+M_1\( \$ v_{q+1}^2-v_q^2 \$ _{W^{\frac12,\frac65},2r}+( \$ v_q^1 \$ _{L^2,2r}+ \frac1{\mathfrak{C}}\$ v_q^2 \$ _{W^{\frac12,\frac65},2r})\lambda_{q}^{-\frac\alpha{24}}\).
\end{align*}
Using the fact that $v_q^1(0)=0$  we obtain
\begin{align}
 \mathbf{E}[\sup_{s\in[0,1]}\|v_{q+1}^1(s)-v_q^1(s)\|_{B_{\frac32,\infty}^{\frac23-\fa-\kappa}}^{2r}]^{\frac{1}{2r}}\lesssim \$ v_{q+1}^2-v_q^2 \$ _{W^{\frac12,\frac65},2r}+( \$ v_q^1 \$ _{L^2,2r}+ \frac1{\mathfrak{C}}\$ v_q^2 \$ _{W^{\frac12,\frac65},2r})\lambda_{q}^{-\frac\alpha{24}}.\notag
\end{align}
 Thus together with \eqref{omegan}, \eqref{3.10}, (\ref{3.16}), \eqref{3.13.1}  and similar as \eqref{4.1} we obtain for any $N\in\mathbb{N}_0$
 \begin{align}
 \mathbf{E}[\sup_{s\in[N,N+1]}\|v_{q+1}^1(s)-v_q^1(s)\|_{B_{\frac32,\infty}^{\frac23-\fa-\kappa}}^{2r}]^{\frac{1}{2r}}&\lesssim \$ v_{q+1}^2-v_q^2 \$ _{W^{\frac12,\frac65},2r}+( \$ v_q^1 \$ _{L^2,2r}+  \frac1{\mathfrak{C}}\$ v_q^2 \$ _{W^{\frac12,\frac65},2r})\lambda_{q}^{-\frac\alpha{24}}\notag\\ 
 &\lesssim \lambda_{q+1}^{-\alpha/2}+(1+ \frac{1+\fm}{\mathfrak{C}})\lambda_{q}^{-\frac\alpha{24}}\lesssim \lambda_{q}^{-\frac\alpha{24}}.\notag
\end{align}
 Here similar as \eqref{4.1}, the implicit constants are independent of $N$. Moreover, by the embedding $B_{\frac32,\infty}^{\frac23-\fa-\kappa}\subset L^2$ we  deduce 
\begin{align}
  \$ v_{q+1}^1-v_q^1 \$ _{L^2,2r}  &\lesssim \$ v_{q+1}^1-v_q^1 \$ _{B_{\frac32,\infty}^{\frac23-\fa-\kappa},2r}\lesssim\lambda_{q}^{-\frac\alpha{24}}\leq \delta_{q+1}^{1/2},\notag
\end{align}
where we used the conditions on the parameters in Section \ref{cops} to deduce $\frac{\alpha}{24}>\beta b$, and chose $a$ large enough to absorb the universal constant. 

Similarly by \eqref{3.10}, \eqref{3.16} and  \eqref{3.13.2}  we deduce 
\begin{align}
   \$ v_{q+1}^1-v_q^1 \$ _{L^2,m}&\lesssim   \$ v_{q+1}^2-v_q^2 \$ _{W^{\frac12,\frac65},m}+( \$ v_q^1 \$ _{L^2,m}+ \$ v_q^2 \$ _{W^{\frac12,\frac65},m})\lambda_{q}^{-\frac\alpha{24}}\lesssim \lambda_{q}^{-\frac\alpha{24}}(7^{q}\cdot \frac{7}2 \bar{M}_L m)^{\frac{21}2\cdot {7}^{q}}.\label{v1l2}
\end{align}

\subsection{Proof of  (\ref{4.4.1})}

Finally we estimate the following bound. These norms may priori blow up during the iteration. By the same argument as \eqref{4.1} we obtain  for $0<\kk<\frac23-2\fa$

\begin{align}
\|v_q^1\|_{C_{[N+1.N+2]}^{\frac16}L^\infty}&+\sup_{s\in[N+1,N+2]}\|v_q^1(s)\|_{C^{\frac13}}\notag\\ 
&\lesssim e^{-\fc}\sup_{s\in[N,N+1]}\|v_q^1(s)\|_{C^{\frac13}}+\sup_{s\in[N,N+2]}(
\|z(s)\|_{C^{-\fa-\kappa}}+\|z^{:2:}(s)\|_{C^{-2\fa-\kappa}})\notag\\
&+\sup_{s\in[N,N+2]}\|v_q^1(s)+v_q^2(s)\|_{C^{-\frac43}}\|\cL_{ f(q)}\cH_{R}z(s)\|_{C^{\frac23}}.\notag
\end{align}
By (\ref{deltaf}) and the fact $\frac23+\fa+\kappa<1$ we get for $s\in[N,N+2]$
\begin{align*}
    \|\cL_{ f(q)}\cH_{R}z(s)\|_{C^{\frac23}}\lesssim \lambda_q^{\frac{\alpha}{2}(\frac23+\fa+\kappa)}\|z(s)\|_{C^{-\fa-\kappa}}\lesssim \lambda_q^{\alpha/2}\|z(s)\|_{C^{-\fa-\kappa}}.
\end{align*}
Then by the embedding $L^2\subset C^{-\frac43},W^{\frac12,\frac65}\subset C^{-\frac43}$, (\ref{3.7}), \eqref{3.7*}, (\ref{3.10}) and (\ref{3.13.2}) we obtain 
\begin{align}
  &  \$ v_q^1 \$ _{C_t^{\frac16}L^\infty,m}+  \$ v_q^1 \$ _{C^{\frac13},m}\notag\\
   &\lesssim  \$ z \$ _{C^{-\fa-\kappa},m}+ \$ z^{:2:} \$ _{C^{-2\fa-\kappa},m}+\lambda_q^{\alpha/2}( \$ v_q^1 \$ _{L^{2},\frac{8m}7}+ \$ v_q^2 \$ _{W^{\frac12,\frac65},\frac{8m}7}) \$ z \$ _{C^{-\fa-\kappa},8m}\notag\\
    &\lesssim \lambda_q^{\alpha/2}(7^{q-1}\cdot 4 \bar{M}_L m)^{\frac{21}2\cdot {7}^{q-1}}m^{1/2}L,\label{4.4.2}
    \end{align}
and
    \begin{align}
  \$ v_q^1 \$ _{C_t^{\frac16}L^\infty,2r}+   \$ v_q^1 \$ _{C^{\frac13},2r}&\lesssim \lambda_q^{\alpha/2}(7^{q-1}\cdot8 \bar{M}_L r)^{\frac{21}2\cdot {7}^{q-1}}r^{1/2}L\lesssim\lambda_q^{2\alpha}r^{1/2}L,\notag
\end{align}
where in the last inequality we used the conditions on the parameters $a>(56 \bar{M}_L r)^d,\alpha b>7$ in Section \ref{cops}  to deduce for $q\geq1$, $(7^{q-1}\cdot8 \bar{M}_L r)^{\frac{21}2\cdot {7}^{q-1}}\leq \lambda_{q-1}^{\frac{21}2}\leq \lambda_q^{\frac32\alpha}$. In the case $q=0$ we used $ \bar{M}_L r\leq a^\alpha$.

\section{Detailed iteration – proof of Proposition \ref{prop:3.1}}\label{proofof3.1}
The proof follows a similar structure as many convex integration schemes. First, we begin by setting the parameters in Section~\ref{cops} and then performing a mollification step in Section~\ref{moll}. In Section~\ref{covq+12}, we introduce the new iteration $v_{q+1}^2$. This is the core of the construction, utilizing the accelerating jets introduced in Appendix \ref{s:appB}. Section~\ref{guina4} contains the inductive estimates for $v_{q+1}^2$.
 Finally, in Section~\ref{dotrs}, we define the new stress $\mathring{R}_{q+1}$  and establish the inductive moment estimate for $\mathring{R}_{q+1}$ in Section~\ref{votiefr}.

\subsection{Choice of parameters}\label{cops}
In the subsequent analysis, additional parameters will be indispensable and their value has to be carefully chosen in order to respect all the compatibility conditions appearing in the estimates below.  For a sufficiently small $\alpha\in(0,1)$ to be chosen below, we take $l:=\lambda_{q+1}^{-\frac{3\alpha}{2}}\lambda_q^{-2}$ and have $l\lambda_q^4\leq\lambda_{q+1}^{-\alpha},l^{-1}\leq\lambda_{q+1}^{2\alpha}$ provided $\alpha b>4$. 

We also introduce a small constant $\gamma=\frac{1}{113}.$
In the sequel, we use the following bounds
$$\alpha<\frac\gamma{54},\ \ \alpha b>\frac{8}{\fa_0},\ \ \alpha>( 48\beta b^2)\vee(\frac{8}{\fa_0}\beta b).$$
Here we recall $\fa_0:=\frac{1-3\fa}6\wedge \frac1{12}$. In addition, we postulate
\begin{align}
 \lambda_{q}^{-\frac{\alpha}{24}+2\beta b^2}\fc\ll1,\ \ \lambda_{q+1}^{-\frac{\fa_0}{4}\alpha+2\beta b}\fC^6\ll 1,\label{a2}
\end{align}
which can be satisfied by choosing $a$ large enough.

The above can be obtained by choosing $\alpha=\frac{1}{55\times 113}$, and choosing $b\in 113\mathbb{N}$ large enough such that
$b> \frac{8}{\fa_0\alpha}$, and finally choosing $\beta$ small such that $\alpha>(48\beta b^2)\vee( \frac{8}{\fa_0}\beta b)$. The last free parameter is $a$ which satisfies
the lower bounds given through (\ref{ab2}) and (\ref{a2}). Let $d$ satisfies $q{7}^q\leq d8^q$ for any $q\in \N$. We then
choose $a\in \mathbb{N}$ which obeys $a>(70 \bar{M}_L r)^d\vee( \bar{M}_L r)^{\frac1{\alpha}}$.   In the sequel, we increase $a$ in order to absorb various implicit
and universal constants. 

\subsection{Mollification}\label{moll}
We intend to replace $v^2_q$ by a mollified velocity field $v_l^2$, and we define
\begin{align}
v_l^1=(v_q^1*_x\phi_l)*_t\varphi_l,\ \ v_l^2=(v_q^2*_x\phi_l)*_t\varphi_l,\ \ \mathring{R}_l=(\mathring{R}_q*_x\phi_l)*_t\varphi_l\notag,
\end{align}
where $\phi_l:=\frac{1}{l^2}\phi(\frac{\cdot}{l})$ is a family of standard mollifiers on $\mathbb{R}^2$, and $\varphi_l:=\frac{1}{l}\varphi(\frac{\cdot}{l})$ is a family of standard mollifiers with support on $(0,1)$. Here the one side mollifier is used to preserve adaptedness, i.e. $v_l^1,v_l^2$ and $\mathring{R}_l$ are $(\mathcal{F}_t)_{ t\in\mathbb{R}}$-adapted. Then by (\ref{3.6}) we obtain
\begin{align}
(\partial_t-\Delta)v_l^2-\fc v_l^1+\nabla p_l+\div N_{com}&=\div(\mathring{R}_l+{R}_{com}),\\
\div v_l^2&=0,\notag
\end{align}\label{4.5}
where
\begin{align*}
N_{com}=(v_q^1&+v_q^2)\otimes\cL_{R}z+\cL_{R}z\otimes(v_q^1+v_q^2)+(v_q^1+v_q^2)\otimes(v_q^1+v_q^2)\\
&+(v_q^1+v_q^2)\succcurlyeq\!\!\!\!\!\!\!\bigcirc\cH_{R}z+\cH_{R}z\preccurlyeq\!\!\!\!\!\!\!\bigcirc(v_q^1+v_q^2),    
\end{align*}
and $R_{com}$ is the trace-free part of 
${N}_{com}-{N}_{com}*_x\phi_l*_t\varphi_l,$
$$p_l:=(p_q^2*_x\phi_l)*_t\varphi_l-\frac{1}{2}\tr(N_{com}-N_{com}*_x\phi_l*_t\varphi_l).$$
We have for any $t\in\mathbb{R}$
\begin{align*}
    \|v_q^2(t)-v_l^2(t)\|_{L^2}\lesssim l\|v_q^2\|_{C_{[t-l,t],x}^1},
\end{align*}
which by (\ref{3.9}) implies that
\begin{align}
     \$ v_q^2-v_l^2 \$ _{L^2,2r}\lesssim l \$ v_q^2 \$ _{{C_{t,x}^1},2r}\leq l\fm\lambda_q^4\leq \fm\lambda_{q+1}^{-\alpha}\leq\frac{1}{4}\fm\delta_{q+1}^{1/2},\label{vq-vl}
\end{align}
where we used the fact that $l\lambda_q^4\leq\lambda_{q+1}^{-\alpha},\alpha>\beta$. In addition, for any $t\in\mathbb{R}$
\begin{align}
\|v_l^2\|_{C_{[t,t+1]}L^2}\leq \sup_{s\in[t-l,t+1]}\|v_q^2(s)\|_{L^2},\ \ \|v_l^2\|_{L^2_{[t,t+1]}L^2}\leq \|v_q^2\|_{L^2_{[t-l,t+1]}L^2}.\label{vl<vq}
\end{align}

\subsection{Construction of $v_{q+1}^2$}\label{covq+12}
Let us now move forward with the construction of the perturbation $w_{q+1}$, which define the next iteration by $v_{q+1}^2 := v_l^2 +w_{q+1}$. For this purpose, we make use of the accelerating jet flows introduced in \cite[Section 3]{CL20}, which are briefly recalled in Appendix \ref{s:appB}. In particular, the building blocks $W_{(\xi)}$ for $\xi\in\Lambda$ are defined in (\ref{Wxi}), and the set $\Lambda$ is introduced in Lemma \ref{lem:B.1}. The necessary estimates are collected in (\ref{2.4}). For the accelerating jet flows we choose the following parameters
$$\sigma=\lambda_{q+1}^{2\gamma},\ \eta=\lambda_{q+1}^{16\gamma},\ \nu=\lambda_{q+1}^{1-8\gamma},\ \mu=\lambda_{q+1},\ \theta=\lambda_{q+1}^{1+5\gamma},$$
Here  we recall $\gamma=\frac1{113}$ and  it is required that $b\in113\mathbb{N}$ to ensure the above parameters are integers.

 As the next step, we shall define certain amplitude functions used in the definition of the perturbations $w_{q+1}$ similarly as \cite[Section 3.1.3]{HZZ22}. Let
$$ \rho:=\sqrt{l^2+|\mathring{R}_l|^2}.$$
Then it follows that $\rho$ is $(\mathcal{F}_ t)_{t\in\mathbb{R}}$-adapted. Then by the estimate of $\rho$ in Appendix \ref{s:appA} we have for $N\geq1,t\in\R$
\begin{align}
\|\rho\|_{C_{{[t,t+1]},x}^N}\lesssim l^{-N-4}\|\mathring{R}_q\|_{L^1_{[t-l,t+1]}L^1}+l^{-6N+1}\|\mathring{R}_q\|_{L^1_{[t-l,t+1]}L^1}^N+1.
\label{4.6}
\end{align}

The amplitude functions are defined as follows:
\begin{align}
a_{(\xi)}(\omega,x,t):=\rho^{1/2}(\omega,x,t)\gamma_{\xi}(\Id-\frac{\mathring{R}_l(\omega,x,t)}{\rho(\omega,x,t)}), \label{4.7}
\end{align}
where $\gamma_\xi$ is introduced in Lemma \ref{lem:B.1}. By $|\mathrm{Id}-\frac{\mathring{R}_l}{\rho}-\Id|\leq\frac{1}{2}$ and Lemma \ref{lem:B.1},  it holds 
\begin{align}
\rho\Id-\mathring{R}_l=\sum_{\xi\in\Lambda}\rho\gamma_{\xi}^2(\Id-\frac{\mathring{R}_l}{\rho})\xi\otimes\xi=\sum_{\xi\in\Lambda}a_{(\xi)}^2\xi\otimes\xi.\label{phoid-r}
\end{align}
Then by the estimate of $a_{(\xi)}$ in Appendix \ref{s:appA} we have for $N\geq0, t\in\mathbb{R}$

\begin{align}
\|a_{(\xi)}\|_{C_{{[t,t+1]},x}^N}&\lesssim l^{-7-6N}(\|\mathring{R}_q\|_{L^1_{[t-l,t+1]}L^1}+1)^{N+1}.\label{4.8}
\end{align}

With these preparations in hand, similarly as \cite[Section 3.4]{CL20} we define the principal part
\begin{align}\label{4.9}
w_{q+1}^{(p)}:=\sum_{\xi\in\Lambda}a_{(\xi)}g_{(\xi)}W_{(\xi)},
\end{align}
where $W_{(\xi)},g_{(\xi)}$ are introduced in Appendix \ref{s:appB}. Then by (\ref{wpingjun}), (\ref{wbujiao}) and (\ref{phoid-r}) we obtain
\begin{align}\label{4.10}
w_{q+1}^{(p)}\otimes w_{q+1}^{(p)}+\mathring{R}_l&=\sum_{\xi\in\Lambda}a_{(\xi)}^2g_{(\xi)}^2
(W_{(\xi)}\otimes W_{(\xi)})-\sum_{\xi\in\Lambda}a_{(\xi)}^2\xi\otimes\xi+\rho \Id\notag\\
&=\sum_{\xi\in\Lambda}a_{(\xi)}^2g_{(\xi)}^2\mathbb{P}_{\neq0}(W_{(\xi)}\otimes W_{(\xi)})+\sum_{\xi\in\Lambda}a_{(\xi)}^2(g_{(\xi)}^2-1)\xi\otimes\xi+\rho \Id,
\end{align}
where we recall the notation $\mathbb{P}_{\neq0}f:=f-\int_{\mathbb{T}^2}\!\!\!\!\!\!\!\!\!\!\!\!\; {}-{}\  f\dif x$.

We define the incompressibility corrector by
\begin{align}
w_{q+1}^{(c)}:=\sum_{\xi\in\Lambda}\left(a_{(\xi)}g_{(\xi)}W_{(\xi)}^{(c)}+\sigma^{-1}\nabla^{\perp}a_{(\xi)}g_{(\xi)}\Psi_{(\xi)}\right),\label{4.11}
\end{align}
where $W_{(\xi)}^{(c)}$ and $\Psi_{(\xi)}$ are introduced in Appendix \ref{s:appB}. By (\ref{2.6}) we have
\begin{align}\label{wp+wc}
  w_{q+1}^{(p)}+w_{q+1}^{(c)}=\sigma^{-1}\sum_{\xi\in\Lambda}\nabla^{\perp}[a_{(\xi)}g_{(\xi)}\Psi_{(\xi)}],
\end{align}
which implies $\div(w_{q+1}^{(p)}+w_{q+1}^{(c)})=0$.

The temporal corrector $w_{q+1}^{(t)}$ is defined as
\begin{align}
    w_{q+1}^{(t)} = w_{q+1}^{(o)} + w_{q+1}^{(a)},\label{t=0+a}
\end{align}
where $w_{q+1}^{(o)}$ is the temporal oscillation corrector
\begin{align}\label{4.12}
w_{q+1}^{(o)}:=-\sigma^{-1}\mathbb{P}_{\rm{H}}\mathbb{P}_{\neq0}\sum_{\xi\in\Lambda}h_{(\xi)}\div(a_{(\xi)}^2\xi\otimes\xi),
\end{align}
 and $w_{q+1}^{(a)}$ is the acceleration corrector
\begin{align}\label{4.13}
w_{q+1}^{(a)}:=-\theta^{-1}\sigma\mathbb{P}_{\rm{H}}\mathbb{P}_{\neq0}\sum_{\xi\in\Lambda}a_{(\xi)}^2g_{(\xi)}|W_{(\xi)}|^2\xi.
\end{align}
Here $\mathbb{P}_{\rm{H}}$ is the Helmholtz projection, and $h_{(\xi)}, g_{(\xi)}$ are given in Appendix \ref{s:appB}.

Since $\rho$ and $\mathring{R}_l$ are $(\mathcal{F}_t)_{t\in\mathbb{R}}$-adapted, we know $a_{(\xi)}$ is $(\mathcal{F}_t)_{ t\in\mathbb{R}}$-adapted. Moreover, as $W_{(\xi)},W_{(\xi)}^{(c)}$, $\Psi_{(\xi)}, g_{(\xi)}$ are deterministic, $w_{q+1}^{(p)},w_{q+1}^{(c)},w_{q+1}^{(t)}$ are also $(\mathcal{F}_t)_{ t\in\mathbb{R}}$-adapted.

Finally, the total perturbation $w_{q+1}$ and the new velocity $v_{q+1}^2$ are defined by
$$w_{q+1}:={w}_{q+1}^{(p)}+{w}_{q+1}^{(c)}+w_{q+1}^{(t)},\ \ v_{q+1}^2:=v_l^2+w_{q+1},$$
which are $(\mathcal{F}_t)_{ t\in\mathbb{R}}$-adapted and divergence-free with mean-zero.

\subsection{Verification of the inductive estimates for $v_{q+1}^2$}\label{guina4}
First we estimate $L^n$-norm for $n\in(1,\infty)$ and then verify the inductive estimates (\ref{3.15})-(\ref{3.16}) and (\ref{3.8.1})-(\ref{3.10}) at the level $q+1$.

We first consider bound $w^{(p)}_{q+1}$ in $L^2$. We know that $W_{(\xi)}$ is $(\mathbb{T}/\sigma)^2$-periodic, so by applying Theorem \ref{ihiot} with
$p=2$, (\ref{4.7}), (\ref{4.8}), (\ref{4.9}) and (\ref{2.4}) we obtain for $s\in[t,t+1]$, $t\in\mathbb{R}$
\begin{align}
\|{w}_{q+1}^{(p)}(s)\|_{L^2}&\lesssim\sum_{\xi\in\Lambda}|g_{(\xi)}(s)|(\|a_{(\xi)}(s)\|_{L^2}\|W_{(\xi)}\|_{L^2}+\sigma^{-1/2}\|a_{(\xi)}(s)\|_{C_x^1}\|W_{(\xi)}\|_{L^2})\notag\\
&\lesssim\sum_{\xi\in\Lambda}|g_{(\xi)}(s)|(\|a_{(\xi)}(s)\|_{L^2}+\sigma^{-1/2}l^{-13}(\|\mathring{R}_q\|_{L^1_{[t-l,t+1]}L^1}+1)^{2}).\notag\\
&\lesssim \sum_{\xi\in\Lambda}|g_{(\xi)}(s)|(\|\rho(s)\|_{L^1}^{1/2}+\lambda_{q+1}^{26\alpha-\gamma}(\|\mathring{R}_q\|_{L^1_{[t-l,t+1]}L^1}^2+1)).\notag
\end{align}
Then we obverse that $g_{(\xi)}$ is $\mathbb{T}/\sigma$-periodic, we again apply Theorem \ref{ihiot} and   (\ref{gwnp}) to obtain  
\begin{align}
\|{w}_{q+1}^{(p)}\|_{L^2_{[t,t+1]}L^2}&\lesssim\sum_{\xi\in\Lambda}\|g_{(\xi)}\|_{L^2_{[t,t+1]}}(\|\rho\|_{L^1_{[t,t+1]}L^1}^{1/2}+\sigma^{-1/2}\left\|\|\rho(\cdot)\|_{L^1}^{1/2}\right\|_{C^{0,1}{[t,t+1]}})\notag\\
&\ \ \ \ \ \ \ +\lambda_{q+1}^{26\alpha-\gamma}(\|\mathring{R}_q\|_{L^1_{[t-l,t+1]}L^1}^2+1)\notag\\
&\lesssim(\|\mathring{R}_q\|_{L^1_{[t-l,t+1]}L^1}+l)^{1/2}+ \lambda_{q+1}^{11\alpha-\gamma}(\|\mathring{R}_q\|_{L^1_{[t-l,t+1]}L^1}+1)\notag\\
&\ \ \ \ \ \ \ +\lambda_{q+1}^{26\alpha-\gamma}(\|\mathring{R}_q\|_{L^1_{[t-l,t+1]}L^1}^2+1)\notag\\
&\leq\frac{M_0}{4}(\|\mathring{R}_q\|_{L^1_{[t-l,t+1]}L^1}^{1/2}+ \delta_{q+1}^{1/2})+ M_0\lambda_{q+1}^{26\alpha-\gamma}(\|\mathring{R}_q\|_{L^1_{[t-l,t+1]}L^1}^2+1),\label{wp21}
\end{align}
where we used $l\leq \delta_{q+1}$ and $M_0$ is an universal constant. In the  second inequality we used

\begin{align*}
    \|\rho\|^{1/2}_{L^1_{[t,t+1]}L^1}&\lesssim (l+\|\mathring{R}_{q}\|_{L^1_{[t-l,t+1]}L^1})^{\frac12},
\end{align*}
and by \eqref{4.6} 
$$\left\|\|\rho(\cdot)\|_{L^1}^{1/2}\right\|_{C^{0,1}{([t,t+1])}}\lesssim l^{-1/2}\|\rho\|_{C_{[t,t+1],x}^1}\lesssim l^{-11/2}(\|\mathring{R}_q\|_{L^1_{[t-l,t+1]}L^1}+1),$$
where $C^{0,1}$ is the Lipschitz norm.

For general $L^n$-norm with $n\in(1,\infty)$, by (\ref{4.8}), (\ref{4.9}) and (\ref{2.4}) we have for $s\in[t,t+1]$, $t\in\mathbb{R}$
\begin{align}
\|{w}_{q+1}^{(p)}(s)\|_{L^n}&\lesssim \sum_{\xi\in\Lambda}\|a_{(\xi)}(s)\|_{L^\infty}|g_{(\xi)}(s)|\|W_{(\xi)}(s)\|_{L^n}\notag\\
&\lesssim (\|\mathring{R}_q\|_{L^1_{[t-l,t+1]}L^1}+1)l^{-7}\sum_{\xi\in\Lambda}|g_{(\xi)}(s)|(\nu\mu)^{1/2-1/n}\notag\\
&\lesssim(\|\mathring{R}_q\|_{L^1_{[t-l,t+1]}L^1}+1)\lambda_{q+1}^{14\alpha+(1/2-1/n)(2-8\gamma)}\sum_{\xi\in\Lambda}|g_{(\xi)}(s)|.\label{4.16}
\end{align}
By (\ref{4.8}), (\ref{4.11}) and (\ref{2.4}) we obtain
\begin{align}
\|{w}_{q+1}^{(c)}(s)\|_{L^n}&\lesssim\sum_{\xi\in\Lambda}|g_{(\xi)}(s)|\left(\|a_{(\xi)}(s)\|_{L^\infty}\|W^{(c)}_{(\xi)}(s)\|_{L^n}+\sigma^{-1}\|\nabla^{\perp} a_{(\xi)}(s)\|_{L^\infty}\|\Psi_{(\xi)}(s)\|_{L^n}\right)\notag\\
&\lesssim(\|\mathring{R}_q\|_{L^1_{[t-l,t+1]}L^1}+1)^2l^{-13}\sum_{\xi\in\Lambda}|g_{(\xi)}(s)|(\nu\mu^{-1}+\sigma^{-1}\mu^{-1})(\nu\mu)^{1/2-1/n}\notag\\
&\lesssim(\|\mathring{R}_q\|_{L^1_{[t-l,t+1]}L^1}+1)^2\lambda_{q+1}^{26\alpha+(1/2-1/n)(2-8\gamma)-8\gamma}\sum_{\xi\in\Lambda}|g_{(\xi)}(s)|.\label{4.17}
\end{align}
By (\ref{4.8}), (\ref{4.12}), (\ref{hlinfty}) and (\ref{2.4}) we get
\begin{align}
\|{w}_{q+1}^{(o)}(s)\|_{L^n}&\lesssim \sigma^{-1}\sum_{\xi\in\Lambda}\|h_{(\xi)}(s)\|_{L^\infty}\|\div(a_{(\xi)}^2(s)\xi\otimes\xi)\|_{L^n}\notag\\
&\lesssim (\|\mathring{R}_q\|_{L^1_{[t-l,t+1]}L^1}+1)^3l^{-20}\sigma^{-1}\lesssim (\|\mathring{R}_q\|_{L^1_{[t-l,t+1]}L^1}+1)^3\lambda_{q+1}^{40\alpha-2\gamma},\label{4.18}
\end{align}
and by (\ref{4.8}), (\ref{4.13}) and (\ref{2.4}) we obtain
\begin{align}
\|{w}_{q+1}^{(a)}(s)\|_{L^n}&\lesssim\sigma\theta^{-1}\sum_{\xi\in\Lambda}|g_{(\xi)}(s)|\|a_{(\xi)}(s)\|_{L^\infty}^2\|W_{(\xi)}(s)\|_{L^{2n}}^2\notag\\
&\lesssim (\|\mathring{R}_q\|_{L^1_{[t-l,t+1]}L^1}+1)^2l^{-14}\sigma\theta^{-1}(\nu\mu)^{1-1/n}\sum_{\xi\in\Lambda}|g_{(\xi)}(s)|\notag\\
&\lesssim(\|\mathring{R}_q\|_{L^1_{[t-l,t+1]}L^1}+1)^2\lambda_{q+1}^{28\alpha+(1-1/n)(2-8\gamma)-3\gamma-1}\sum_{\xi\in\Lambda}|g_{(\xi)}(s)|.\label{4.19}
\end{align}
Thus we have
\begin{align*}
    \|{w}_{q+1}^{(c)}(s)+{w}_{q+1}^{(t)}(s)\|_{L^2}
    &\lesssim \lambda_{q+1}^{40\alpha-2\gamma}(\|\mathring{R}_q\|_{L^1_{[t-l,t+1]}L^1}+1)^3(\sum_{\xi\in\Lambda}|g_{(\xi)}(s)|+1).
\end{align*}
Integrating in time and  (\ref{gwnp}) imply
\begin{align}
\|w^{(c)}_{q+1}+w^{(t)}_{q+1}\|_{L^2_{[t,t+1]}L^2}\lesssim\lambda_{q+1}^{40\alpha-2\gamma}(\|\mathring{R}_q\|_{L^1_{[t-l,t+1]}L^1}+1)^3.\label{wc+wt22}
\end{align}
Then by the above estimate and (\ref{wp21}) we obtain
\begin{align}
\|w_{q+1}\|_{L^2_{[t,t+1]}L^2}
\leq \frac{M_0}{4}(\|\mathring{R}_q\|_{L^1_{[t-l,t+1]}L^1}^{1/2}+\delta_{q+1}^{1/2})+M_0\lambda_{q+1}^{40\alpha-\gamma}(\|\mathring{R}_q\|_{L^1_{[t-l,t+1]}L^1}^3+1),\notag
\end{align}
which together with (\ref{3.11}) implies

\begin{align}
 \$ w_{q+1} \$ _{L^2_tL^2,2r}
&\leq \frac{M_0}{4}( \$ \mathring{R}_q \$ _{L^1_tL^1,r}^{1/2}+\delta_{q+1}^{1/2})+M_0\lambda_{q+1}^{40\alpha-\gamma}(8 \$ \mathring{R}_q \$ _{L^1_tL^1,6r}^3+1)\notag\\
&\leq \frac{M_0}4\delta_{q+1}^{1/2} (M_L^{1/2}+1)+ 9M_0\lambda_{q+1}^{40\alpha-\gamma}(7^q\cdot6\bar{M}_Lr)^{9\cdot {7}^q}\notag\\
&\leq \frac{M_0}4\delta_{q+1}^{1/2}(M_L^{1/2}+1)+9M_0\lambda_{q+1}^{40\alpha-\gamma}\lambda_q^9\leq \frac{3M_0}4\delta_{q+1}^{1/2}M_L^{1/2},\label{4.20.1}
\end{align}
where we used the conditions on the parameters $a>(42\bar{M}_Lr)^d$ to deduce $(7^q\cdot6\bar{M}_Lr)^{{7}^q}\leq\lambda_q$ and used $40\alpha-\gamma+9/b<-\beta$ in the last inequality. Similarly
\begin{align}
 \$ w_{q+1} \$ _{L^2_tL^2,m}
&\leq \frac{M_0}4( \$ \mathring{R}_q \$ _{L^1_tL^1,\frac m2}^{1/2}+\delta_{q+1}^{1/2})+ M_0\lambda_{q+1}^{40\alpha-\gamma}( 8\$ \mathring{R}_q \$ _{L^1_tL^1,3m}^3+1)\notag\\
&\leq \frac{M_0}{2}(7^q\cdot3\bar{M}_Lm)^{9\cdot {7}^q}.\label{4.20.2}
\end{align}

With these bounds, we have all in hand to complete the estimate of $v_{q+1}^2$. We split the details into several subsections.
\subsubsection{Proof of (\ref{3.15}), (\ref{3.8.1}) and (\ref{3.8.2}) at the level $q+1$} By (\ref{vq-vl}), (\ref{4.20.1})  and $\fm\leq M_L^{1/2}$
 we obtain
\begin{align}
 \$ v_{q+1}^2-v_q^2  \$ _{L^2_tL^2,2r}
\leq \frac{3M_0}{4} \delta_{q+1}^{1/2}M_L^{1/2}+\frac14\delta_{q+1}^{1/2}\fm\leq  M_0\delta_{q+1}^{1/2}M_L^{1/2}.\notag
\end{align}
Hence  (\ref{3.15}) and (\ref{3.8.1}) at the level $q+1$ follow. 
Similarly by (\ref{vl<vq}) and (\ref{4.20.2}) 
\begin{align*}
 \$ v_{q+1}^2 \$ _{L^2_tL^2,m}&\leq  \$ w_{q+1} \$ _{L^2_tL^2,m}+ \$ v_l^2 \$ _{L^2_tL^2,m}\leq M_0(7^{q}\cdot 3 \bar{M}_L m)^{9\cdot {7}^q}.
\end{align*}
Hence (\ref{3.8.2}) holds at the level $q+1$.

\subsubsection{Proof of (\ref{3.9}) at the level $q+1$}
As the next step, by  (\ref{4.8}), (\ref{wp+wc}) and (\ref{gwnp}) we have for $t\in\R$
\begin{align}
\|w_{q+1}^{(p)}+w_{q+1}^{(c)}\|_{C_{[t,t+1],x}^1}&=\|\sigma^{-1}\sum_{\xi\in\Lambda}\nabla^{\perp}[a_{(\xi)}g_{(\xi)}\Psi_{(\xi)}]\|_{C_{[t,t+1],x}^1}\notag\\
&\lesssim\sigma^{-1} \sum_{\xi\in\Lambda}\|a_{(\xi)}\|_{C_{[t,t+1],x}^2}\|\nabla^{\perp}[g_{(\xi)}\Psi_{(\xi)}]\|_{C_{[t,t+1],x}^1}\notag\\
&\lesssim (\|\mathring{R}_q\|_{L^1_{[t-l,t+1]}L^1}+1)^3l^{-19}(\eta^{\frac32}\sigma+\sigma\mu\eta^{1/2}+\theta\eta\mu)(\nu\mu)^{1/2}\notag\\
&\lesssim(\|\mathring{R}_q\|_{L^1_{[t-l,t+1]}L^1}+1)^3\lambda_{q+1}^{38\alpha+3+17\gamma},\label{4.21}
\end{align}
where we used (\ref{2.4}) and (\ref{B.10}) to bound $\Psi_{(\xi)}$. 

By (\ref{4.8}), (\ref{gwnp}) and (\ref{parth}) we have
\begin{align}
\|w_{q+1}^{(o)}\|_{C_{[t,t+1],x}^1}&\lesssim \sigma^{-1}\sum_{\xi\in\Lambda}(\|\partial_th_{(\xi)}\|_{L^\infty}+\|h_{(\xi)}\|_{L^\infty})\notag\\
&\ \ \ \ \ \ \ \ \ \ \ \ \ \ \times\Big{(}\|\div(a_{(\xi)}^2\xi\otimes\xi)\|_{C_{[t,t+1]}W^{1+\alpha,\infty}}+\|\div(a_{(\xi)}^2\xi\otimes\xi)\|_{C_{[t,t+1]}^1W^{\alpha,\infty}}\Big{)}\notag\\
&\lesssim (\|\mathring{R}_q\|_{L^1_{[t-l,t+1]}L^1}+1)^5l^{-27}\eta\lesssim(\|\mathring{R}_q\|_{L^1_{[t-l,t+1]}L^1}+1)^5\lambda_{q+1}^{54\alpha+16\gamma}.\label{4.23}
\end{align} 
Here we have a extra $\alpha$ since $\mathbb{P}_{\rm{H}}\mathbb{P}_{\neq0}$ is not a bounded operator on $C^0$. Then we used (\ref{4.8}), (\ref{gwnp}), (\ref{2.4}) and (\ref{B.9}) to deduce
\begin{align}
\|w_{q+1}^{(a)}\|_{C_{[t,t+1],x}^1}&\lesssim\sigma\theta^{-1}\sum_{\xi\in\Lambda}(\|\partial_tg_{(\xi)}\|_{L^\infty}+\|g_{(\xi)}\|_{L^\infty})\notag\\
&\ \ \ \ \ \ \ \ \ \ \ \ \ \ \ \times\left(\|a_{(\xi)}^2|W_{(\xi)}|^2\xi\|_{C_{[t,t+1]}W^{1+\alpha,\infty}}+\|a_{(\xi)}^2|W_{(\xi)}|^2\xi\|_{C_{[t,t+1]}^1W^{\alpha,\infty}}\right)\notag\\
&\lesssim (\|\mathring{R}_q\|_{L^1_{[t-l,t+1]}L^1}+1)^4l^{-21}\sigma^2\theta^{-1}\eta^{\frac32}(\sigma\nu\mu^2+\nu\mu^2\theta\eta^{1/2})(\sigma\mu)^\alpha\notag\\
&\lesssim(\|\mathring{R}_q\|_{L^1_{[t-l,t+1]}L^1}+1)^4\lambda_{q+1}^{42\alpha+3+28\gamma}.\label{4.24}
\end{align}
Here we have a extra $\alpha$ since $\mathbb{P}_{\rm{H}}\mathbb{P}_{\neq0}$ is not a bounded operator on $C^0$. Thus combining (\ref{3.11}), (\ref{4.21})-(\ref{4.24}) we obtain the first inequality in (\ref{3.9}) at the level $q+1$, i.e.
\begin{align}
 \$ v_{q+1}^2 \$ _{C_{t,x}^1,m}&\leq  \$ v_{l}^2 \$ _{C_{t,x}^1,m}+ \$ w_{q+1} \$ _{C_{t,x}^1,m}\notag\\
&\lesssim ( \$ \mathring{R}_q \$ _{L^1_{t}L^1,5m}^5+1)\lambda_{q+1}^{54\alpha+3+28\gamma}+ \$ v_q^2 \$ _{C^1_{t,x},m}\leq \lambda_{q+1}^{7/2}(7^q\cdot 5 \bar{M}_L m)^{15\cdot {7}^q},
\end{align}
where we used the conditions $\alpha<\frac{1}{216},\gamma<\frac{1}{112}$. In particular by the choice of parameters  $a\geq (70 \bar{M}_L r)^d,\alpha b\geq30$ we deduce 
\begin{align}
 \$ v_{q+1}^2 \$ _{C_{t,x}^1,2r}\leq \lambda_{q+1}^{7/2}(7^q\cdot 10 \bar{M}_L r)^{15\cdot {7}^q}\leq \lambda_{q+1}^{7/2}\lambda_q^{15}\leq\fm\lambda_{q+1}^4,\notag
\end{align}
 which implies the second inequality in (\ref{3.9}) at the level $q+1$.

\subsubsection{Proof of (\ref{3.16}) and \eqref{3.10} at the level $q+1$}
We conclude this part with further $W^{1,n}$-norm for $n\in(1,\infty)$. By (\ref{4.8}), (\ref{wp+wc})  and (\ref{2.4}) we obtain for any $s\in[t,t+1], t\in\mathbb{R}$
\begin{align}
\|w_{q+1}^{(p)}(s)+w_{q+1}^{(c)}(s)\|_{W^{1,n}}&\lesssim \sigma^{-1}\sum_{\xi\in\Lambda}\|\nabla^{\perp}[a_{(\xi)}g_{(\xi)}\Psi_{(\xi)}](s)\|_{W^{1,n}}\notag\\
&\lesssim (\|\mathring{R}_q\|_{L^1_{[t-l,t+1]}L^1}+1)^3l^{-19}\sigma\mu(\nu\mu)^{1/2-1/n}\sum_{\xi\in\Lambda}|g_{(\xi)}(s)|
,\notag\\
&\lesssim(\|\mathring{R}_q\|_{L^1_{[t-l,t+1]}L^1}+1)^3\lambda_{q+1}^{38\alpha+(1/2-1/n)(2-8\gamma)+1+2\gamma}\sum_{\xi\in\Lambda}|g_{(\xi)}(s)|,\label{4.25}
\end{align}

\begin{align}
\|{w}_{q+1}^{(o)}(s)\|_{W^{1,n}}&\lesssim \sigma^{-1}\sum_{\xi\in\Lambda}\|h_{(\xi)}\|_{L^\infty}\|\div(a_{(\xi)}^2\xi\otimes\xi)\|_{C_{[t,t+1]}W^{1,n}}\notag\\
&\lesssim (\|\mathring{R}_q\|_{L^1_{[t-l,t+1]}L^1}+1)^4\sigma^{-1}l^{-26}\lesssim(\|\mathring{R}_q\|_{L^1_{[t-l,t+1]}L^1}+1)^4\lambda_{q+1}^{52\alpha-2\gamma},\label{4.27}
\end{align}
and
\begin{align}
\|{w}_{q+1}^{(a)}(s)\|_{W^{1,n}}&\lesssim\sigma\theta^{-1}\sum_{\xi\in\Lambda}|g_{(\xi)}(s)|\|a_{(\xi)}^2\|_{C_{[t,t+1],x}^1}\|W_{(\xi)}\|_{C_{[t,t+1]}W^{1,2n}}\|W_{(\xi)}\|_{C_{[t,t+1]}L^{2n}}\notag\\
&\lesssim (\|\mathring{R}_q\|_{L^1_{[t-l,t+1]}L^1}+1)^3l^{-20}\sigma\theta^{-1}\sigma\mu(\nu\mu)^{1-1/n}\sum_{\xi\in\Lambda}|g_{(\xi)}(s)|\notag\\
&\lesssim(\|\mathring{R}_q\|_{L^1_{[t-l,t+1]}L^1}+1)^3\lambda_{q+1}^{40\alpha+(1-1/n)(2-8\gamma)-\gamma}\sum_{\xi\in\Lambda}|g_{(\xi)}(s)|.\label{4.28}
\end{align}

Moreover together with (\ref{gwnp}) we have
\begin{align*}
    \|{w}_{q+1}-{w}_{q+1}^{(o)}\|_{C_{[t,t+1]}W^{1,\frac65}}&\lesssim(\|\mathring{R}_q\|_{L^1_{[t-l,t+1]}L^1}+1)^3\lambda_{q+1}^{40\alpha+\frac{1}{3}+\frac{38}{3}\gamma}.
\end{align*}
By  (\ref{4.16}), \eqref{4.17}, (\ref{4.19}) and \eqref{gwnp} we have for $t\in\mathbb{R}$
\begin{align}
 \|w_{q+1}-w_{q+1}^{(o)}\|_{C_{[t,t+1]}L^{\frac65}}&\lesssim (\|\mathring{R}_q\|_{L^1_{[t-l,t+1]}L^1}+1)^3\lambda_{q+1}^{14\alpha-\frac23+\frac{32}3\gamma}\notag\\ 
 &\lesssim (\|\mathring{R}_q\|_{L^1_{[t-l,t+1]}L^1}+1)^3\lambda_{q+1}^{-\gamma},\label{4.20*}
\end{align}
where we used the conditions on the parameters to deduce $14\alpha<\gamma$ and $\gamma<\frac1{19}$. 
Then together with \eqref{4.18}, (\ref{4.27}) and interpolation we  imply
\begin{align*}
\|{w}_{q+1}&\|_{C_{[t,t+1]}W^{\frac12,\frac65}}\\ 
& \lesssim  \|{w}_{q+1}-w_{q+1}^{(o)}\|_{C_{[t,t+1]}L^{\frac65}}^{1/2}\|{w}_{q+1}-w_{q+1}^{(o)}\|_{C_{[t,t+1]}W^{1,\frac65}}^{1/2}+\|{w}_{q+1}^{(o)}\|_{C_{[t,t+1]}L^{\frac65}}^{1/2}\|{w}_{q+1}^{(o)}\|_{C_{[t,t+1]}W^{1,\frac65}}^{1/2}\\
&\lesssim(\|\mathring{R}_q\|_{L^1_{[t-l,t+1]}L^1}+1)^{7/2}(\lambda_{q+1}^{27\alpha-\frac{1}{6}+\frac{35}3\gamma}+\lambda_{q+1}^{46\alpha-2\gamma})\\
&\lesssim(\|\mathring{R}_q\|_{L^1_{[t-l,t+1]}L^1}+1)^{7/2}\lambda_{q+1}^{-\gamma},
\end{align*}
where we used the conditions on the parameters to deduce $27\alpha-\frac{1}{6}+\frac{35}3\gamma<-\gamma,46\alpha<\gamma$. Together with (\ref{3.9}), (\ref{3.11}) we obtain
\begin{align}
 \$ v_{q+1}^2-v_q^2 \$ _{W^{\frac12,\frac65},2r}&\lesssim  \$ v_{q}^2-v_l^2 \$ _{W^{\frac12,\frac65},2r}+ \$ w_{q+1} \$ _{W^{\frac12,\frac65},2r}\notag\\
&\lesssim l^{\frac12} \fm\lambda_q^4+\lambda_{q+1}^{-\gamma}(7^{q}\cdot \frac{7}2 \bar{M}_L r)^{\frac{21}2\cdot {7}^q}\leq \lambda_{q+1}^{-\alpha/2}\leq \delta_{q+1}^{1/2},\notag
\end{align}
\begin{align}
 \$ v_{q+1}^2-v_q^2 \$ _{W^{\frac12,\frac65},m}&\lesssim  \$ v_{q}^2-v_l^2 \$ _{W^{\frac12,\frac65},m}+ \$ w_{q+1} \$ _{W^{\frac12,\frac65},m}\notag\\
&\lesssim l^{\frac12}\lambda_q^{\frac72}(7^{q-1}\cdot 5 \bar{M}_L m)^{15\cdot {7}^{q-1}}+\lambda_{q+1}^{-\gamma}(7^{q}\cdot \frac{7}2 \bar{M}_L m)^{\frac{21}2\cdot {7}^q}\notag\\ 
&\leq \lambda_{q+1}^{-\alpha/2}(7^{q}\cdot \frac{7}2 \bar{M}_L m)^{\frac{21}2\cdot {7}^q},\label{v1265}
\end{align}
where  we used the conditions on the parameters  to deduce  $\lambda_{q+1}^{-\gamma}(7^{q}\cdot \frac{7}2 \bar{M}_L r)^{\frac{21}2\cdot{7}^q}\leq \lambda_{q+1}^{-\gamma} \lambda_q^{11}\ll \lambda_{q+1}^{-\alpha/2}, l^{\frac12} \fm\lambda_q^4\leq l^{\frac12}\lambda_q^5\ll \lambda_{q+1}^{-\alpha/2}$.  We chose $a$ large enough to absorb the implicit constant. Hence together with \eqref{ab2}, (\ref{3.16}) and \eqref{3.10} at the level $q+1$ hold.

\subsection{Definition of the Reynolds stress}\label{dotrs}

Combine (\ref{3.6}) at the level $q+1$, (\ref{4.5}) and the definition of $w_{q+1}$ we get
\begin{align}
\div\mathring{R}_{q+1}&-\nabla p_{q+1}^2+\nabla p_l\notag\\
&=-\Delta w_{q+1}+\partial_t(w_{q+1}^{(p)}+w_{q+1}^{(c)})+\div((v_l^2+v_{q+1}^1)\otimes w_{q+1}+w_{q+1}\otimes (v_l^2+v_{q+1}^1))\notag\\
&\ \ \ \ \ \ \ \ \ \ \ \ \ \ \ \ \ \ \ \ \ \ \ \ \ \ \ \ \ \ \ \ \ \ \ \ \ \ \ \ \ \ \ \ \ \ \ \ \ \ \ \ \ \ \  \ \ \ \ \ \ \ \ \ \ \ \ \ (:=\div R_{lin}+\nabla p_{lin})\notag\\
&+\div((w_{q+1}^{(c)}+w_{q+1}^{(t)})\otimes w_{q+1}+w_{q+1}^{(p)}\otimes (w_{q+1}^{(c)}+w_{q+1}^{(t)}))\ (:=\div R_{cor}+\nabla p_{cor})\notag\\
&+\partial_tw_{q+1}^{(t)}+\div(w_{q+1}^{(p)}\otimes w_{q+1}^{(p)}+\mathring{R}_l)\ (:=\div R_{osc}+\nabla p_{osc})\notag\\
&+\div((v_{q+1}-v_q)\otimes\cL_{R}z+\cL_{R}z\otimes(v_{q+1}-v_q)\notag\\
&\ \ \ \ \ \ \ \ \ \ \  \ \ \ +(v_{q+1}-v_q)\succcurlyeq\!\!\!\!\!\!\!\bigcirc\cH_{R}z+\cH_{R}z\preccurlyeq\!\!\!\!\!\!\!\bigcirc(v_{q+1}-v_q))(:=\div R_{com1}+\nabla p_{com1})\notag\\
&+\div((v_l^2+v_{q+1}^1)\otimes(v_l^2+v_{q+1}^1)-(v_q^1+v_q^2)\otimes(v_q^1+v_q^2))(:=\div R_{com2}+\nabla p_{com2})\notag\\
&+\fc(v_l^1-v_{q+1}^1)(:=\div R_{com3}+\nabla p_{com3})\notag\\
&+\div{R}_{com},\notag
\end{align}
where $R_{com}$ is introduced in Section \ref{moll}.

Applying the inverse divergence operator $\mathcal{R}$, by (\ref{rdeltav}) we could define
\begin{align*}
R_{lin}:&=-\nabla w_{q+1}-\nabla^Tw_{q+1}+ \mathcal{R}(\partial_t(w_{q+1}^{(p)}+w_{q+1}^{(c)}))+(v_l^2+v_{q+1}^1)\mathring{\otimes} w_{q+1}+w_{q+1}\mathring{\otimes} (v_l^2+v_{q+1}^1),\\
R_{cor}:&=(w_{q+1}^{(c)}+w_{q+1}^{(t)})\mathring{\otimes} w_{q+1}+w_{q+1}^{(p)}\mathring{\otimes} (w_{q+1}^{(c)}+w_{q+1}^{(t)}),\\
R_{com2}:&=(v_l^2+v_{q+1}^1)\mathring{\otimes}(v_l^2+v_{q+1}^1)-(v_q^1+v_q^2)\mathring{\otimes}(v_q^1+v_q^2),\\
R_{com3}:&=\fc\mathcal{R}(v_l^1-v_{q+1}^1),
\end{align*}
$R_{com1}$ is the trace-free part of the matrix
\begin{align*}
(v_{q+1}-v_q){\otimes}\cL_{R}z+\cL_{R}z{\otimes}(v_{q+1}-v_q) +(v_{q+1}-v_q)\ {\succcurlyeq\!\!\!\!\!\!\!\bigcirc}\cH_{R}z+\cH_{R}z\ {\preccurlyeq\!\!\!\!\!\!\!\bigcirc}(v_{q+1}-v_q).
\end{align*}

 The definition of $R_{osc}:= R_{osc,x}+R_{osc,a}+R_{osc,t}$ follows similarly as
in \cite[Section 5.6]{LZ23}, we omit the details and define

\begin{align*}
    R_{osc,x}:&=\sum_{\xi\in\Lambda}g^2_{(\xi)}\mathcal{B}\left(\nabla a_{(\xi)}^2,\mathbb{P}_{\neq0}(W_{(\xi)}\otimes W_{(\xi)})\right),\\
  R_{osc,a}:&=-\theta^{-1}\sigma\sum_{\xi\in\Lambda}\mathcal{B}\left(\partial_t(a_{(\xi)}^2g_{(\xi)}),|W_{(\xi)}|^2\xi\right),\\
R_{osc,t}:&=-\sigma^{-1}\sum_{\xi\in\Lambda}\mathcal{R}\left(h_{(\xi)}\div\partial_t(a_{(\xi)}^2\xi\otimes\xi)\right).
\end{align*}
 Finally we define the Reynolds stress at the level $q+1$ by
$$\mathring{R}_{q+1}=R_{lin}+R_{cor}+R_{osc}+
R_{com}+R_{com1}+R_{com2}+R_{com3}.$$
We observe that by construction, $\mathring{R}_{q+1}$ is $(\mathcal{F}_t)_{t\in\mathbb{R}}$-adapted.

\subsection{Verification of the inductive estimate for $\mathring{R}_{q+1}$}\label{votiefr}
To conclude the proof of Proposition \ref{prop:3.1} we shall verify (\ref{3.11}) at the level $q+1$.  To this end, we estimate each terms in the definition of $\mathring{R}_{q+1}$ separately.

First we estimate the linear error $R_{lin}$. By (\ref{4.25})-(\ref{4.28})  we obtain for $s\in[t,t+1], t\in\mathbb{R}$ and $\epsilon>0$ small enough
\begin{align}
\|\nabla w_{q+1}(s)+\nabla^T w_{q+1}(s)\|_{L^{1}}&\lesssim\|{w}_{q+1}(s)\|_{W^{1,1+\epsilon}}\notag\\
&\lesssim(\|\mathring{R}_q
\|_{L^1_{[t,t+1]}L^1}+1)^4 (
 \lambda_{q+1}^{40\alpha+6\gamma+2\epsilon}\sum_{\xi\in\Lambda}|g_{(\xi)}(s)|+\lambda_{q+1}^{52\alpha-2\gamma}),\notag
\end{align}
which together with (\ref{gwnp}) implies 
\begin{align}
\|\nabla w_{q+1}+\nabla^T w_{q+1}\|_{L^1_{[t,t+1]}L^{1}}&\lesssim(\|\mathring{R}_q
\|_{L^1_{[t,t+1]}L^1}+1)^4 (
 \lambda_{q+1}^{40\alpha-2\gamma+2\epsilon}+\lambda_{q+1}^{52\alpha-2\gamma})\notag\\
 &\lesssim (\|\mathring{R}_q
\|_{L^1_{[t,t+1]}L^1}+1)^4 \lambda_{q+1}^{-\gamma},\notag
\end{align}
where we used the conditions on the parameters to deduce $2\epsilon<\alpha,52\alpha<\gamma$.

By (\ref{wp+wc}) and (\ref{B.10}) we obtain
\begin{align*}
    \mathcal{R}(\partial_t(w_{q+1}^{(p)}+w_{q+1}^{(c)}))
    &=\sigma^{-1}\sum_{\xi\in\Lambda}\mathcal{R}\nabla^{\perp}[\partial_t(a_{(\xi)}g_{(\xi)})\Psi_{(\xi)}]+\sigma^{-2}\theta\sum_{\xi\in\Lambda}\mathcal{R}\nabla^{\perp}[a_{(\xi)}g_{(\xi)}^2(\xi\cdot\nabla)\Psi_{(\xi)}].
\end{align*}
Now by the fact that $\mathcal{R}\nabla^{\perp}$ is $L^p\to L^p$ bounded for any $p>1$, together with (\ref{4.8}), (\ref{2.4}) and (\ref{2.4*}) we obtain for $s\in[t,t+1], t\in\mathbb{R}$ and $\epsilon>0$ small enough
\begin{align}
&\|\mathcal{R}(\partial_t(w_{q+1}^{(p)}+w_{q+1}^{(c)}))(s)\|_{L^1}\notag\\
&\lesssim \sigma^{-1}\sum_{\xi\in\Lambda}\|\partial_t(a_{(\xi)}g_{(\xi)})\Psi_{(\xi)}(s)\|_{L^{1+\epsilon}}+\sigma^{-2}\theta\sum_{\xi\in\Lambda}\|a_{(\xi)}g_{(\xi)}^2(\xi\cdot\nabla)\Psi_{(\xi)}(s)\|_{L^{1+\epsilon}}\notag\\
&\lesssim\sum_{\xi\in\Lambda}\|a_{(\xi)}\|_{C_{[t,t+1],x}^1}\Big(\sigma^{-1}\|\Psi_{(\xi)}\|_{C_{[t,t+1]}L^{1+\epsilon}}(|\partial_tg_{(\xi)}(s)|+|g_{(\xi)}(s)|)\notag\\
&\ \ \ \ \ \ \ \ \ \ \ \ \ \ \ \ +\sigma^{-2}\theta\|(\xi\cdot\nabla)\Psi_{(\xi)}\|_{C_{[t,t+1]}L^{1+\epsilon}}g^2_{(\xi)}(s)\Big{)}\notag\\
&\lesssim (\|\mathring{R}_q
\|_{L^1_{[t,t+1]}L^1}+1)^2l^{-13}\sum_{\xi\in\Lambda}\left(\sigma^{-1
}\nu^{-\frac12}\mu^{-\frac32}(|\partial_tg_{(\xi)}(s)|+|g_{(\xi)}(s)|)+\sigma^{-1}\theta\nu^{\frac12}\mu^{-\frac32}g^2_{(\xi)}(s)
\right)(\nu\mu)^{\epsilon}\notag\\
&\lesssim (\|\mathring{R}_q
\|_{L^1_{[t,t+1]}L^1}+1)^2 \sum_{\xi\in\Lambda}\left(\lambda_{q+1}^{26\alpha-2+2\gamma+2\epsilon}(|\partial_tg_{(\xi)}(s)|+|g_{(\xi)}(s)|)+\lambda_{q+1}^{26\alpha-\gamma+2\epsilon}g^2_{(\xi)}(s)\right),\notag
\end{align}
which together with (\ref{gwnp}) implies 
\begin{align}
\|\mathcal{R}(\partial_t(w_{q+1}^{(p)}+w_{q+1}^{(c)}))\|_{L^1_{[t,t+1]}L^1}&\lesssim (\|\mathring{R}_q
\|_{L^1_{[t,t+1]}L^1}+1)^2(\lambda_{q+1}^{26\alpha-2+12\gamma+2\epsilon} +\lambda_{q+1}^{26\alpha-\gamma+2\epsilon})\notag\\
&\lesssim(\|\mathring{R}_q
\|_{L^1_{[t,t+1]}L^1}+1)^2\lambda_{q+1}^{-\frac12\gamma},\notag
\end{align}
where we chose $2\epsilon<\alpha$ and used the conditions on the parameters to deduce $54\alpha<\gamma$.

By \eqref{4.18} and (\ref{4.20*}) we obtain
\begin{align}
\|(v_l^2+v_{q+1}^1)&\mathring{\otimes} w_{q+1}+w_{q+1}\mathring{\otimes} (v_l^2+v_{q+1}^1)\|_{L^1_{[t,t+1]}L^{1}}\notag\\
&\lesssim
\|v_l^2+v_{q+1}^1\|_{C_{[t,t+1]} L^\infty}\|w_{q+1}\|_{C_{[t,t+1]}L^\frac65}\notag\\
&\lesssim(\|\mathring{R}_q\|_{L^1_{[t-l,t+1]}L^1}+1)^3\lambda_{q+1}^{-\gamma}(\|v_q^2\|_{C_{[t-l,t+1]} L^\infty}+\|v_{q+1}^1\|_{C_{[t,t+1]} L^\infty}).\notag
\end{align}
Then together with all the estimates above  we deduce 
\begin{align}
\|{R}_{lin}\|_{L^1_{[t,t+1]}L^1}&\lesssim (\|\mathring{R}_q\|_{L^1_{[t-l,t+1]}L^1}^4+1)\lambda_{q+1}^{-\frac\gamma2}\notag\\
&\ \ \ \  +(\|\mathring{R}_q\|_{L^1_{[t-l,t+1]}L^1}^3+1)\lambda_{q+1}^{-\gamma}(\|v_q^2\|_{C_{[t-l,t+1]} L^\infty}+\|v_{q+1}^1\|_{C_{[t,t+1]} L^\infty}).\notag
\end{align}

By (\ref{3.9}), (\ref{3.11}) and (\ref{4.4.1}) we obtain
\begin{align}
  \$ {R}_{lin} \$ _{L^1_tL^1,r}&\lesssim ( \$ \mathring{R}_q \$ _{L^1_tL^1,4r}^4+1)\lambda_{q+1}^{-\frac\gamma2} +( \$ \mathring{R}_q \$ _{L_t^1L^1,6r}^3+1)\lambda_{q+1}^{-\gamma}( \$ v_q^2 \$ _{C_{t,x}^1,2r}+ \$ v_{q+1}^1 \$ _{L^\infty,2r})\notag\\
 &\lesssim \lambda_{q+1}^{-\frac\gamma2}(7^{q}\cdot 4 \bar{M}_L r)^{12\cdot {7}^q}+(7^{q}\cdot 6 \bar{M}_L r)^{9\cdot {7}^q}\lambda_{q+1}^{-\gamma}(\fm\lambda_q^4+\lambda_{q+1}^{2\alpha}r^{1/2}L)\notag\\
 &\lesssim \lambda_{q+1}^{-\frac\gamma2}(\lambda_q^{12}+\lambda_q^{13}+\lambda_q^{9}\lambda_{q+1}^{2\alpha}) M_L^{1/2}\lesssim \lambda_{q+1}^{3\alpha-\frac\gamma2} M_L^{1/2}
 \leq\frac{1}{7}\delta_{q+2}M_L, \notag 
\end{align}
 where we used the conditions on the parameters $a>(42 \bar{M}_L r)^d,\alpha b>9$ to deduce $(7^{q}\cdot 6 \bar{M}_L r)^{{7}^q}\leq \lambda_q,\lambda_q^{9}\leq \lambda_{q+1}^{\alpha}$, and we used $\fm+r^{1/2}L\leq M_L^{1/2}, 3\alpha-\frac\gamma2<-2\beta b$. Similarly  together with \eqref{4.4.2} we have
 \begin{align}
  \$ {R}_{lin} \$ _{L_t^1L^1,m}&\lesssim ( \$ \mathring{R}_q \$ _{L^1_tL^1,4m}^4+1)\lambda_{q+1}^{-\frac\gamma2}+( \$ \mathring{R}_q \$ _{L^1_tL^1,7m}^3+1)\lambda_{q+1}^{-\gamma}( \$ v_q^2 \$ _{C_{t,x}^1,2m}+ \$ v_{q+1}^1 \$ _{L^\infty,\frac74m})\notag\\
 &\lesssim \lambda_{q+1}^{-\frac\gamma2}(7^{q}\cdot 4 \bar{M}_L m)^{12\cdot {7}^q}+\lambda_{q+1}^{-\gamma}(7^{q}\cdot 7 \bar{M}_L m)^{9\cdot {7}^q}\Big{(}\lambda_q^{7/2}(7^{q-1}\cdot 10 \bar{M}_L m)^{15\cdot {7}^{q-1}}\notag\\
 &\ \ \ \ \  \ +\lambda_{q+1}^{\alpha/2}(7^{q}\cdot7 \bar{M}_L m)^{\frac{21}2\cdot {7}^q}m^{1/2}L\Big{)}\leq \frac{1}{7}(7^{q+1}\cdot  \bar{M}_L m)^{3\cdot {7}^{q+1}},\notag
\end{align}
 where we used the conditions on the parameters $\alpha b>7,\gamma>\alpha/2$ to deduce $\lambda_q^{7/2}\leq \lambda_{q+1}^{\alpha/2}<\lambda_{q+1}^{\gamma}$ and chose $a$ large enough to absorb the  implicit constant. 

The corrector error ${R}_{cor}$ is estimated using   (\ref{wc+wt22})
\begin{align}
\|{R}_{cor}\|_{L^1_{[t,t+1]}L^1}&\leq\|w_{q+1}^{(c)}+w_{q+1}^{(t)}\|_{L^2_{[t,t+1]}L^{2}}(\|w_{q+1}\|_{L^2_{[t,t+1]}L^{2}}
+\|w_{q+1}^{(p)}\|_{L^2_{[t,t+1]}L^{2}})\notag\\
&\lesssim(\|\mathring{R}_q\|_{L^1_{[t-l,t+1]}L^1}+1)^3\lambda_{q+1}^{40\alpha-2\gamma}(\|w_{q+1}\|_{L^2_{[t,t+1]}L^{2}}
+\|w_{q+1}^{(p)}\|_{L^2_{[t,t+1]}L^{2}}).\notag
\end{align}

 Then (\ref{3.11}), (\ref{4.20.1}) and (\ref{4.20.2}) yield
\begin{align}
 \$ {R}_{cor} \$ _{L^1_tL^1,r}&\lesssim \lambda_{q+1}^{40\alpha-2\gamma}( \$ \mathring{R}_q \$ _{L^1_tL^1,6r}^3+1)( \$ w_{q+1} \$ _{L^2_tL^{2},2r}+ \$ w_{q+1}^{(p)} \$ _{L^2_tL^{2},2r})   \notag\\
&\lesssim \lambda_{q+1}^{40\alpha-2\gamma}((7^{q}\cdot 6 \bar{M}_L r)^{9\cdot {7}^q}+1) \delta_{q+1}^{1/2}M_L^{1/2}
\lesssim \lambda_{q+1}^{40\alpha-2\gamma}\lambda_q^{9}M_L^{1/2}\leq\frac{1}{7}\delta_{q+2}M_L,\notag\\
 \$ {R}_{cor} \$ _{L^1_tL^1,m}&\lesssim \lambda_{q+1}^{40\alpha-2\gamma}( \$ \mathring{R}_q \$ _{L^1_tL^1,6m}^3+1)( \$ w_{q+1} \$ _{L^2_tL^{2},2m}+ \$ w_{q+1}^{(p)} \$ _{L^2_tL^{2},2m})  \notag\\
&\lesssim\lambda_{q+1}^{40\alpha-2\gamma}((7^{q}\cdot 6 \bar{M}_L m)^{9\cdot {7}^q}+1)(7^q\cdot6 \bar{M}_L m)^{9\cdot {7}^q}\leq \frac{1}{7}(7^{q+1}\cdot  \bar{M}_L m)^{3\cdot {7}^{q+1}},\notag
\end{align}
where we used the conditions on the parameters $ a>(42 \bar{M}_L r)^d$ to deduce $(7^{q}\cdot 6 \bar{M}_L r)^{ {7}^q}\leq\lambda_q$,  and used  $
40\alpha-2\gamma+\frac9b<-2\beta b$. In the last inequality we chose $a$ large enough to absorb the  implicit constant.

We continue with oscillation error  ${R}_{osc}=R_{osc,x}+R_{osc,a}+R_{osc,t}$. In order to bound the first term, we apply (\ref{bb1}), (\ref{bb}), (\ref{4.8}), (\ref{gwnp}) and (\ref{2.4}) to deduce
\begin{align}
\|R_{osc,x}\|_{L^1_{[t,t+1]}L^1}&\lesssim\sum_{\xi\in\Lambda}\|g^2_{(\xi)}\|_{L^1_{[t,t+1]}}\|\mathcal{B}(\nabla a_{(\xi)}^2,\mathbb{P}_{\neq0}(W_{(\xi)}\otimes W_{(\xi)}))\|_{C_{[t,t+1]}L^1}\notag\\
&\lesssim \sum_{\xi\in\Lambda}\|g_{(\xi)}\|_{L^2_{[t,t+1]}}^2\|\nabla a_{(\xi)}^2\|_{C_{[t,t+1]}C_x^1}\|\mathcal{R}(W_{(\xi)}\otimes W_{(\xi)})\|_{C_{[t,t+1]}L^{1}}\notag\\
&\lesssim (\|\mathring{R}_q\|_{L^1_{[t-l,t+1]}L^1}+1)^4l^{-26}\sigma^{-1}\|W_{(\xi)}\|_{C_{[t,t+1]}L^2}^2\lesssim \lambda_{q+1}^{52\alpha-2\gamma}(\|\mathring{R}_q\|_{L^1_{[t-l,t+1]}L^1}+1)^4,\notag
\end{align}
\begin{align}
\|{R}_{osc,a}\|_{L^1_{[t,t+1]}L^1}&\lesssim \theta^{-1}\sigma\sum_{\xi\in\Lambda}\|\mathcal{B}\left(\partial_t(a_{(\xi)}^2g_{(\xi)}),|W_{(\xi)}|^2\xi\right)\|_{L^1_{[t,t+1]}L^1}\notag\\
&\lesssim\theta^{-1}\sigma\sum_{\xi\in\Lambda}\|\partial_t(a_{(\xi)}^2g_{(\xi)})\|_{L^1_{[t,t+1]}C_x^1}\|\mathcal{R}(|W_{(\xi)}|^2\xi)\|_{C_{[t,t+1]}L^1}\notag\\
&\lesssim \theta^{-1}\sum_{\xi\in\Lambda}\|a_{(\xi)}^2\|_{C_{[t,t+1],x}^2}\(\|\partial_t g_{(\xi)}\|_{L^1_{[t,t+1]}}+\| g_{(\xi)}\|_{L^1_{[t,t+1]}}\)\|W_{(\xi)}\|_{C_{[t,t+1]}L^2}^2
,\notag\\
&\lesssim \theta^{-1}\sigma\eta^{1/2}l^{-26}(\|\mathring{R}_q\|_{L^1_{[t-l,t+1]}L^1}+1)^4\lesssim \lambda_{q+1}^{52\alpha-1+5\gamma}(\|\mathring{R}_q\|_{L^1_{[t-l,t+1]}L^1}+1)^4.\notag
\end{align}

By (\ref{4.8}) and (\ref{hlinfty}) we have
\begin{align}
\|{R}_{osc,t}\|_{L^1_{[t,t+1]}L^1}
&\lesssim\sigma^{-1}\sum_{\xi\in\Lambda}\|a_{(\xi)}^2\|_{C_{[t,t+1],x}^1}\|h_{(\xi)}\|_{L_{[t,t+1]}^\infty}\lesssim\lambda_{q+1}^{40\alpha-2\gamma}(\|\mathring{R}_q\|_{L^1_{[t-l,t+1]}L^1}+1)^3.\notag
\end{align}
These bounds yield
\begin{align}
   \|{R}_{osc}\|_{L^1_{[t,t+1]}L^1}
   \leq\lambda_{q+1}^{-\gamma}(\|\mathring{R}_q\|_{L^1_{[t-l,t+1]}L^1}+1)^4,\notag
\end{align}
where we used the conditions on the parameters to deduce $52\alpha<\gamma<\frac{1}{7}$. Then by (\ref{3.11})  we obtain
\begin{align}
 \$ {R}_{osc} \$ _{L^1_tL^1,r}
   \leq\lambda_{q+1}^{-\gamma}( \$ \mathring{R}_q \$ _{L^1_tL^1,4r}^4+1)\leq\lambda_{q+1}^{-\gamma} (7^{q}\cdot 4 \bar{M}_L r)^{12\cdot {7}^q}\leq\lambda_{q+1}^{-\gamma} \lambda_q^{12}\leq\frac{1}{7}\delta_{q+2}M_L,\notag \\
   \$ {R}_{osc} \$ _{L^1_tL^1,m}
   \leq\lambda_{q+1}^{-\gamma}( \$ \mathring{R}_q \$ _{L^1_tL^1,4m}^4+1)\leq \lambda_{q+1}^{-\gamma}(7^{q}\cdot 4 \bar{M}_L m)^{12\cdot {7}^q}\leq\frac{1}{7} (7^{q+1}\cdot  \bar{M}_L m)^{3\cdot {7}^{q+1}},\notag
\end{align}
where we used the conditions on the parameters $a>(28 \bar{M}_L r)^d$ to deduce $ (7^{q}\cdot 4 \bar{M}_L r)^{{7}^q}\leq \lambda_q$ and used $\gamma+\frac{12}b<-2\beta b$.

For ${R}_{com1}$, by Lemma \ref{lem:2.3} with $0<\kappa<\frac13-\fa-\fa_0$  and Proposition \ref{deltaz} we obtain for $t\in\R$
\begin{align}
\|{R}_{com1}\|_{L^1_{[t,t+1]}L^1}&\lesssim
(\|v_{q+1}^2-v_q^2\|_{C_{[t,t+1]}W^{\frac12,\frac65}}+\|v_{q+1}^1-v_q^1\|_{C_{[t,t+1]}B^{\frac23-\fa-\kappa}_{\frac32,\infty}})\notag\\
&\ \ \ \ \ \ \ \ \ \ \times(\|\cL_{R}z\|_{C_{[t,t+1]}L^\infty}+\|z\|_{C_{[t,t+1]}C^{-\fa-\kappa}}),\notag\\
&\lesssim(\|v_{q+1}^2-v_q^2\|_{C_{[t,t+1]}W^{\frac12,\frac65}}+\|v_{q+1}^1-v_q^1\|_{C_{[t,t+1]}B^{\frac23-\fa-\kappa}_{\frac32,\infty}}) \fC^3(1+\|z\|_{C_{[t,t+1]}C^{-\fa-\kappa}}^4),\notag
\end{align}
which together with (\ref{3.16}), (\ref{3.12}) implies
\begin{align}
 \$ {R}_{com1} \$ _{L^1_tL^1,r}&\lesssim   ( \$ v_{q+1}^2-v_q^2 \$ _{W^{\frac12,\frac65},2r}+ \$ v_{q+1}^1-v_q^1 \$ _{B^{\frac23-\fa-\kappa}_{\frac32,\infty},2r})\fC^3( \$ z \$ _{C^{-\fa-\kappa},8r}^4+1)\notag\\
&\lesssim (\lambda_{q+1}^{-\frac\alpha2}+\lambda_{q}^{-\frac\alpha{24}}) \fC^3r^2L^4\lesssim  \lambda_{q}^{-\frac\alpha{24}} \fC^3r^2L^4\leq \frac{1}{7}\delta_{q+2}M_L,\notag
\end{align}
where we used the conditions on the parameters to deduce $\alpha>48\beta b^2$. Similarly by \eqref{3.13.2} and \eqref{v1265}
\begin{align}
 \$ {R}_{com1} \$ _{L^1_tL^1,m}&\lesssim( \$ v_{q+1}^2-v_q^2 \$ _{W^{\frac12,\frac65},2m}+ \$ v_{q+1}^1 \$ _{B^{\frac23-\fa-\kappa}_{\frac32,\infty},2m}+ \$ v_q^1 \$ _{B^{\frac23-\fa-\kappa}_{\frac32,\infty},2m}) \fC^3( \$ z \$ _{C^{-\fa-\kappa},8m}^4+1)\notag\\
&\lesssim(7^{q}\cdot 7 \bar{M}_L m)^{\frac{21}2\cdot {7}^q}\fC^4m^2L^4\leq \frac{1}{7} (7^{q+1}\cdot  \bar{M}_L m)^{3\cdot{7}^{q+1}}.\notag
\end{align}
where we chose $a$ large enough to absorb the implicit constant. 

For ${R}_{com2}$, by (\ref{3.8.1}),  (\ref{3.13.1}), (\ref{3.12})  and \eqref{vq-vl} we obtain
\begin{align}
 \$ {R}_{com2} \$ _{L^1_tL^1,r}&\lesssim ( \$ v_{q+1}^1 \$ _{L^2,2r}+ \$ v_{q}^1 \$ _{L^2,2r}+ \$ v_q^2 \$ _{L^2_tL^2,2r})( \$ v_{q+1}^1-v_q^1 \$ _{L^2,2r}+l \$ v_q^2 \$ _{C_{t,x}^1,2r})\notag\\
&\lesssim  M_L^{1/2}
(\lambda_{q}^{-\frac\alpha{24}}+l \fm\lambda_{q}^4)\lesssim  M_L
\lambda_{q}^{-\frac\alpha{24}}\leq \frac{1}{7}\delta_{q+2}M_L,\notag
\end{align} 
where we used the conditions on the parameters to deduce $l\lambda_q^4\leq\lambda_{q+1}^{-\alpha}$, $\alpha>48\beta b^2$  and $\fm\leq M_L^{1/2}$. 
 Similarly by \eqref{3.8.2}, \eqref{3.9}, \eqref{3.13.2} and \eqref{v1l2} we have

 \begin{align}
 \$ {R}_{com2} \$ _{L^1_tL^1,m}&\lesssim ( \$ v_{q+1}^1 \$ _{L^2,2m}+ \$ v_{q}^1 \$ _{L^2,2m}+ \$ v_q^2 \$ _{L^2_tL^2,2m})( \$ v_{q+1}^1-v_q^1 \$ _{L^2,2m}+l \$ v_q^2 \$ _{C_{t,x}^1,2m})\notag\\
&\lesssim 
(\lambda_{q}^{-\frac\alpha{24}}+l\lambda_{q}^{7/2})(7^{q}\cdot 7 \bar{M}_L m)^{21\cdot {7}^{q}}\leq  \frac{1}{7} (7^{q+1}\cdot \bar{M}_Lm)^{3\cdot{7}^{q+1}}.\notag
\end{align} 

We continue with ${R}_{com3}$. By (\ref{bb1}), \eqref{3.12}, \eqref{4.4.1}, \eqref{v1l2} and \eqref{4.4.2} we have
\begin{align*}
    \$ {R}_{com3} \$ _{L^1_tL^1,r}&\lesssim \fc \$ v_l^1-v_{q+1}^1 \$ _{L^2,r}\lesssim \fc( \$ v_{q+1}^1-v_q^1 \$ _{L^2,r}+ \$ v_l^1-v_q^1 \$ _{L^2,r})\\
   &\lesssim \fc(\lambda_{q}^{-\frac\alpha{24}}+l^{\frac16}\lambda_q^{2\alpha}r^{1/2}L)\lesssim \fc \lambda_{q}^{-\frac\alpha{24}}M_L\leq\frac{1}{7}\delta_{q+2}M_L,\notag\\
     \$ {R}_{com3} \$ _{L^1_tL^1,m}&\lesssim \fc( \$ v_{q+1}^1-v_q^1 \$ _{L^2,m}+ \$ v_l^1-v_q^1 \$ _{L^2,m})\notag\\
    &\lesssim \fc(\lambda_{q}^{-\frac\alpha{24}}(7^{q}\cdot \frac{7}2 \bar{M}_L m)^{\frac{21}2\cdot {7}^{q}}+l^{\frac16}\lambda_q^{\alpha/2}(7^{q-1}\cdot 4 \bar{M}_L m)^{\frac{21}2\cdot {7}^{q-1}}m^{1/2}L)\notag\\ 
    &\leq\frac{1}{7}  (7^{q+1}\cdot  \bar{M}_L m)^{3\cdot{7}^{q+1}},
\end{align*}
where we used the conditions on the parameters  to deduce $l^{\frac16}\lambda_q^{2\alpha}\leq \lambda_q^{-\frac13}\lambda_q^{2\alpha}\leq\lambda_{q}^{-\frac\alpha{24}}$, $\alpha>48\beta b^2$ and $\fc\lambda_{q}^{-\frac\alpha{24}+2\beta b^2}\ll1$.

Finally we consider each term in ${R}_{com}$  separately, which is defined in Section \ref{moll}. For $I_1=(v_q^1+v_q^2)\otimes\cL_{R}z+\cL_{R}z\otimes(v_q^1+v_q^2)$, by Proposition \ref{deltaz}
 we have
\begin{align}
&\sup_{s\in[t,t+1]}\|(I_1-I_1*_x\phi_l*_t\varphi_l)(s)\|_{L^1}\notag\\
&\lesssim l^{\frac{\fa_0}2}(\|v_q^1+v_q^2\|_{C_{[t-l,t+1]}^{\frac{\fa_0}2}L^2}+\|v_q^1+v_q^2\|_{C_{{[t-l,t+1]}}H^{\frac{\fa_0}2}})\|\cL_{R}z\|_{C^{\frac{\fa_0}2}_{[t-1,t+1]}C^{\frac13}}\notag\\
&\lesssim l^{\frac{\fa_0}2}(\|v_q^1+v_q^2\|_{C_{[t-l,t+1]}^{\frac{\fa_0}2}L^2}+\|v_q^1+v_q^2\|_{C_{{[t-l,t+1]}}H^{\frac{\fa_0}2}}) \notag\\ 
&\ \ \ \ \ \ \ \ \ \ \times \fC^{6}(1+\sup_{s\in[t-1,t+1]}\|z(s)\|_{C^{-\fa-\kappa}}) ^{6} \|z\|_{C_{[t-1,t+1]}^{\frac{\fa_0}{2}}C^{-\fa-\fa_0-\kk}}.\notag
\end{align}

By \eqref{3.7}, \eqref{3.9}, \eqref{4.4.1}, \eqref{4.4.2} and  $0<\fa_0\leq\frac1{12}$ we obtain
\begin{align*}
     \$ I_1-I_1*_x\phi_l*_t\varphi_l \$ _{L^1_tL^1,r}&\lesssim l^{\frac{\fa_0}{2}}( \fm\lambda_q^4+\lambda_q^{2\alpha}r^{1/2}L) \fC^{6}(r^{\frac12}L)^{7}\lesssim l^{\frac{\fa_0}{2}}\lambda_q^4\fC^6 M_L\leq\frac{1}{21}\delta_{q+2}M_L,
\end{align*}
\begin{align*}
     \$ I_1&-I_1*_x\phi_l*_t\varphi_l \$ _{L^1_tL^1,m}\\
     &\lesssim l^{\frac{\fa_0}{2}}\(\lambda_q^{7/2}(7^{q-1}\cdot 10 \bar{M}_L m)^{15\cdot {7}^{q-1}}+\lambda_q^{\alpha/2}(7^{q-1}\cdot 8 \bar{M}_L m)^{\frac{21}2\cdot {7}^{q-1}}m^{1/2}L\)\fC^6m^{7/2}L^7\\
    &\leq  \frac{1}{21} (7^{q+1}\cdot  \bar{M}_L m)^{3\cdot {7}^{q+1}},
\end{align*}
where we used the conditions on the parameters to deduce $ l^{\frac{\fa_0}{2}}\lambda_q^4\fC^6\ll \lambda_{q+1}^{-\frac{\fa_0}{4}\alpha}\fC^6\ll \delta_{q+2}$ by choosing $\alpha b>\frac{8}{\fa_0}$ and $\alpha>\frac{8}{\fa_0}\beta b$.

For $I_2=(v_q^1+v_q^2)\otimes(v_q^1+v_q^2)$ we have
\begin{align}
\sup_{s\in[t,t+1]}\|(I_2&-I_2*_x\phi_l*_t\varphi_l)(s)\|_{L^1}\notag\\
&\lesssim
l^{\fa_0}(\|v_q^1+v_q^2\|_{C_{[t-l,t+1]}^{\fa_0}L^2}+\|v_q^1+v_q^2\|_{C_{{[t-l,t+1]}}H^{\fa_0}})(\|v_q^1\|_{C_{[t-l,t+1]}L^2}+\|v_q^2\|_{C_{[t-l,t+1],x}^1}).\notag
\end{align}
Thus by  (\ref{3.9}), \eqref{4.4.1}, \eqref{4.4.2}  and  $0<\fa_0\leq\frac1{12}$  we obtain
\begin{align*}
     \$ I_2-I_2*_x\phi_l*_t\varphi_l \$ _{L^1_tL^1,r}\lesssim l^{\fa_0}(\fm\lambda_q^4+\lambda_q^{2\alpha}r^{1/2}L)^2\lesssim l^{\fa_0}\lambda_q^8M_L\leq\frac{1}{21}\delta_{q+2}M_L,
\end{align*}
\begin{align*}
     \$ I_2-I_2*_x\phi_l*_t\varphi_l \$ _{L^1_tL^1,m}&\lesssim l^{\fa_0}\(\lambda_q^{7/2}(7^{q-1}\cdot 10 \bar{M}_L m)^{15\cdot {7}^{q-1}}+\lambda_q^{\alpha/2}(7^{q-1}\cdot 8 \bar{M}_L m)^{\frac{21}2\cdot {7}^{q-1}}m^{1/2}L\)^2\\
    &\leq  \frac{1}{21} (7^{q+1}\cdot  \bar{M}_L m)^{3\cdot {7}^{q+1}},
\end{align*}
where we used the conditions on the parameters  to deduce $l^{\fa_0}\lambda_q^8\ll \delta_{q+2}$ by choosing $\alpha b>\frac{8}{\fa_0}$ and $\alpha>\frac{4}{\fa_0}\beta b$.

For $I_3=(v_q^1+v_q^2)\succcurlyeq\!\!\!\!\!\!\!\bigcirc\cH_{R}z+\cH_{R}z\preccurlyeq\!\!\!\!\!\!\!\bigcirc(v_q^1+v_q^2)$, by Proposition \ref{deltaz}  we have for $0<\kk<\frac13-\fa-\fa_0$
\begin{align}
&\sup_{s\in[t,t+1]}\|(I_3-I_3*_x\phi_l*_t\varphi_l)(s)\|_{L^1}\notag\\
&\lesssim l^{\frac{\fa_0}2}\Big{(}(\|v_q^1\|_{C_{[t-l,t+1]}B^{\frac23-\fa-\kappa}_{\frac32,\infty}}+\|v_q^2\|_{C_{[t-l,t+1]}H^{1}})\|z\|_{C_{[t-1,t+1]} C^{-\fa-\kappa}}\notag\\
&\ \ +(\|v_q^1\|_{C_{[t-l,t+1]}^{\frac{1}{6}-\frac{\kappa+\fa}2}B^{\frac13}_{\frac32,\infty}}+\|v_q^2\|_{C_{[t-l,t+1]}^{\frac{1}{6}-\frac{\kappa+\fa}2}H^{\frac13}})\|z\|_{C_{[t-1,t+1]} C^{-\fa-\kappa}}\notag\\
&\ \  +(\|v_q^1\|_{C_{[t-l,t+1]}B^{\frac23-\fa-\kappa}_{\frac32,\infty}}+\|v_q^2\|_{C_{[t-l,t+1]}H^{1}})\|\cH_{R}z\|_{C_{[t-1,t+1]}^{\frac{\fa_0}2}C^{-\fa-\fa_0-\kk}}\Big{)}\notag \\
&\lesssim l^{\frac{\fa_0}2}\left(\|v_q^1\|_{C_{[t-l,t+1]}B^{\frac23-\fa-\kappa}_{\frac32,\infty}}+\|v_q^1\|_{C_{[t-l,t+1]}^{\frac{1}{6}-\frac{\kappa+\fa}2}B^{\frac13}_{\frac32,\infty}}+\|v_q^2\|_{C_{[t-l,t+1],x}^1}\right) \notag\\ 
&\  \ \ \ \  \ \ \times    \fC^{6}(1+\sup_{s\in[t-1,t+1]}\|z(s)\|_{C^{-\fa-\kappa}}) ^{6} \|z\|_{C_{[t-1,t+1]}^{\frac{\fa_0}{2}}C^{-\fa-\fa_0-\kk}}.\notag
\end{align}

Then by (\ref{3.7}), (\ref{3.9}), (\ref{3.13.1}), (\ref{3.13.2}), (\ref{4.2}) and the fact that $ 0<\kappa<\frac1{3}-\fa-\fa_0$ we obtain 
\begin{align*}
     \$ I_3-I_3*_x\phi_l*_t\varphi_l \$ _{L^1_tL^1,r}\lesssim l^{\frac{\fa_0}2}(\fm\lambda_q^4+1) \fC^{6}(r^{\frac12}L)^{7}\lesssim l^{\frac{\fa_0}2}\lambda_q^4 \fC^{6}M_L\leq\frac{1}{21}\delta_{q+2}M_L,
\end{align*}

\begin{align}
     \$ I_3&-I_3*_x\phi_l*_t\varphi_l \$ _{L^1_tL^1,m}\notag\\
     &\lesssim  l^{\frac{\fa_0}2}\(\lambda_q^{7/2}(7^{q-1}\cdot 10 \bar{M}_L m)^{15\cdot {7}^{q-1}}+(7^{q-1}\cdot 7 \bar{M}_L m)^{\frac{21}2\cdot {7}^{q-1}}m^{1/2}L\) \fC^{6}(m^{\frac12}L)^{7}\notag\\
    &\leq \frac{1}{21} (7^{q+1}\cdot  \bar{M}_L m)^{3\cdot {7}^{q+1}},\notag
\end{align}
 where we used the conditions on the parameters to deduce $  l^{\frac{\fa_0}2}\lambda_q^4\fC^{6}\ll \lambda_{q+1}^{-\frac{\fa_0}4\alpha}\fC^6\ll \delta_{q+2}$ by choosing $\alpha b>\frac{8}{\fa_0}$ and $\alpha>\frac{8}{\fa_0}\beta b$.

Summarizing all the above estimates we obtain  (\ref{3.11}) at the level $q+1$:
\begin{align*}
 \$ \mathring{R}_{q+1} \$ _{L^1_tL^1,r}\leq\delta_{q+2}M_L,\ \
 \$ \mathring{R}_{q+1} \$ _{L^1_tL^1,m}\leq  (7^{q+1}\cdot  \bar{M}_L m)^{3\cdot{7}^{q+1}}.
\end{align*}

We finish the proof of Proposition \ref{prop:3.1}.

\section{Proof of Theorem \ref{thm:exist}}\label{sec6}
In this section, we will prove the main result of this paper, namely Theorem \ref{thm:exist}. More specifically, in the previous Theorem \ref{thm:guina}, we have osuccessfully derived a solution to \eqref{1.1} over the interval $[0,\infty)$ with consistent moment estimates. Subsequently, we  follow the similar argument as \cite[Section 4]{HZZ22} to construct infinitely many stationary solutions as well as ergodic stationary solutions  by choosing different $v_0$.

We recall the trajectory space $\mathcal{T}= \(C(\mathbb{R}; W^{\frac13,\frac65})\cap L^2_{\rm loc}(\mathbb{R};L^2) \) 
  \times  C(\mathbb{R};C^{-1-\kappa})$ and let $S_t,t\in\mathbb{R}$ be the shifts on trajectories:
\begin{align*}
S_t(v,B)(\cdot) = (v(\cdot + t), B(\cdot + t) - B(t)),\ \ ( v,B)\in\mathcal{T}.
\end{align*} 

 In the following we additionally assume
\begin{align}
 0<\aleph< \frac{1}{C_1}\(\frac{  1 }{82 (2\pi)^{2-2\fa} }-  \frac{(2\pi)^{2\fa} }{80(4\pi^2+1)}\),\label{defaleph}
\end{align}
where $C_1$ is the implicit constant from \eqref{vkphi} below. Otherwise we replace $\aleph$ by the minimum of $\aleph$ and $\frac1{2C_1}(\frac{  1 }{82 (2\pi)^{2-2\fa} }-  \frac{(2\pi)^{2\fa} }{80(4\pi^2+1)})$.

The first result of this section is existence of stationary solutions.

\bt\label{thm:6.1} 
 Assume \eqref{defaleph}. Let $u$ be a solution obtained in Theorem \ref{thm:guina}  
and satisfying (\ref{eq:guina1})-(\ref{eq:guina3}) with given $ \fa\in[0,\frac13), r\geq1$  and divergence-free, mean-zero $v_0\in C^1(\mathbb{T}^2;\R^2)$.
Then there exists a sequence $T_n\to\infty$
and a stationary solution $((\tilde{\Omega}, \tilde{\mathcal{F}},(\mathcal{F}_t)_{t{ \in\mathbb{R}}}, \tilde{\mathbf{P}}), \tilde{u}, \tilde{B})$ to (\ref{1.1}) in the sense of Definition \ref{def:weaksol} such that
\begin{align*}
\frac{1}{T_n}\int_0^{T_n}\mathcal{L}[S_t(u-z,B)]\dif t\to \mathcal{L}[\tilde{u}-\tilde{z},\tilde{B}]
\end{align*}
weakly in the sense of probability measures on $\mathcal{T}$ as $n \to \infty$, where $\tilde{z}(t)=\int_{-\infty}^te^{(t-s)(\Delta-I)}(-\Delta)^{\fa/2}\dif \tilde{B}_s$.
Moreover, it holds true that for some $\kk>0$ small enough
\begin{align}
 \$ \tilde{u}-\tilde{z}-v_0 \$ _{W^{(\frac23-\fa-\kappa)\wedge\frac12,\frac65},2r}\leq\sqrt{\aleph},\label{6.1}
\end{align}
and for some $\zeta\in(0,1)$, every $N\in\mathbb{N}$
\begin{align}
\mathbf{E}[
\|\tilde{u}-\tilde{z}\|_{L^{2}([-N,N];H^{\zeta})}^{2r}+\|\tilde{u}-\tilde{z}\|_{C^\zeta ([-N,N],W^{ (\frac23-\fa-\kappa)\wedge\frac12,\frac65})}^{2r}]\lesssim N. \label{6.2}
\end{align}

 Moreover,  by choosing different $v_0$, there exist infinitely many stationary solutions  whose distributions are not Gaussian $N(0,\frac12(-\Delta)^{\fa-1})$.
\et

\begin{proof}
For the solution $(u, B)$ of the equation \eqref{1.1} constructed in Theorem \ref{thm:guina}, let us define $v:=u-z$ for convenience, where $z(t)=\int_{-\infty}^te^{(t-s)(\Delta-I)}(-\Delta)^{\fa/2}\dif B_s$.
We define the ergodic averages of the solution $(u, B)$ as the probability measures on the space
of trajectories $\mathcal{T}$
\begin{align*}
\nu_{ \tau,T}=\frac{1}{T}\int_0^{T}\mathcal{L}[S_{t+\tau}(v,B)]\dif t,\ \ \tau,T\geq0.
\end{align*}
By \eqref{eq:guina1} we obtain  for some $\kk>0$  small enough, $\zeta\in(0,1)$ 
\begin{align*}
\sup_{s\in\mathbb{R}}\mathbf{E}\|v\|^{2r}_{C^\zeta_{[-N+s,N+s]}W^{(\frac23-\fa-\kappa)\wedge\frac12,\frac65}}&\leq \sup_{s\in\mathbb{R}}\sum_{i=-N}^{N-1}\mathbf{E}\|v\|^{2r}_{C^\zeta_{[i+s,i+1+s]}W^{(\frac23-\fa-\kappa)\wedge\frac12,\frac65}}\\
&\lesssim N \$ v \$ _{C_t^\zeta W^{(\frac23-\fa-\kappa)\wedge\frac12,\frac65},2r}^{2r}\lesssim N,
\end{align*}
and similarly  
\begin{align*}
\sup_{s\in\mathbb{R}}\mathbf{E}[\|v(\cdot + s)\|^{2r}_{L^{2}([-N,N];H^{\zeta})}]\lesssim N.
\end{align*}
For $R_N > 0, N \in \mathbb{N}$, we define 

\begin{align*}
K_M := \cap_{N=M}^\infty\big{\{}
g_1; &\| g_1\|_{ C^\zeta([-N,N];W^{(\frac23-\fa-\kappa)\wedge\frac12,\frac65})} +\|g_1\|_{L^{2}([-N,N];H^{\zeta})} \leq R_N\big{\}}.
\end{align*}
 Using $(\frac23-\fa-\kappa)\wedge\frac12>\frac13$ we obtain the following compact embedding
\begin{align*}
      C^\zeta([-N,N];W^{(\frac23-\fa-\kappa)\wedge\frac12,\frac65})
     \cap L^{2}([-N,N];H^{\zeta})\subset C([-N,N];W^{\frac13,\frac65})\cap L^2([-N,N];L^2).
\end{align*}
Then it follows that $K_M$ is relatively compact
 in $\mathcal{T}$.  As a consequence, we deduce that the time shifts $S_tv, t\in\mathbb{R}$, are
tight on $\mathcal{T}$. Since  $S_tB$ is Wiener processes for every $t\in\mathbb{R}$, the law of $S_tB$ does not change
 and is tight.  Accordingly, for any $\epsilon > 0$ there is a compact set $\bar{K}_\epsilon$ in $\mathcal{T}$ such that
$$\sup_{t\in\mathbb{R}}\mathbf{P}(S_t(v,B) \in \bar{K}_\epsilon^c) < \epsilon.$$
This implies
$$\nu_{\tau,T} ( \bar{K}_\epsilon^c) = \frac1T\int_0^T\mathbf{P}(S_{t+\tau}(v,B)\in\bar{K}_\epsilon^c )\dif t < \epsilon.$$
Then by the similar argument as in \cite[Theorem 4.2]{HZZ22c}, for $\tau > 0$ the measures $\nu_{\tau,T-\tau}$ are tight on $\mathcal{T}$. Hence there is a weakly converging subsequence  $\nu_{\tau,T_n-\tau}\to\nu$ weakly in $\mathcal{P} (\mathcal{T} )$ for some $\nu \in\mathcal{P} (\mathcal{T} )$  and we know that $\nu$ is independent of the choice of $\tau$. Furthermore, we take a sequence
$\tau_m \to\infty$ and consider $\nu_{\tau_m,T_n-\tau_m} , m, n \in\mathbb{N}$. Denote by $\rm d$ the metric on $\mathcal{P} (\mathcal{T} )$ metrizing the weak
convergence. For $m\in\mathbb{N}$ we could find $n(m)$ such that ${\rm d}(\nu_{\tau_m ,T_{n(m)}-\tau_m},\nu) < \frac1m$. Hence, it follows
that $\nu_{\tau_m,T_{n(m)}-\tau_m} \to\tau$ weakly in $\mathcal{P} (\mathcal{T} )$ as $m \to\infty$.

Define for every $T\geq0$ the set 
 \begin{align*}
 &   A_T=\Big{\{} (v,B)\in\mathcal{T}; {\rm for} \ \forall\psi\in C^\infty(\mathbb{T}^2), \div \psi=0,-T\leq s\leq t<\infty,\\
&  \langle v(t), \psi\rangle + \int_s^t \langle
\div(v\otimes v+v\otimes z+z\otimes v+z^{:2:}), \psi\rangle \dif r 
= \langle v(s), \psi\rangle+  \int_s^t \langle
 z, \psi\rangle \dif r + \int_s^t  \langle \Delta v,\psi\rangle \dif r, \\
 &{\rm where\ }
 z(t)=\int_{-\infty}^te^{(t-s)(\Delta-I)}(-\Delta)^{\fa/2}\dif B_s, z^{:2:}\ {\rm is\ defined\ through}\ \eqref{z:2:}. \Big{\}}
\end{align*}
Since $(u,B)$ in the statement of the theorem satisfies  (\ref{1.1}), we have for all $t\geq0$  and $\tau\geq T$
$$\mathcal{L}[S_{t+\tau}(u-z, B)](A_T) = 1.$$
Hence, for $m$ large enough $\nu_{\tau_m,T_{n(m)}-\tau_m}(A_T) = 1$. By Jakubowski–Skorokhod representation theorem, there is a probability space $(\tilde{\Omega}, \tilde{\mathcal{F}},\tilde{\mathbf{P}})$ with a sequence of random variables $(\tilde{v}_m, \tilde{B}_m),m\in\mathbb{N}$, such that $\mathcal{L}[\tilde{v}_m, \tilde{B}_m]=\nu_{\tau_m,T_{n(m)}-\tau_m}$ and $(\tilde{v}_m+\tilde{z}_m, \tilde{B}_m)$ satisfy equation (\ref{1.1}) on $[ -T,\infty)$, where  $\tilde{z}_m=\int_{-\infty}^te^{(t-s)(\Delta-I)}(-\Delta)^{\fa/2}\dif \tilde{B}_m$.
 Moreover, there is a random process $(\tilde{v}, \tilde{B})$ having the law $\mathcal{L}[\tilde{v},\tilde{B}] = \nu$ so that  $(\tilde{v}_m, \tilde{B}_m)\to (\tilde{v}, \tilde{B})$
 $\tilde{\mathbf{P}}$-a.s. in $\mathcal{T}$. 
 
We denote $\tilde{u}:=\tilde{v}+\tilde{z}$,  where $\tilde{z}(t):=\int_{-\infty}^te^{(t-s)(\Delta-I)}(-\Delta)^{\fa/2}\dif \tilde{B}_s$. Then
 \eqref{6.1} and \eqref{6.2} follow from a lower-semicontinuity argument
and the related bound for the approximations.

 Then we  prove  that   $(\tilde{u},\tilde{B})$ is a solution to \eqref{1.1} on $\mathbb{R}$. In fact, we only need to prove the convergence of $\tilde{z}_n$ and $\tilde{z}_n^{:2:}$. By
 \begin{align}
  \tilde{z}_n(t) =&(-\Delta)^{\fa/2}\tilde{B}_n(t) + \int^t_{-\infty}(\Delta -I)e^{(t-s)(\Delta-I)}(-\Delta)^{\fa/2}\tilde{B}_n\dif s,   \label{zndecom}
 \end{align}
 we know that $\tilde{z}_n \to \tilde{z}$ in $C(\mathbb{R};H^{-4})\ \mathbf{\tilde{P}}$-a.s. 
 To prove the convergence of $\tilde{z}_{n}^{:2:}$,  let $\tilde{z}_{n,\epsilon},\tilde{z}_{\epsilon}$ be the corresponding mollification of $\tilde{z}_{n},\tilde{z}$ and $\tilde{z}_{n,\epsilon}^{:2:},\tilde{z}_{\epsilon}^{:2:}$ are defined by
 
 \begin{align}
     z^{:2:}_{n,\epsilon}:= z_{n,\epsilon} \otimes z_{n,\epsilon}-\tilde{\mathbf{E}}[z_{n,\epsilon} \otimes z_{n,\epsilon}], \ \
        z^{:2:}_{\epsilon}:= z_{\epsilon} \otimes z_{\epsilon}-\tilde{\mathbf{E}}[z_{\epsilon} \otimes z_{\epsilon}].\notag
     \end{align}
   As $\tilde{z}_n\overset{\rm d}{=}\tilde{z},\tilde{z}_{n,\epsilon}\overset{\rm d}{=}\tilde{z}_\epsilon$, for any $N>0$,  by Kolmogorov’s continuity criterion   we obtain
 $\tilde{z}_{n,\epsilon}^{:2:}\to \tilde{z}_{n}^{:2:}, \tilde{z}_{\epsilon}^{:2:}\to \tilde{z}^{:2:}$ in $L^2(\tilde{\Omega};C([-N
 ,N];C^{-2\fa-\kappa}))$ as $\epsilon\to0$  uniformly in $n$, i.e. for any  $\delta>0$, there exists $\epsilon_0>0$ such that for any $n\in\mathbb{N}$, 
 $$\tilde{\mathbf{E}}[\sup_{t\in[-N,N]}\|\tilde{z}_{n,\epsilon_0}^{:2:}(t)- \tilde{z}_{n}^{:2:}(t)\|_{C^{-2\fa-\kappa}}]<\delta/3,\ \ \tilde{\mathbf{E}}[\sup_{t\in[-N,N]}\|\tilde{z}_{\epsilon_0}^{:2:}(t)- \tilde{z}^{:2:}(t)\|_{C^{-2\fa-\kappa}}]<\delta/3.$$
 By \eqref{zndecom} and mollification estimates, there exists a $N_0\in\N$ such that for $n\geq N_0$, we have $$\tilde{\mathbf{E}}[\sup_{t\in[-N,N]}\|\tilde{z}_{n,\epsilon_0}^{:2:}(t)-\tilde{z}_{\epsilon_0}^{:2:}(t)\|_{C^{-2\fa-\kappa}}]<\delta/3,$$ which implies that for $n\geq N_0$.
 $$\tilde{\mathbf{E}}[\sup_{t\in[-N,N]}\|\tilde{z}_{n}^{:2:}(t)- \tilde{z}^{:2:}(t)\|_{C^{-2\fa-\kappa}}]<\delta.$$
 Then by passing  the limit in the equation we obtain \eqref{3.2*} holds on $[-T,\infty)$ for any $T>0$.

The shift invariance follows from the same argument as in \cite[Lemma 5.2]{BFH20}. Namely, it
holds for $G \in C_b(\mathcal{T})$ and $r \in\mathbb{R}$
\begin{align*}
    \int_{\mathcal{T}}G\circ S_r(v, B)\dif\nu(v, B) =   \int_{\mathcal{T}}G(v, B)\dif\nu(v, B).
    \end{align*}

To prove the non-uniqueness, by \eqref{3.7} and  \eqref{6.1}, we have
 \begin{align}
 \$ \tilde{u}  -v_0 \$ _{W^{-\fa-\kappa,\frac32},2}\lesssim  \$ \tilde{u}- \tilde{z}  -v_0 \$ _{W^{(\frac23-\fa-\kappa)\wedge\frac12,\frac65},2r}+ \$  \tilde{z} \$_{C^{-\fa-\kappa},2}\lesssim\sqrt{\aleph}+L\leq  C_0,\label{u-v_0}
\end{align}
where $C_0$ is independent of the choice of $v_0$. Define
\begin{align*}
    \mathcal{V}:=\{f\in C^1(\mathbb{
    T}^2;\mathbb{R}^2)|\div f=0,\langle f,1\rangle=\langle f,\varphi\rangle=0, {\rm where}\ \varphi(x_1,x_2)=(0,e^{2\pi ix_1})\}.
\end{align*}
Then  there exists a sequence $f_k\in \mathcal{V},k\in\mathbb{N}$ such that $\|f_i-f_j\|_{W^{-\fa-\kappa,\frac32}}> 2C_0$ for any $i,j\in\mathbb{N}$. For any $k\in\mathbb{N}$, we start iteration from  $v_0=f_k$ and obtain corresponding stationary solution denoted by  $u_k$. Then by \eqref{u-v_0}, $u_k$ is located in the ball $B(f_k,C_0)$ in the space $L^2(\Omega;C(\mathbb{R};W^{-\fa-\kappa,\frac32}))$. By the choice of $f_k$,   $u_k$ are different stationary solutions. 

Moreover, the distributions are not Gaussian $N(0,\frac12(-\Delta)^{\fa-1})$.  Otherwise, we take inner product with $\varphi(x_1,x_2)=(0,e^{2\pi ix_1})$. By Young's inequality  we have
\begin{align}
    \mathbf{E}\langle u_k(0),\varphi\rangle^2&=\mathbf{E}[\langle v_k(0),\varphi\rangle+\langle z(0),\varphi\rangle]^2\leq 41\mathbf{E}\langle v_k(0),\varphi\rangle^2+\frac{41}{40}\mathbf{E}\langle z(0),\varphi\rangle^2,\notag
\end{align}
where we denote $v_k=u_k-z$.
By H\"older's inequality, Sobolev embedding, \eqref{6.1} and the choice of $f_k$  we have 
\begin{align}
\mathbf{E}\langle v_k(0),\varphi\rangle^2=\mathbf{E}\langle v_k(0)-f_k,\varphi\rangle^2\leq C_1\|\varphi\|_{L^\infty}^2 \$ v_k-f_k \$ _{W^{(\frac23-\fa-\kappa)\wedge\frac12,\frac65},2r}^2\leq C_1 \aleph,\label{vkphi}
\end{align}
where $C_1$ is the implicit constant. By direct calculation we have
\begin{align}
\mathbf{E}\langle z(0), \varphi\rangle^2=\int_{-\infty}^0e^{2s(4\pi^2+1)}{(2\pi)^{2\fa}}\dif s=\frac{(2\pi)^{2\fa}}{2(4\pi^2+1)}, \notag
\end{align}
which implies that
\begin{align*}
\frac  {1 }{2(2\pi)^{2-2\fa}}=\frac12\langle (-\Delta)^{\fa-1}\varphi,\varphi\rangle= \mathbf{E}\langle u_k(0),\varphi\rangle^2\leq  \frac{41 (2\pi)^{2\fa}}{80(4\pi^2+1)}+41C_1\aleph,
\end{align*}
which is a contradiction by the choice of $\aleph$ in \eqref{defaleph}. We finish the proof of Theorem \ref{thm:6.1}.
\end{proof}

The first claim in Theorem \ref{thm:exist} follows. By a general  Krein–Milman argument as in \cite{BFH20, FFH21, HZZ22,HZZ22c,HLZZ23} and using Theorem \ref{thm:6.1} we obtain the second claim in Theorem \ref{thm:exist}.

\bt\label{thm:erg}
 Assume \eqref{defaleph} and let $\fa\in[0,\frac{1}{3}), r\geq1$.  For any divergence-free, mean-zero condition $v_0\in C^1(\mathbb{T}^2;\R^2)$, there exist $C>0$, some $\kk>0$ small enough and an ergodic stationary solution $((\Omega, \mathcal{F},(\mathcal{F}_t)_{t \in\mathbb{R}}, \mathbf{P}), u, B)$ satisfying
\begin{align}
 \$ u- z-v_0 \$ _{W^{(\frac23-\fa-\kappa)\wedge\frac12,\frac65},2r}\leq\sqrt{\aleph}\label{6.1*}
\end{align}
 with  $z(t)=\int_{-\infty}^te^{(t-s)(\Delta-I)}(-\Delta)^{\fa/2}\dif B_s$. For some $\zeta\in(0,1)$  and every $N\in\mathbb{N}$
\begin{align}
\mathbf{E}[
 \|u- z\|_{L^{2}([-N,N];H^{\zeta})}^{2r}+\|u- z\|_{C^\zeta ([-N,N],W^{(\frac23-\fa-\kappa)\wedge\frac12,\frac65})}^{2r}]\leq CN. \label{6.2*}
\end{align}
There are infinitely many ergodic stationary solutions by choosing different $v_0$ and the laws are not Gaussian $N(0,\frac12(-\Delta)^{\fa-1})$.
\et
\begin{proof}
We obverse that the
set of all laws of stationary solutions satisfying \eqref{6.1*}, \eqref{6.2*} is non-empty, convex, tight and closed which follows from the same argument as the proof of Theorem \ref{thm:6.1}. Hence there exists an extremal point, which is an ergodic stationary solution by a classical argument.

  By a similar argument as Theorem \ref{thm:6.1},  we start iteration from $v_0 = f_k$ and \eqref{6.1*} implies these ergodic stationary solutions are different and the laws are not Gaussian $N(0,\frac12(-\Delta)^{\fa-1})$. 
\end{proof}

\appendix
 \renewcommand{\appendixname}{Appendix~\Alph{section}}
  \renewcommand{\theequation}{A.\arabic{equation}}

  \section{Estimates of $\rho$ and $a_{(\xi)}$}
  \label{s:appA}
For completeness, we include here the detailed proof of the estimates (\ref{4.6}), (\ref{4.8}) employed in Section \ref{covq+12}.
We first give the estimates of $\rho$ for any $t\in\mathbb{R}$. In view of the definition of $\rho$  we obtain 
\begin{align}
    \|\rho\|_{L^\infty}\lesssim  l+\|\mathring{R}_l\|_{L^\infty}.\notag
\end{align}
Furthermore, by mollification estimates, the embedding $W^{4,1}\subset L^\infty$ (in fact the time-space dimension is three), for $N\geq 0$
$$\|\mathring{R}_l\|_{C_{[t,t+1],x}^N}\lesssim l^{-4-N}\|\mathring{R}_q\|_{L^1_{[t-l,t+1]}L^1},$$
which in particular leads to
$$\|\rho\|_{C_{{[t,t+1]},x}^0}\lesssim l+l^{-4}\|\mathring{R}_q\|_{L^1_{[t-l,t+1]}L^1}.$$
Next we estimate the $C_{t,x}^N$-norm for $N\geq1$. We apply the chain rule in \cite[Proposition C.1]{BDLIS15} to $f(z)=\sqrt{l^2+z^2},|D^mf(z)|\lesssim l^{-m+1}$ to obtain
\begin{align}
\|\rho\|_{C_{{[t,t+1]},x}^N}&\lesssim\|\sqrt{l^2+|\mathring{R}_l|^2}\|_{C_{{[t,t+1]},x}^0}+\|Df\|_{C^0}\|\mathring{R}_l\|_{C_{{[t,t+1]},x}^N}+\|Df\|_{C^{N-1}}\|\mathring{R}_l\|_{C_{{[t,t+1]},x}^1}^N\notag\\
&\lesssim l^{-N-4}\|\mathring{R}_q\|_{L^1_{[t-l,t+1]}L^1}+l^{-6N+1}\|\mathring{R}_q\|_{L^1_{[t-l,t+1]}L^1}^N +1,\label{A.1s}
\end{align}
which implies (\ref{4.6}).

Next we estimate $a_{(\xi)}$
in $C_{[t,t+1],x}^N$-norm. By Leibniz rule we get
$$\|a_{(\xi)}\|_{C_{{[t,t+1]},x}^N}\lesssim\sum_{m=0}^N\|\rho^{1/2}\|_{C_{{[t,t+1]},x}^m}\|\gamma_\xi(\Id-\frac{\mathring{R}_l}{\rho})\|_{C_{{[t,t+1]},x}^{N-m}}.$$
Applying \cite[Proposition C.1]{BDLIS15} to $f(z)=z^{1/2}$, $|D^mf(z)|\lesssim|z|^{1/2-m}$, for $m=1,...,N$, and using (\ref{A.1s}) 
we obtain  for $m\geq1$
\begin{align*}
\|\rho^{1/2}\|_{C_{{[t,t+1]},x}^m}&\lesssim \|\rho^{1/2}\|_{C_{{[t,t+1]},x}^0}+l^{-1/2}\|\rho\|_{C_{{[t,t+1]},x}^m}+l^{1/2-m}\|\rho\|_{C_{{[t,t+1]},x}^1}^m\\
&\lesssim l^{-2}\|\mathring{R}_q\|^{1/2}_{L^1_{[t-l,t+1]}L^1}+l^{-9/2-m}\|\mathring{R}_q\|_{L^1_{[t-l,t+1]}L^1}+l^{-6m+1/2}(\|\mathring{R}_q\|_{L^1_{[t-l,t+1]}L^1}^m+1)\\
&\lesssim l^{-6m+1/2}(1+\|\mathring{R}_q\|_{L^1_{[t-l,t+1]}L^1}^m).
\end{align*}

Next we estimate $\gamma_\xi(\Id-\frac{\mathring{R}_l}{\rho})$, by \cite[Proposition C.1]{BDLIS15} we need to estimate
$$\|\frac{\mathring{R}_l}{\rho}\|_{C_{{[t,t+1]},x}^{N-m}}+\|\frac{\nabla_{t,x}\mathring{R}_l}{\rho}\|_{C_{{[t,t+1]},x}^0}^{N-m}+\|\frac{\mathring{R}_l}{\rho^2}\|_{C_{{[t,t+1]},x}^0}^{N-m}\|\rho\|_{C_{{[t,t+1]},x}^1}^{N-m}.$$
We use $\rho\geq l$ to have 
$$\|\frac{\nabla_{t,x}\mathring{R}_l}{\rho}\|_{C_{{[t,t+1]},x}^0}^{N-m}\lesssim l^{-(N-m)}l^{-5(N-m)}\|\mathring{R}_q\|^{N-m}_{L^1_{[t-l,t+1]}L^1}\lesssim l^{-6(N-m)}\|\mathring{R}_q\|_{L^1_{[t-l,t+1]}L^1}^{N-m},$$
and in view of $|\frac{\mathring{R}_l}{\rho}|\leq 1$
$$\|\frac{\mathring{R}_l}{\rho^2}\|_{C_{{[t,t+1]},x}^0}^{N-m}\lesssim\|\frac{1}{\rho}\|_{C_{{[t,t+1]},x}^0}^{N-m}\lesssim l^{-(N-m)},$$
and by (\ref{A.1s}) 
$$\|\rho\|_{C_{{[t,t+1]},x}^1}^{N-m}\lesssim  l^{-5(N-m)}\|\mathring{R}_q\|^{N-m}_{L^1_{[t-l,t+1]}L^1}+1.$$
Moreover, to estimate the first term we write
\begin{align*}
\|\frac{\mathring{R}_l}{\rho}\|_{C_{{[t,t+1]},x}^{N-m}}&\lesssim  \sum_{k=0}^{N-m}\|\mathring{R}_l\|_{C_{{[t,t+1]},x}^k}\|\frac{1}{\rho}\|_{C_{{[t,t+1]},x}^{N-m-k}},
\end{align*}
where for $N - m - k = 0$ we have
$$\|\frac{1}{\rho}\|_{C_{[t,t+1],x}^0}\lesssim l^{-1},$$
and for $k=0,..,N-m-1$
using (\ref{A.1s})
\begin{align*}
\|\frac{1}{\rho}\|&_{C_{{[t,t+1]},x}^{N-m-k}}\lesssim \|\frac{1}{\rho}\|_{C_{{[t,t+1]},x}^{0}}+l^{-2}\|{\rho}\|_{C_{{[t,t+1]},x}^{N-m-k}}+l^{-(N-m-k)-1}\|{\rho}\|_{C_{{[t,t+1]},x}^1}^{N-m-k}\\
&\lesssim  l^{-(N-m-k)-1}+l^{-(N-m-k)-6}\|\mathring{R}_q\|_{L^1_{[t-l,t+1]}L^1}+l^{-6(N-m-k)-1}\|\mathring{R}_q\|_{L^1_{[t-l,t+1]}L^1}^{N-m-k}.
\end{align*}
Thus we obtain
\begin{align*}
\|\frac{\mathring{R}_l}{\rho}\|_{C_{{[t,t+1]},x}^{N-m}}&\lesssim   l^{-10-(N-m)}\|\mathring{R}_q\|_{L^1_{[t-l,t+1]}L^1}^2+l^{-5-6(N-m)}(\|\mathring{R}_q\|^{N-m}_{L^1_{[t-l,t+1]}L^1}+1).
\end{align*}
Finally the above bounds leads for $N-m\geq1$
\begin{align}
 \|\gamma_\xi(\Id-\frac{\mathring{R}_l}{\rho})\|_{C_{{[t,t+1]},x}^{N-m}}\lesssim l^{-5-6(N-m)}(\|\mathring{R}_q\|^{N-m}_{L^1_{[t-l,t+1]}L^1}+\|\mathring{R}_q\|^{2}_{L^1_{[t-l,t+1]}L^1} +1). \label{gammaxi}
\end{align}
Combining this with the bounds for $\rho^{1/2}$ above yields for $N\geq2$
\begin{align}
\|a_{(\xi)}\|_{C_{{[t,t+1]},x}^N}&\lesssim l^{-7-6N}(\|\mathring{R}_q\|_{L^1_{[t-l,t+1]}L^1}+1)^{N+1}.\label{A.8s}
\end{align}
It is easy to see that
\begin{align*}
\|a_{(\xi)}\|_{C_{{[t,t+1]},x}^0}&\lesssim l^{-2}(\|\mathring{R}_q\|_{L^1_{[t-l,t+1]}L^1}+1)^{1/2}.
\end{align*}
By interpolation, (\ref{A.8s}) holds for $N =0, 1$.

 \renewcommand{\appendixname}{Appendix~\Alph{section}}
  \renewcommand{\theequation}{B.\arabic{equation}}

 \section{Introduction to accelerating jet flows}
  \label{s:appB}

In this part we recall the construction of accelerating jets from \cite[Section 3.1-3.3]{CL20}. We point
out that the construction is entirely deterministic, that is, none of the functions below depends on
$\omega$. Let us begin with the following geometric lemma. Recall that $\mathcal{S}^{2\times 2}$ is the set of symmetric $2\times 2$ matrics.

\bl$($\cite[Lemma 3.1]{CL20}$)$\label{lem:B.1}
Denote by $\bar{B}_{1/2}(\mathrm{Id})$ the closed ball of radius $1/2$ around the identity matrix Id, in the
space of $\mathcal{S}^{2\times 2}$. There exists $\Lambda\in \mathbb{S}^{1}\cap \mathbb{Q}^2$ such that for each $\xi\in\Lambda$ there exists
a $C^\infty$-function $\gamma_\xi$: $\bar{B}_{1/2}(\mathrm{Id})\to \mathbb{R}$ such that
$$R=\sum_{\xi\in\Lambda}\gamma_\xi^2(R)(\xi\otimes\xi)$$
for every symmetric matrix satisfying $|R-\mathrm{Id}|\leq 1/2$.
\el

We choose a collection of distinct points $p_{(\xi)} \in \mathbb{T}^2$
for $\xi\in\Lambda$ and a number $\mu_0 >2$ such
that
$$\cup_{\xi\in\Lambda}B_{\frac{2}{\mu_0}}(p_{(\xi)})\subset[0,1]^2,$$
and for $\xi\neq\xi'\in\Lambda$
$$\dif(p_{(\xi)},p_{(\xi')})>\frac{2}{\mu_0},
$$
the points $p_{(\xi)}$ will be the centers of $W_{(\xi)}$ and $W_{(\xi)}^{(c)}$.\\
For $\xi\in\Lambda$, let us introduce unit vectors $\xi^{\perp}:=(\xi_2,-\xi_1),$
and their associated coordinates
$$x_{(\xi)}=(x-p_{(\xi)})\cdot\xi,\ \ y_{(\xi)}=(x-p_{(\xi)})\cdot\xi^{\perp}\ \ {\rm{for}}\ \ x\in \mathbb{T}^2 .$$
We first construct a family of stationary jets. To this end, we introduce the parameters:
$\nu,\mu\in\mathbb{N}$ with $\mu_0<\nu\leq \mu.$ Now we choose compactly supported nontrivial $\varphi, \psi\in C_c^\infty((-\frac{1}{\mu_0} , \frac{1}{\mu_0}
))$ and define non-periodic potentials $\tilde{\Psi}_{(\xi)}\in C_c^\infty(\mathbb{R}^2)$ and vector fields $\tilde{W}_{(\xi)},\tilde{W}_{(\xi)}^{(c)}\in C_c^{\infty}(\mathbb{R}^2)$
$$\tilde{\Psi}_{(\xi)}(x)=c_{(\xi)}\mu^{-1}(\nu\mu)^{1/2}\varphi(\nu x_{(\xi)})\psi(\mu y_{(\xi)}),$$
$$\tilde{W}_{(\xi)}(x)=-c_{(\xi)}(\nu\mu)^{1/2}\varphi(\nu x_{(\xi)})\psi'(\mu y_{(\xi)})\xi,$$
$$\tilde{W}^{(c)}_{(\xi)}(x)=c_{(\xi)}\nu\mu^{-1}(\nu\mu)^{1/2}\varphi'(\nu x_{(\xi)})\psi(\mu y_{(\xi)})\xi^{\perp},$$
where $c_{(\xi)}>0$ are normalizing constants such that
\begin{align}
\int_{\mathbb{T}^2}\!\!\!\!\!\!\!\!\!\!\!\!\; {}-{}  \tilde{W}_{(\xi)}\otimes \tilde{W}_{(\xi)}\dif x=\xi\otimes\xi.\notag
\end{align}
We periodize so that $\tilde{\Psi}_{(\xi)},\tilde{W}_{(\xi)}$ and $\tilde{W}^{(c)}_{(\xi)}$ are viewed as periodic functions on $\mathbb{T}^2$.
Finally note that each $\tilde{W}_{(\xi)}$ has disjoint support
$$\supp \tilde{W}_{(\xi)} \cap \supp \tilde{W}_{(\xi')} = \emptyset\ \ if\ \ \xi\neq\xi',$$
and for any $\xi\in\Lambda$, the following identities hold
\begin{align}
\div(\tilde{W}_{(\xi)}\otimes \tilde{W}_{(\xi)})=\xi\cdot\nabla|\tilde{W}_{(\xi)}|^2\xi,\notag
\end{align}
$$\nabla^{\perp}\tilde{\Psi}_{(\xi)}=-\partial_{y_{(\xi)}}\tilde{\Psi}_{(\xi)}\xi+\partial_{x_{(\xi)}}\tilde{\Psi}_{(\xi)}\xi^{\perp}=\tilde{W}_{(\xi)}+\tilde{W}^{(c)}_{(\xi)}.$$

Next, we introduce a simple method to avoid the collision of the
support sets of different $\tilde{W}_{(\xi)}$. Let us first choose temporal functions $g_{(\xi)}$ and $h_{(\xi)}$ to oscillate the
building blocks $\tilde{W}_{(\xi)}$ intermittently in time. Let $G \in C_c^\infty(0, 1)$ be such that
$$\int_0^1G^2(t)\dif t=1,\ \ \int_0^1G(t)\dif t=0.$$
For any $\eta \geq 1$, we define $\tilde{g}_{(\xi)}: \mathbb{T}\to \mathbb{R}$ as the 1-periodic extension of $\eta^{1/2}G(\eta(t -t_\xi))$, where $t_\xi$ are chosen so that $\tilde{g}_{(\xi)}$ have disjoint supports for different $\xi$. In other words,
$$\tilde{g}_{(\xi)}(t)=\sum_{n\in\mathbb{Z}}\eta^{1/2}G(\eta(n+t -t_\xi)).$$
We will also oscillate the velocity perturbation at a large
frequency $\sigma\in\mathbb{N}$. So we define
$$g_{(\xi)}(t)=\tilde{g}_{(\xi)}(\sigma t).$$
Then we have
\begin{align}\label{gwnp}
    \|g_{(\xi)}\|_{W^{n,p}([0,1])}\lesssim(\sigma\eta)^n\eta^{1/2-1/p}.
\end{align}
For the corrector term we define $h_{(\xi)}:\mathbb{T}\to\mathbb{R}$ by
\begin{align}
h_{(\xi)}(t)=\int_0^{\sigma t}(\tilde{g}^2_{(\xi)}(s)-1)\dif s.\notag
\end{align}

In view of the zero-mean condition for ${\tilde{g}}^2_{(\xi)}(t)-1$, these $h_{(\xi)}$ are $\mathbb{T}/\sigma$-periodic  and we have
\begin{align}
    \|h_{(\xi)}\|_{L^\infty}\leq1.\label{hlinfty}
\end{align}
Moreover, we have the identity
\begin{align}
\partial_t(\sigma^{-1}h_{(\xi)})=g^2_{(\xi)}-1.    \label{parth}
\end{align}

Now we will let the stationary flows $\tilde{W}_{(\xi)}$ travel along $\mathbb{T}^2$ in time, relating the velocity of the
moving support sets to the intermittent oscillator $g_{(\xi)}$. More precisely, we define
\begin{align}\label{Wxi}
W_{(\xi)}(x,t):&=\tilde{W}_{(\xi)}(\sigma x+\phi_{(\xi)}(t)\xi),\notag\\
W_{(\xi)}^{(c)}(x,t):&=\tilde{W}_{(\xi)}^{(c)}(\sigma x+\phi_{(\xi)}(t)\xi),\\
\Psi_{(\xi)}(x,t):&=\tilde{\Psi}_{(\xi)}(\sigma x+\phi_{(\xi)}(t)\xi),\notag
\end{align}
where $\phi_{(\xi)}(t)$ is defined by
$$\phi'_{(\xi)}=\theta g_{(\xi)}(t)$$
for some $\theta\in\mathbb{N}$. Hence we have
\begin{align}
\int_{\mathbb{T}^2}\!\!\!\!\!\!\!\!\!\!\!\!\; {}-{}  {W}_{(\xi)}\otimes{W}_{(\xi)}\dif x=\xi\otimes\xi,\label{wpingjun}
\end{align}
and each $g_{(\xi)}{W}_{(\xi)}$ has disjoint support
\begin{align}
    \supp (g_{(\xi)}{W}_{(\xi)}) \cap \supp (g_{(\xi')}{W}_{(\xi')}) = \emptyset\ \ if\ \ \xi\neq\xi'.\label{wbujiao}
\end{align}

Then, we claim that for $N\geq0,p\in [1, \infty]$ the
following holds
\begin{align}
\|\nabla^NW_{(\xi)}\|_{C_tL^p}+\nu^{-1}\mu\|\nabla^NW^{(c)}_{(\xi)}\|_{C_tL^p}+\mu\|\nabla^N\Psi_{(\xi)}\|_{C_tL^p}\lesssim(\sigma\mu)^N (\nu\mu)^{1/2-1/p},\label{2.4}
\end{align}
\begin{align}
\|(\xi\cdot\nabla)\Psi_{(\xi)}\|_{C_tL^p}\lesssim\mu^{-1}\sigma\nu (\nu\mu)^{1/2-1/p},\label{2.4*}
\end{align}
where the implicit constants may depend on $N,p$, but are independent of $\nu,\sigma,\theta,\eta $ and $\mu$. These estimates can be easily deduced from the definitions.

Moreover we have
\begin{align}
\sigma^{-1}\nabla^{\perp}\Psi_{(\xi)}&=W_{(\xi)}+W^{(c)}_{(\xi)},\label{2.6}\\
\partial_t|W_{(\xi)}|^2\xi&=\sigma^{-1}\theta g_{(\xi)}\div(W_{(\xi)}\otimes W_{(\xi)}),\label{B.9}\\
\partial_t\Psi_{(\xi)}&=\sigma^{-1}\theta g_{(\xi)}(\xi\cdot\nabla)\Psi_{(\xi)}.\label{B.10}
\end{align}

 \renewcommand{\appendixname}{Appendix~\Alph{section}}
  \renewcommand{\theequation}{C.\arabic{equation}}
  \section{Some technical tools}
  \label{s:appC}
\subsection{Improved H\"older’s inequality on $\mathbb{T}^d$}\label{ihiot2}
This theorem allows us to quantify the decorrelation in the usual H\"older’s
inequality when we increase the oscillation of one function. The proof follows from the  same argument as \cite[Lemma 2.1]{MS18} and \cite[Theorem C.1]{LZ23}.
\bt\label{ihiot}
Let $p \in [1, \infty],u\in\mathbb{R}$ and $a\in C^1(\mathbb{R}^d; \mathbb{R}), f\in L^p(\mathbb{T}^d ; \mathbb{R})$. Then for any $\sigma\in\mathbb{N}$,
$$| \|af(\sigma\cdot)\|_{L^p([u,u+1]^d)}- \|a\|_{L^p([u,u+1]^d)}\|f\|_{L^p(\mathbb{T}^d)}|\lesssim \sigma^{-1/p}\|a\|_{C^{0,1}([u,u+1]^d)}\|f\|_{L^p(\mathbb{T}^d)}.$$
Here $C^{0,1}$ is the Lipschitz norm.
In particular, if $d=1$ and $f$ is mean-zero we have
$$|\int_u^{u+1}a(t)f(\sigma t)\dif t|\lesssim\frac{1}{\sigma}\|a\|_{C^{0,1}_{[u,u+1]}}\|f\|_{L^1
(\mathbb{T})}.$$
\et

\subsection{Estimates for the heat operator}\label{eftho}
In this section we give the basic estimates for the heat semigroup $P_t:=e^{t\Delta}$. Let $T\geq1$. 

\bl$($\cite[Lemma 2.8]{ZZZ20}$)$\label{C.2}
For any $\theta > 0$, $\alpha \in \mathbb{R}, p, q \in [1, \infty]$ and all $t \in(0,T]$,
$$ \|P_tf\|_{ B^{\theta+\alpha}_{p,q}} \lesssim T^{\theta/2}t^{-\theta/2}\| f\|_{ B^\alpha_{p,q}}.$$
Here the implicit constant is independent of $T$.

 For any $\beta\in(0,2), t\in[0,1]$, 
$$ \|P_tf-f\|_{L^p}\lesssim t^{\frac\beta2}\|f\|_{B^{\beta}_{p,\infty}}.$$
\el

Let $\mathcal{I}_cf = \int_{0}^\cdot e^{-c(\cdot -s)} P_{\cdot -s}f\dif s$ 
for some $c>0$. Then we have
\bl$($\cite[Lemma 2.9]{ZZZ20}$)$\label{C.3}
Let $\alpha \in(0, 2),\beta\in[0,2], p\in[1,\infty]$.
Then 
$$\| \mathcal{I}_cf\|_{C_{[0,2]}B^\alpha_{p,\infty}}+\|\mathcal{I}_cf\|_{C_{[0,2]}^{\frac\alpha2}L^p} \lesssim c^{\beta/2-1} \| f\|_{C_{[0,2]}B^{\alpha-\beta}_{p,\infty}}.$$
\el

\end{document}